\pdfoutput=1
\RequirePackage{ifpdf}
\ifpdf 
\documentclass[pdftex]{sigma}
\else
\documentclass{sigma}
\fi

\usepackage[all]{xy}

\newcommand\tstrut{\rule{0pt}{2.6ex}}

\newcommand{\ext}{\raise1pt\hbox{\small$\bigwedge$}}
\newcommand{\tr}{\mbox{\bf tr}}
\newcommand{\dirac}{/\kern-5pt\partial}

\newcommand{\lra}{\longrightarrow}

\numberwithin{equation}{section}

\newtheorem{Theorem}{Theorem}[section]
\newtheorem{Corollary}[Theorem]{Corollary}
\newtheorem{Lemma}[Theorem]{Lemma}
\newtheorem{Proposition}[Theorem]{Proposition}
 { \theoremstyle{definition}
\newtheorem{Definition}[Theorem]{Definition}
\newtheorem{Remark}[Theorem]{Remark}
\newtheorem{Example}[Theorem]{Example}
}

\begin{document}

\newcommand{\arXivNumber}{1506.07681}

\renewcommand{\PaperNumber}{072}

\FirstPageHeading

\ShortArticleName{Spinorially Twisted Spin Structures, II}

\ArticleName{Spinorially Twisted Spin Structures.~II:\\ Twisted Pure Spinors, Special Riemannian Holonomy\\ and Clifford Monopoles}

\Author{Rafael HERRERA~$^\dag$ and Noemi SANTANA~$^\ddag$}

\AuthorNameForHeading{R.~Herrera and N.~Santana}

\Address{$^\dag$~Centro de Investigaci\'on en Matem\'aticas, A.P.~402, Guanajuato, Gto., C.P.~36000, M\'exico}
\EmailD{\href{mailto:rherrera@cimat.mx}{rherrera@cimat.mx}}

\Address{$^\ddag$~Universidad Marista Valladolid, Jos\'e Juan Tablada 1111, Santa Maria de Guido,\\
\hphantom{$^\ddag$}~58090 Morelia, Mich., M\'exico}
\EmailD{\href{mailto:noemi_santana@hotmail.com}{noemi{\textunderscore}santana@hotmail.com}}

\ArticleDates{Received February 12, 2019, in final form September 07, 2019; Published online September 22, 2019}

\Abstract{We introduce a notion of twisted pure spinor in order to characterize, in a unified way, all the special Riemannian holonomy groups just as a classical pure spinor characterizes the special K\"ahler holonomy. Motivated by certain curvature identities satisfied by manifolds admitting parallel twisted pure spinors, we also introduce the Clifford monopole equations as a natural geometric generalization of the Seiberg--Witten equations. We show that they restrict to the Seiberg--Witten equations in 4 dimensions, and that they admit non-trivial solutions on manifolds with special Riemannian holonomy.}


\Classification{53C10; 53C25; 53C27; 58J60}

\section{Introduction}

The Berger--Simons' theorem \cite{Berger, Simons} states that the holonomy group of an irreducible non-locally symmetric oriented Riemannian manifold is contained in one of the groups in Table~\ref{table: Berger list} (see~\cite{Joyce, Salamon} for extensive accounts on the theory of Riemannian holonomy).
\begin{table}[h!]\centering
\begin{tabular}{|c|l|}\hline
group & geometry\\\hline
${\rm SO}(m)$ & generic\\
${\rm U}(m)$ & K\"ahler\\
${\rm SU}(m)$ & Calabi--Yau\\
${\rm Sp}(m){\rm Sp}(1)$ & quaternion-K\"ahler\\
${\rm Sp}(m)$ & hyper-K\"ahler\\
${\rm Spin}(7)$ & exceptional\\
$G_2$ & exceptional\\\hline
\end{tabular}
\caption{Special Riemannian holonomy groups and geometries.}\label{table: Berger list}
\end{table}

Manifolds whose holonomies are contained in
\begin{gather*} {\rm SU}(m), \quad {\rm Sp}(m), \quad {\rm Spin}(7),\quad G_2
\end{gather*}
are known to be Ricci-flat, Spin and to carry parallel spinors for their
classical (untwisted) Spin structures~\cite{Hitchin,Wang}. In the cases of ${\rm SU}(m)$, ${\rm Spin}(7)$ and $G_2$, such classical parallel
spinors characterize the holonomies~\cite{Lawson}.
There exist, however, manifolds with holonomies contained in~${\rm U}(n)$ and ${\rm Sp}(n){\rm Sp}(1)$ which are not Spin, but do admit
natural ${\rm Spin}^c$ and ${\rm Spin}^q$ structures respectively.

The purpose of this note is to introduce a suitable notion of twisted pure spinor for manifolds admitting
spinorially twisted Spin structures \cite{Espinosa-Herrera} in order to give a unified
treatment of special Riemannnian holonomies.
While there have been other unifying efforts involving the normed division algebras
(e.g.,~\cite{Leung}),
our approach is centered on Clifford algebras and twisted spinors.

The main motivation for our work has been the relationship between (classical) pure spinors and complex structures. More precisely,
let $\mathbb{R}^n$ and $\Delta_n$ denote the (real and complex) standard representations of ${\rm SO}(n)$ and ${\rm Spin}(n)$ respectively.
\'E. Cartan defined (classical) pure spinors in terms of maximal isotropic subspaces~\cite{Cartan}.
Equivalently, a spinor $\phi\in \Delta_n$ is pure if for every vector $X\in\mathbb{R}^n$ there exists a vector
$Y\in\mathbb{R}^n$ satisfying the following equation~\cite{Lawson}
\begin{gather}
X\cdot\phi = {\rm i}Y\cdot \phi,\label{eq: first equation}
\end{gather}
where ``$\cdot$'' denotes Clifford multiplication.
This condition says that the two subspaces $\mathbb{R}^n\cdot\phi$ and ${\rm i}\mathbb{R}^n\cdot\phi$ of
$\Delta_n$ coincide, which
allows the transfer of the effect of multiplication by ${\rm i}=\sqrt{-1}$ within the complex space
$\Delta_n$ to the real vector space $\mathbb{R}^n$.
Indeed, setting $Y=J(X)$ in~\eqref{eq: first equation} and a little algebraic manipulation show that $J$
defines a complex structure on $\mathbb{R}^n$. The isotropy group of a pure spinor is isomorphic to ${\rm SU}(n/2)$ and,
therefore, a spin manifold admits a parallel pure spinor field if and only if it is special K\"ahler~\cite{Lawson}.
Furthermore, A.~Moroianu proved that a~${\rm Spin}^c$ manifold admits a (similarly defined) parallel pure spinor field if
and only if it is K\"ahler~\cite{Moroianu}, a~result which includes non-Spin non-Ricci-flat K\"ahler manifolds.

Now recall that the tangent spaces of quaternionic K\"ahler manifolds and 8-manifolds with ${\rm Spin}(7)$ holonomy
are representation spaces of $\mathfrak{sp}(1)\cong \mathfrak{spin}(3)$ and $\mathfrak{spin}(7)$ respectively,
which can be viewed as restrictions of representations of the even Clifford algebras ${\rm Cl}_3^0$ and ${\rm Cl}_7^0$ to such subalgebras.
Thus, we conjectured that the special Riemannian holonomies must be determined by twisted spinors which,
somehow, induce a {\em transfer} of algebraic
structure from an even Clifford algebra to the bundle of endomorphisms of the tangent spaces of the manifold (see \cite{Santana}).
More precisely,
let $M$ be a smooth Riemannian manifold and $F$ be an auxiliary Riemannian vector bundle of rank~$r$.
Let $(e_1,\dots ,e_n)$ and $(f_1,\dots ,f_r)$ be local orthonormal frames of $TM$ and $F$ respectively,
$S(TM)$ and $S(F)$ be the locally defined
spinor vector bundles of $M$ and $F$, and suppose $m\in\mathbb{N}$ is such that the bundle
$S(TM)\otimes S(F)^{\otimes m}$ is globally defined.
A spinor field $\phi\in\Gamma\big(S(TM)\otimes S(F)^{\otimes m}\big)$ determines maps
\begin{gather*}
 T_xM \lra T_xM\cdot\phi_x \subset S(T_xM)\otimes S(F_x)^{\otimes m} ,\\
 T_xM \lra T_xM\cdot \kappa_{r*}^m(f_kf_l)\cdot \phi_x \subset S(T_xM)\otimes S(F_x)^{\otimes m} ,
\end{gather*}
at $x\in M$, for all $1\leq k<l\leq r$, where $\kappa_{r*}^m$ is the induced representation of~$\mathfrak{spin}(r)$ on~$S(F)^{\otimes m}$.
Given a pair $k<l$, we have a projection map
\begin{gather*}
 T_xM\cdot \kappa_{r*}^m(f_kf_l)\cdot \phi_x \lra T_xM\cdot \phi_x,\\
 X\cdot \kappa_{r*}^m(f_kf_l)\cdot \phi_x ,\longmapsto, \sum_{j=1}^n
\operatorname{Re}\big\langle X\cdot \kappa_{r*}^m(f_kf_l)\cdot \phi_x,e_j\cdot \phi_x\big\rangle e_j\cdot\phi_x,
\end{gather*}
which, in turn, gives the map
\begin{gather*}
 T_xM \lra T_xM,\\
 X \longmapsto \sum_{j=1}^n \operatorname{Re}\big\langle X\cdot \kappa_{r*}^m(f_kf_l)\cdot \phi_x,e_j\cdot \phi_x\big\rangle e_j.
\end{gather*}
The conjectured transfer of algebraic structure from the even Clifford algebra ${\rm Cl}_r^0$ to $\operatorname{End}(TM)$ must then be encoded in these maps.
Thus, we will define twisted pure spinors in such a way that the local $2$-forms and endomorphisms
\begin{gather*}
 \eta_{kl}^\phi(X,Y) = \operatorname{Re}\langle X\wedge Y\cdot \kappa_{r*}^m(f_kf_l)\cdot \phi,\phi\rangle ,\\
 \hat\eta_{kl}^\phi(X) = (X\lrcorner\eta_{kl}^\phi)^\sharp,
\end{gather*}
$X,Y\in T_xM$, induce a non-trivial representation of ${\rm Cl}_r^0$ on $T_xM$, where $\sharp$ denotes metric dualization of a 1-form.
Moreover, by assuming the spinor field to be parallel (given a choice of connection~$\theta$ on~$F$),
the induced almost even-Clifford Hermitian structure will also be parallel,
and we will be able to identify the special Riemannian holonomies from Berger's list.
Since all of our considerations hinge on the existence of such special spinors, we give
explicit representatives for the ranks $r=3 , 7$.

Note that there has been interest on various Clifford-type structures on manifolds (beyond the quaternionic ones) for quite some time (see \cite{AGH,AH,AHS,Barberis,Burdujan1,Burdujan2,Console, Ferus,Friedrich1,GH,Moroianu-Hadfield, Joyce0,Lazaroiu1,Lazaroiu2,Moroianu-Pilca,Moroianu-Semmelmann,Nikolayevsky,Parton0,
Parton1, Parton2, Piccinni, Spindel}). The parallel even-Clifford Hermitian structure resulting from our spinorial construction corresponds (with a minor difference) to the notion of a parallel even Clifford structure introduced in~\cite{Moroianu-Semmelmann}. In particular, Moroianu and Semmelmann studied the relationship between parallel even Clifford structures and special Riemannian holonomy groups, with the exception of~$G_2$. While carrying out the relevant curvature calculations (as in~\cite{Moroianu-Semmelmann}), we noticed that the existence of a parallel twisted pure spinor implies an identity between the curvature of the connection $\theta$
and a 2-form with values in $\ext^2 F$ associated to the spinor (see below), much in the same way as the self-dual part of the ${\rm U}(1)$-connection is related to the
2-form associated to a positive spinor in the Seiberg--Witten monopole equations on 4-manifolds (see \cite{Morgan, Nicolaescu, Witten}). Thus, we introduce the {\em Clifford monopole equations}.

Let $M$ be a ${\rm Spin}^r$ manifold with auxiliary bundle $P_{{\rm SO}(r)}$ endowed with a connection $\theta$, $F$ the associated Riemannian rank $r$ vector bundle
and $m\in\mathbb{N}$ be such that the twisted Dirac operator $\dirac^\theta\colon \Gamma\big(S(M)\otimes S(F)^{\otimes m}\big)\lra \Gamma\big(S(M)\otimes S(F)^{\otimes m}\big)$ is well defined. The {\em Clifford monopole equations} are
\begin{gather*}
 \dirac^\theta\phi = 0,\qquad
\Theta = E\big(\eta^\phi\big),
\end{gather*}
where
$ \Theta \in \Gamma\big(\ext^2T^*M\otimes\ext^2 F\big)$ is the curvature of $\theta$,
\begin{gather*}
 \eta^\phi = \sum_{1\leq k<l \leq r} \eta_{kl}^\phi\otimes f_{kl} \in \Gamma\big(\ext^2T^*M\otimes\ext^2 F\big)
\end{gather*}
is the 2-form with values in $\ext^2 F$ associated to $\phi$, and $E$ is an endomorphism of 2-forms. We will argue that, for suitable choices of the parameters, such equations
are a natural geometric generalization of the Seiberg--Witten equations on 4-manifolds. Indeed, we will show that they restrict to the Seiberg--Witten equations on 4-manifolds, and will exhibit non-trivial solutions (involving parallel twisted pure spinors) on manifolds with special Riemannian holonomy. Preliminary work (with A.~Quintero and to be published elsewhere) indicates the existence of a~smooth compact moduli space which, according to the Mathai--Quillen--Atiyah--Jeffrey for\-ma\-lism, will give raise to a topological quantum field theory. Such topological field theory, at least in dimensions $8$, might turn out to be a topological twist of an $N=2$ supersymmetric theory.

The paper is organized as follows. In Section~\ref{preliminaries}, we recall Clifford algebras, twisted spin groups, representations and structures. In Section~\ref{sec:algebraic-setup}, we define twisted pure spinors, deduce their relevant properties and show explicit representatives. In Section~\ref{sec:holonomy}, we characterize the special Riemannian holonomies by the existence of parallel twisted pure spinor fields. In Section~\ref{sec:Clifford-monopoles}, we show that the Clifford monopole equations restrict to the Seiberg--Witten equations on 4-manifolds, and exhibit solutions on manifolds with special Riemannian holonomy.

\section{Preliminaries}\label{preliminaries}

In this section, we recall various concepts from \cite{Espinosa-Herrera, Friedrich,Lawson} that will be used throughout.

\subsection{Clifford algebras, twisted spin groups and representations}

\subsubsection{Clifford algebra}

Let ${\rm Cl}_n$ denote the Clifford algebra generated by the orthonormal vectors $e_1, e_2, \dots, e_n\in \mathbb{R}^n$ subject to the relations
\begin{gather*}
e_i e_j + e_j e_i = -2\delta_{ij},
\end{gather*}
and $\mathbb{C}{\rm l}_n={\rm Cl}_n\otimes_{\mathbb{R}}\mathbb{C}$ its complexification. The even Clifford subalgebra ${\rm Cl}_n^0$ is defined as the invariant $(+1)$-subspace of the involution of ${\rm Cl}_r$
induced
by the map $-{\rm Id}_{\mathbb{R}^r}$.
We have
\begin{gather*} \mathbb{C}{\rm l}_n\cong
 \begin{cases}
 \operatorname{End}\big(\mathbb{C}^{2^k}\big), & \mbox{if $n=2k$,}\\
 \operatorname{End}\big(\mathbb{C}^{2^k}\big)\oplus\operatorname{End}\big(\mathbb{C}^{2^k}\big), & \mbox{if $n=2k+1$.}
 \end{cases}
\end{gather*}
The space of spinors is defined as
\begin{gather*} \Delta_n:=\mathbb{C}^{2^k}=\underbrace{\mathbb{C}^2\otimes \dots \otimes \mathbb{C}^2}_{k\ {\rm times}}.\end{gather*}
The map
\begin{gather*} \kappa\colon \ \mathbb{C}{\rm l}_n \lra \operatorname{End}(\Delta_n)\end{gather*}
is defined to be either the above mentioned isomorphism for $n$ even, or the isomorphism followed by the projection onto the first summand for~$n$ odd. In order to make $\kappa$ explicit, consider the following matrices
\begin{gather*} {\rm Id} = \left(\begin{matrix}
1 & 0\\
0 & 1
 \end{matrix}\right),\qquad
g_1 = \left(\begin{matrix}
{\rm i} & 0\\
0 & -{\rm i}
 \end{matrix}\right),\qquad
g_2 = \left(\begin{matrix}
0 & {\rm i}\\
{\rm i} & 0
 \end{matrix}\right),\qquad
T = \left(\begin{matrix}
0 & -{\rm i}\\
{\rm i} & 0
 \end{matrix}\right).
\end{gather*}
In terms of the generators $e_1, \dots, e_n$ of the Clifford algebra, $\kappa$ can be
described explicitly as follows,
\begin{gather*}
e_1 \mapsto {\rm Id}\otimes {\rm Id}\otimes \dots\otimes {\rm Id}\otimes {\rm Id}\otimes g_1, \\
e_2 \mapsto {\rm Id}\otimes {\rm Id}\otimes \dots\otimes {\rm Id}\otimes {\rm Id}\otimes g_2, \\
e_3 \mapsto {\rm Id}\otimes {\rm Id}\otimes \dots\otimes {\rm Id}\otimes g_1\otimes T, \\
e_4 \mapsto {\rm Id}\otimes {\rm Id}\otimes \dots\otimes {\rm Id}\otimes g_2\otimes T, \\
\cdots\cdots \cdots\cdots\cdots\cdots\cdots\cdots\cdots\cdots\cdots\cdots \\
e_{2k-1} \mapsto g_1\otimes T\otimes \dots\otimes T\otimes T\otimes T, \\
e_{2k} \mapsto g_2\otimes T\otimes\dots\otimes T\otimes T\otimes T,
\end{gather*}
and, if $n=2k+1$,
\begin{gather*} e_{2k+1}\mapsto {\rm i} T\otimes T\otimes\dots\otimes T\otimes T\otimes T.\end{gather*}
The vectors
\begin{gather*} u_{+1}={1\over \sqrt{2}}(1,-{\rm i})\qquad\mbox{and}\qquad u_{-1}={1\over \sqrt{2}}(1,{\rm i}),\end{gather*}
form a unitary basis of $\mathbb{C}^2$ with respect to the standard Hermitian product. Thus,
\begin{gather*}\mathcal{B}:=\{u_{(\varepsilon_1,\dots,\varepsilon_k)}=u_{\varepsilon_1}\otimes\dots\otimes
u_{\varepsilon_k}\, |\, \varepsilon_j=\pm 1,\, j=1,\dots,k\}\end{gather*}
is a unitary basis of $\Delta_n=\mathbb{C}^{2^k}$ with respect to the naturally induced Hermitian product.

\begin{Remark} We will denote inner and Hermitian products (as well as Riemannian and Hermitian metrics) by the same symbol $\langle \cdot,\cdot\rangle $ trusting that the context will make clear which product is being used.
\end{Remark}

By means of $\kappa$ we have Clifford multiplication
\begin{align*}
\mu_n\colon \ \mathbb{R}^n\otimes \Delta_n &\lra \Delta_n,\\
x \otimes \phi &\mapsto \mu_n(x\otimes \phi)=x\cdot\phi :=\kappa(x)(\phi) .
\end{align*}
$\mu_n$ is skew-symmetric with respect to the Hermitian product
\begin{gather}
\langle x\cdot\phi_1 , \phi_2\rangle =\langle \mu_n(x\otimes \phi_1) , \phi_2\rangle =-\langle \phi_1 , \mu_n(x\otimes \phi_2)\rangle
=-\langle \phi_1 , x\cdot \phi_2\rangle , \label{clifford-skew-symmetric}
\end{gather}
and can be extended to a map
\begin{align*}
\mu_n\colon \ \ext^*\big(\mathbb{R}^n\big)\otimes \Delta_n &\lra \Delta_n,\\
\omega \otimes \psi &\mapsto \omega\cdot\psi.
\end{align*}

When $n$ is even, we define the following involution
\begin{gather*}
\Delta_n \longrightarrow \Delta_n, \\
\psi \mapsto (-{\rm i})^{n\over 2}{\rm vol}_n\cdot \psi.
\end{gather*}
The $\pm 1$ eigenspace of this involution is denoted $\Delta_n^\pm$. These spaces have equal dimension and are irreducible representations of ${\rm Spin}(n)$.

There exist real or quaternionic structures on the spin representations. A quaternionic structure $\alpha$ on $\mathbb{C}^2$ is given by
\begin{gather*}\alpha\left(\begin{matrix}
z_1\\
z_2
 \end{matrix}
\right) = \left(\begin{matrix}
-\overline{z}_2\\
\overline{z}_1
 \end{matrix}\right),\end{gather*}
and a real structure $\beta$ on $\mathbb{C}^2$ is given by
\begin{gather*} \beta\left(\begin{matrix}
z_1\\
z_2
 \end{matrix}
\right) = \left(\begin{matrix}
\overline{z}_1\\
\overline{z}_2
 \end{matrix}\right).\end{gather*}
The real and quaternionic structures $\gamma_n$ on $\Delta_n=\big(\mathbb{C}^2\big)^{\otimes [n/2]}$ are built as follows
\begin{alignat*}{4}
& \gamma_n= (\alpha\otimes\beta)^{\otimes 2k} \qquad &&\text{if $n=8k,8k+1$} \quad && \text{(real)},& \\
& \gamma_n= \alpha\otimes(\beta\otimes\alpha)^{\otimes 2k} \qquad &&\text{if $n=8k+2,8k+3$}\quad && \text{(quaternionic)},& \\
& \gamma_n= (\alpha\otimes\beta)^{\otimes 2k+1} \qquad & &\text{if $n=8k+4,8k+5$}\quad &&\text{(quaternionic)},& \\
& \gamma_n= \alpha\otimes(\beta\otimes\alpha)^{\otimes 2k+1} \qquad &&\text{if $n=8k+6,8k+7$}\quad &&\text{(real)}.&
\end{alignat*}

\subsubsection{Spin group and representation}
The Spin group ${\rm Spin}(n)\subset {\rm Cl}_n$ is the subset
\begin{gather*} {\rm Spin}(n) =\big\{x_1x_2\cdots x_{2l-1}x_{2l} \,|\, x_j\in\mathbb{R}^n, \, |x_j|=1,\,l\in\mathbb{N}\big\},\end{gather*}
endowed with the product of the Clifford algebra. It is a Lie group and its Lie algebra is
\begin{gather*} \mathfrak{spin}(n)=\operatorname{span} \{e_ie_j\, | \,1\leq i< j \leq n\}.\end{gather*}
The restriction of $\kappa$ to ${\rm Spin}(n)$ defines the Lie group representation
\begin{gather*} \kappa_n:=\kappa|_{{\rm Spin}(n)}\colon \ {\rm Spin}(n)\lra {\rm GL}(\Delta_n),\end{gather*}
which is, in fact, special unitary. We have the corresponding Lie algebra representation
\begin{gather*} \kappa_{n_*}\colon \ \mathfrak{spin}(n)\lra \mathfrak{gl}(\Delta_n).\end{gather*}
Both representations can be extended to tensor powers $\kappa_{n*}^m\colon \mathfrak{spin}(n)\longrightarrow \operatorname{End}\big(\Delta_n^{\otimes m}\big)$, $m\in \mathbb{N}$, in the usual way. Recall that the Spin group ${\rm Spin}(n)$ is the universal double cover of ${\rm SO}(n)$, $n\ge 3$. For $n=2$ we consider ${\rm Spin}(2)$ to be the connected double cover of ${\rm SO}(2)$. The covering map will be denoted by
\begin{gather*} \lambda_n\colon \ {\rm Spin}(n)\rightarrow {\rm SO}(n)\subset {\rm GL}\big(\mathbb{R}^n\big).\end{gather*}
Its differential is given
by $\lambda_{n_*}(e_ie_j) = 2E_{ij}$, where $E_{ij}=e_i^*\otimes e_j - e_j^*\otimes e_i$ is the standard basis of the skew-symmetric matrices, and $e^*$ denotes the metric dual of the vector~$e$. Furthermore, we will abuse the notation and also denote by~$\lambda_n$ the induced representation on the exterior algebra~$\ext^*\mathbb{R}^n$.
Note that Clifford multiplication~$\mu_n$ is an equivariant map of~${\rm Spin}(n)$ representations.

Now, we summarize some results about real representations of ${\rm Cl}_r^0$ in Table~\ref{table: irreducible representations even Clifford algebras} (see~\cite{Lawson}). Here~$d_r$ denotes the dimension of an irreducible representation of ${\rm Cl}^0_r$ and $v_r$ the number of distinct irreducible representations. Let $\tilde\Delta_r$ denote the irreducible representation of ${\rm Cl}_r^0$ for $r\not\equiv0$ $(\text{mod}~4) $ and $\tilde\Delta^{\pm}_r$ denote the irreducible representations for $r\equiv0$ $(\text{mod}~4)$.
\begin{table}[h!]\centering
$
\begin{array}{|c|c|c|c|c|c|}
\hline
r(\text{mod}~8)&d_r&{\rm Cl}_r^0&\tilde\Delta_r\, / \,\tilde\Delta_r^\pm& v_r \tstrut\\
\hline
1&2^{\lfloor\frac r2\rfloor}&\mathbb R(d_r)&\mathbb{R}^{d_r}&1 \tstrut\\
\hline
2&2^{\frac r2}&\mathbb C(d_r/2)&\mathbb{C}^{d_r/2}&1 \tstrut\\
\hline
3&2^{\lfloor\frac r2\rfloor+1}&\mathbb H(d_r/4)&\mathbb{H}^{d_r/4}&1 \tstrut\\
\hline
4&2^{\frac r2}&\mathbb H(d_r/4)\oplus \mathbb H(d_r/4)&\mathbb{H}^{d_r/4}&2 \tstrut\\
\hline
5&2^{\lfloor\frac r2\rfloor+1}&\mathbb H(d_r/4)&\mathbb{H}^{d_r/4}&1 \tstrut\\
\hline
6&2^{\frac r2}&\mathbb C(d_r/2)&\mathbb{C}^{d_r/2}&1 \tstrut\\
\hline
7&2^{\lfloor\frac r2\rfloor}&\mathbb R(d_r)&\mathbb{R}^{d_r}&1 \tstrut\\
\hline
8&2^{\frac r2-1}&\mathbb R(d_r)\oplus \mathbb R(d_r)&\mathbb{R}^{d_r}&2 \tstrut\\
\hline
\end{array}
$
\caption{Irreducible representations of even Clifford algebras.}\label{table: irreducible representations even Clifford algebras}
\end{table}

Note that
the representations are complex for $r\equiv 2,6$ $(\text{mod}~8)$ and quaternionic for $r\equiv 3,4,5$
$(\text{mod}~8)$.

\subsubsection{Spinorially twisted spin groups and representations}

By using the unit-length complex numbers ${\rm U}(1)$ or the unit-length quaternions ${\rm Sp}(1)$, the Spin group has been ``twisted'' as follows
\begin{gather*}
{\rm Spin}^c(n) = ({\rm Spin}(n) \times {\rm U}(1))/\{\pm (1,1)\} = {\rm Spin}(n) \times_{\mathbb{Z}_2} {\rm U}(1),\\
 {\rm Spin}^q(n) = ({\rm Spin}(n) \times {\rm Sp}(1))/\{\pm (1,1)\} = {\rm Spin}(n) \times_{\mathbb{Z}_2} {\rm Sp}(1).\end{gather*}
These groups give rise to the following short exact sequences
\begin{gather*}
1\lra \mathbb{Z}_2\lra {\rm Spin}^c(n)\lra {\rm SO}(n)\times {\rm U}(1)\lra 1,\\
1\lra \mathbb{Z}_2\lra {\rm Spin}^q(n)\lra {\rm SO}(n)\times {\rm SO}(3)\lra 1,\end{gather*}
respectively, which lead to the notions of ${\rm Spin}^c$ and ${\rm Spin}^q$ structures \cite{Friedrich,Lawson,Nagase}. Noticing that ${\rm U}(1)={\rm Spin}(2)$ and ${\rm Sp}(1)={\rm Spin}(3)$, we are led to define the twisted Spin group ${\rm Spin}^r(n)$ as follows
\begin{gather*}
{\rm Spin}^{r}(n) = ({\rm Spin}(n) \times {\rm Spin}(r))/\{\pm (1,1)\} = {\rm Spin}(n) \times_{\mathbb{Z}_2} {\rm Spin}(r),
\end{gather*}
where $r\in\mathbb{N}$ and $r\geq 2$. ${\rm Spin}^r(n)$ also fits into a short exact sequence
\begin{gather*} 1\lra \mathbb{Z}_2\lra {\rm Spin}^r(n)\xrightarrow{\lambda_n\times\lambda_r} {\rm SO}(n)\times {\rm SO}(r)\lra 1,\end{gather*}
where
\begin{align*}
 \lambda_n\times\lambda_r\colon \ {\rm Spin}^r(n)&\longrightarrow {\rm SO}(n)\times {\rm SO}(r),\\
 [g,h] &\mapsto (\lambda_n(g),\lambda_r(h)).
\end{align*}
We will call $r$ the {\em rank} of the twisting. Note that the groups ${\rm Spin}^{2}(n)={\rm Spin}^c(n)$ and ${\rm Spin}^{3}(n)={\rm Spin}^q(n)$. The Lie algebra of ${\rm Spin}^r(n)$ is
\begin{gather*} \mathfrak{spin}^r(n) = \mathfrak{spin}(n) \oplus \mathfrak{spin}(r).\end{gather*}
Consider the representations
\begin{align*}
\kappa_{n,r}^m:=\kappa_n\otimes\kappa_r^m\colon \ {\rm Spin}^r(n)&\longrightarrow {\rm GL}\big(\Delta_n\otimes \Delta_r^{\otimes m}\big),\\
 [g,h] &\mapsto \kappa_{n,r}^m([g,h])=\kappa_n(g)\otimes\kappa_r^m(h),
\end{align*}
where $m\in\mathbb{N}$, which are unitary with respect to the Hermitian metric. We will also use the notation
\begin{gather*} [g,h]\cdot(\psi\otimes\varphi) := \kappa_n(g)\otimes\kappa_r^m(h)(\psi\otimes\varphi) =
(\kappa_n(g)(\psi))\otimes\big(\kappa_r^m(h)(\varphi)\big).\end{gather*}
An element $\phi$ of $\Delta_n\otimes\Delta_r^{\otimes m}$ will be called a {\em twisted spinor}, or simply a {\em spinor}.

Also consider the map
\begin{align*}
 \mu_n\otimes\mu_r\colon \ \big(\ext^*\mathbb{R}^n\otimes_\mathbb{R} \ext^*\mathbb{R}^r\big) \otimes_\mathbb{R} (\Delta_n\otimes \Delta_r) &\longrightarrow \Delta_n\otimes\Delta_r,\\
(w_1 \otimes w_2)\otimes (\psi\otimes \varphi) &\mapsto (w_1\otimes w_2)\cdot (\psi\otimes \varphi)= (w_1\cdot\psi) \otimes (w_2\cdot \varphi).
\end{align*}
As in the untwisted case, $\mu_n\otimes\mu_r$ is an equivariant homomorphism of ${\rm Spin}^r(n)$ representations. Note that we can also take tensor products with more copies of $\Delta_r$ as follows
\begin{align*}
\mu_r^a:= {\rm Id}_{\Delta_r}^{\otimes a-1}\otimes \mu_r\otimes {\rm Id}_{\Delta_r}^{\otimes m-a}\colon \
\ext^*\mathbb{R}^r
 \otimes_\mathbb{R} \Delta_r^{m} &\longrightarrow \Delta_r^{m},\\
\beta\otimes ( \varphi_1\otimes\cdots\otimes \varphi_{m}) &\mapsto \varphi_1\otimes\cdots\otimes (\mu_r(\beta\otimes\varphi_a))\otimes\cdots\otimes\varphi_{m}\\
& \hphantom{\mapsto}{} \ = \varphi_1\otimes\cdots\otimes (\beta\cdot\varphi_a)\otimes\cdots\otimes \varphi_{m},
\end{align*}
with Clifford multiplication taking place only in the $a$-th factor. We will also write
\begin{gather*}
\mu_r^a(\beta\otimes \varphi_1\otimes\cdots\otimes \varphi_{m})
 = \mu_r^a(\beta)\cdot (\varphi_1\otimes\cdots\otimes \varphi_{m}) .
\end{gather*}
Notice that if $(f_1,\dots,f_r)$ is an orthonormal frame of $\mathbb{R}^r$,
\begin{gather}
\kappa_{r*}^m(f_kf_l) (\varphi_1\otimes\cdots \otimes\varphi_m)=
\big(\mu_r^1(f_kf_l)\cdot\varphi_1\big)\otimes\cdots \otimes\varphi_m +\cdots\nonumber\\
\hphantom{\kappa_{r*}^m(f_kf_l) (\varphi_1\otimes\cdots \otimes\varphi_m)=}{} + \varphi_1\otimes\cdots \otimes(\mu_r^m(f_kf_l)\cdot\varphi_m).\label{kappa-en-mus}
\end{gather}

\subsection{Spinorially twisted Spin structures}

\subsubsection{Spin structures on oriented Riemannian vector bundles}

Let $F$ be an oriented Riemannian vector bundle over a smooth manifold $M$, with \mbox{$r\!=\!\operatorname{rank}(F)\!\geq\! 3$}. Let $P_{{\rm SO}(r)}$ denote the orthonormal frame bundle of $F$. A Spin structure on $F$ is a~principal ${\rm Spin}(r)$-bundle $P_{{\rm Spin}(r)}$ together with a 2 sheeted covering
\begin{gather*} \Lambda\colon \ P_{{\rm Spin}(r)}\lra P_{{\rm SO}(r)},\end{gather*}
such that $\Lambda(pg)=\Lambda(p)\lambda_r(g)$ for all $p\in P_{{\rm Spin}(r)}$, and all $g\in {\rm Spin}(r)$, where $\lambda_r\colon {\rm Spin}(r)\lra {\rm SO}(r)$ denotes the universal covering map. In the case when $r=\operatorname{rank}(F)=2$, we set $\lambda_2\colon {\rm Spin}(2)\lra {\rm SO}(2)$ to be the connected 2-fold covering of~${\rm SO}(2)$. When $r=1$ a Spin structure is only a 2-fold covering of the base manifold $M$.

Given a Spin structure $P_{{\rm Spin}(r)}$ one can associate a spinor bundle
\begin{gather*} S(F) = P_{{\rm Spin}(r)}\times_{\kappa_r} \Delta_r,\end{gather*}
where $\Delta_r$ denotes the standard complex representation of ${\rm Spin}(r)$. In fact, one can also associate spinor bundles whose fibers are tensor powers of $\Delta_r$,
\begin{gather*} S(F)^{\otimes m}=P_{{\rm Spin}(r)}\times_{\kappa_r^m} \Delta_r^{\otimes m},\end{gather*}
where $m\in\mathbb{N}$.

\subsubsection{Spinorially twisted spin structures on oriented Riemannian manifolds}

\begin{Definition}Let $M$ be an oriented $n$-dimensional Riemannian manifold, $P_{{\rm SO}(n)}$ be its principal bundle of orthonormal frames and $r\in\mathbb{N}$, $r\geq 2$. A~{\em ${\rm Spin}^r$ structure} on $M$ consists of an auxiliary principal ${\rm SO}(r)$-bundle $P_{{\rm SO}(r)}$ and either
\begin{itemize}\itemsep=0pt
\item a principal ${\rm Spin}^r(n)$-bundle $P_{{\rm Spin}^{r}(n)}$ together with an equivariant $2:1$ covering map
\begin{gather*} \Lambda\colon \ P_{{\rm Spin}^{r}(n)}\lra P_{{\rm SO}(n)} \tilde{\times} P_{{\rm SO}(r)},\end{gather*}
where $\tilde{\times}$ denotes the fibre-product, such that $\Lambda(pg)=\Lambda(p)(\lambda_n\times\lambda_r)(g)$ for all $p\in P_{{\rm Spin}^{r}(n)}$ and
$g\in {\rm Spin}^{r}(n)$, where $\lambda_n\times\lambda_r\colon {\rm Spin}^r(n)\lra {\rm SO}(n)\times {\rm SO}(r)$ denotes the canonical $2$-fold cover;
\item or a Spin structure on $TM$.
\end{itemize}

A manifold $M$ admitting a ${\rm Spin}^r$ structure will be called a {\em ${\rm Spin}^r$ manifold}.
\end{Definition}

\begin{Remark}There are three possibilities:
\begin{itemize}\itemsep=0pt
 \item $M$ is a non-Spin ${\rm Spin}^r$ manifold so that the structure group to be considered is ${\rm Spin}^r(n)$.
 \item $M$ and $F$ are both Spin so that the structure group to be considered is ${\rm Spin}(n)\times {\rm Spin}(r)$.
 \item $M$ is Spin and $F$ is not Spin so that the structure groups to be considered is ${\rm Spin}(n)\times {\rm SO}(r)$.
\end{itemize}
A ${\rm Spin}^r$ manifold has associated vector bundles
\begin{gather*} S(M,F,m)=P_{{\rm Spin}^r(n)}\times_{\kappa_n\otimes\kappa_r^{m}}\Delta_n\otimes\Delta_r^{\otimes m},\end{gather*}
where $m\in\mathbb{N}$ may be odd, arbitrary or even, respectively.
\end{Remark}

\begin{Remark} One can also consider the case when $F$ is only locally defined, but $\ext^2 F$ is globally defined, as in the case of some almost quaternionic manifolds.
\end{Remark}

\begin{Remark} There is also a projective case when $r$ is even, i.e., when the auxiliary bundle has structure group $\mathbb{P}{\rm SO}(r)={\rm SO}(r)/\{\pm {\rm Id}_{r\times r}\}$. There are analogous observations for the structure group and the bundles $S(M,F,m)$ in this case.
\end{Remark}

\begin{Example}[\cite{Espinosa-Herrera}] Let us consider the real Grassmannians of oriented $k$-dimensional subspaces of $\mathbb{R}^{k+l}$
\begin{gather*}\mathbb{G}{\rm r}_k\big(\mathbb{R}^{k+l}\big)={{\rm SO}(k+l)\over {\rm SO}(k)\times {\rm SO}(l)}.\end{gather*}
Let $r = ak + bl$, $a,b\in\mathbb{N}$. There exists a homomorphism ${\rm SO}(k)\times {\rm SO}(l)\rightarrow {\rm Spin}^{r}(kl)$ providing a~${\rm Spin}^{r}(kl)$-structure on the real Grassmannian $\mathbb{G}{\rm r}_k\big(\mathbb{R}^{k+l}\big)$ if
\begin{gather*}
 a \equiv l \pmod 2,\\
b \equiv k \pmod 2.
\end{gather*}
\end{Example}

\subsubsection{Covariant derivatives on twisted Spin bundles}

Let $M$ be a ${\rm Spin}^r$ $n$-dimensional manifold and $F$ its auxiliary Riemannian vector bundle of rank~$r$. Assume $F$ is endowed with a covariant derivative $\nabla^F$ (or equivalently, that $P_{{\rm SO}(r)}$ is endowed with a connection 1-form $\theta$) and denote by~$\nabla$ the Levi-Civita covariant derivative on~$M$. These two derivatives induce the spinor covariant derivative
\begin{gather*}\nabla^{\theta}\colon \ \Gamma(S(M,F,m))\lra \Gamma(T^*M\otimes S(M,F,m))\end{gather*}
given locally by
\begin{gather*}\nabla^{\theta} (\psi\otimes \varphi) ={\rm d}(\psi\otimes\varphi)
+\left[{1\over 2}\sum_{1\leq i<j\leq n}\omega_{ji}\otimes e_ie_j\cdot\psi\right] \otimes\varphi\\
\hphantom{\nabla^{\theta} (\psi\otimes \varphi) =}{}
+\psi\otimes\left[{1\over 2}\sum_{1\leq k<l\leq r}\theta_{lk}\otimes \kappa_{r*}^m(f_kf_l)\cdot
\varphi\right] ,\end{gather*}
where $\psi\otimes\varphi\in\Gamma(S(M,F,m))$, $(e_1,\dots,e_n)$ and $(f_1,\dots,f_r)$ are a local orthonormal frames of~$TM$ and~$F$ respectively, $\omega_{ij}$ and $\theta_{kl}$ are the local connection 1-forms. From now on, we will often omit the upper and lower bounds on the indices, by declaring $i$ and $j$ to be the indices for the local vectors of $M$, and $k$ and $l$ to be the indices for the local sections of~$F$. For any tangent vectors $X,Y\in \Gamma(TM)$,{\samepage
\begin{gather}
R^{\theta}(X,Y)(\psi\otimes\varphi)= \left[{1\over 2} \sum_{i<j} \Omega_{ji}(X,Y) e_ie_j\cdot\psi\right]\otimes\varphi\nonumber\\
\hphantom{R^{\theta}(X,Y)(\psi\otimes\varphi)=}{}+ \psi\otimes\left[{1\over 2} \sum_{k<l} \Theta_{lk}(X,Y) \kappa_{r*}^m(f_kf_l)\cdot\varphi\right],\label{eq: basic curvature identity}
\end{gather}
where $\Omega_{ji}$ and $\Theta_{lk}$ are local curvature 2-forms.}

For $X, Y\in\Gamma(TM)$ vector fields and $\phi\in\Gamma(S(M,F,m))$ a spinor field, we also have the compatibility of the covariant derivative with Clifford multiplication,
\begin{gather*}
\nabla^{\theta}_X(Y\cdot\phi) = (\nabla_XY)\cdot\phi + Y\cdot\nabla_X^{\theta}\phi.
\end{gather*}

\subsection{Almost even-Clifford hermitian structures}\label{subsec: almost even-Clifford hermitian structures}
\begin{Definition} Let $n\in \mathbb{N}$ and $(f_1,\dots,f_r)$ be an orthonormal basis of $\mathbb{R}^r$.
\begin{itemize}\itemsep=0pt
\item A {\em linear even-Clifford structure of rank $r$} on $\mathbb{R}^n$ is a homomorphism of associative algebras with unit
\begin{gather*}\Psi\colon \ {\rm Cl}_r^0\longrightarrow \operatorname{End}\big(\mathbb{R}^n\big).\end{gather*}
\item A {\em linear even-Clifford Hermitian structure of rank~$r$} on the inner product space $\big(\mathbb{R}^n,\langle \,,\,\rangle \big)$ is a linear even-Clifford structure of rank~$r$ such that each bivector $f_kf_l$, $1\leq k< l \leq r$, is mapped to a skew-symmetric
endomorphism.
\end{itemize}
\end{Definition}

\begin{Remark}\quad
\begin{itemize}\itemsep=0pt
\item This is also known as the structure of a (left) ${\rm Cl}_r^0$-module.
\item Note that, for $1\leq k<l\leq r$,
\begin{gather*}
(\Psi(f_kf_l))^2=-{\rm Id}_{\mathbb{R}^n}.
\end{gather*}
\item Given a linear even-Clifford structure of rank $r$ on $\mathbb{R}^n$, we can average the standard inner product $\langle \,,\,\rangle $ on $\mathbb{R}^n$ as follows: for $X,Y\in\mathbb{R}^n$,
\begin{gather*}(X,Y)=\sum_{k=1}^{[r/2]} \left[\sum_{1\leq i_1<\dots<i_{2k}<r}
\langle \Psi(f_{i_1}\cdots f_{i_{2k}})(X),\Psi(f_{i_1}\cdots f_{i_{2k}})(Y)\rangle \right], \end{gather*}
so that the linear even-Clifford structure is Hermitian with respect to the averaged inner product.
\item Given a linear even-Clifford Hermitian structure of rank $r$, the subalgebra $\mathfrak{spin}(r)$ is mapped injectively into the skew-symmetric endomorphisms $\operatorname{End}^-\big(\mathbb{R}^n\big)$.
\end{itemize}
\end{Remark}

\begin{Definition} Let $r\geq 2$.
\begin{itemize}\itemsep=0pt
\item A {\em rank $r$ almost even-Clifford structure} on a smooth manifold $M$ is a smoothly varying choice of rank $r$ linear even-Clifford structures on the tangent spaces of~$M$.
\item A smooth manifold carrying an almost even-Clifford structure will be called an {\em almost even-Clifford manifold}.
\item A {\em rank $r$ almost even-Clifford Hermitian structure} on a~Riemannian manifold~$M$ is a~smoo\-th\-ly varying choice of linear even-Clifford Hermitian structures on the tangent spaces of~$M$.
\item A Riemannian manifold carrying an almost even-Clifford Hermitian structure will be called an {\em almost even-Clifford Hermitian manifold}.
\end{itemize}
\end{Definition}

\begin{Remark} Our terminology differs from that of~\cite{Moroianu-Semmelmann}. We have added the words ``almost'' and ``Hermitian'' since, in principle, there is no integrability condition on the structure and the compatibility with a Riemannian metric is a separate condition.
\end{Remark}

\section{Twisted spinors}\label{sec:algebraic-setup}

Throughout this section, we will let $(e_1,\dots, e_n)$ and $(f_1,\dots, f_r)$ be orthonormal bases for~$\mathbb{R}^n$ and~$\mathbb{R}^r$ respectively. A~linear basis for ${\rm Cl}_r^0$ is given by the products $f_{i_1}f_{i_2}\cdots f_{i_{2s}}$, where $\{i_1,i_2,\dots,i_{2s}\}\subset \{1,\dots , r\}$. In order to simplify notation, we will often write
$f_{kl}:=f_kf_l$.

\subsection{2-forms and skew-symmetric endomorphisms associated to a spinor}

\begin{Lemma}\label{lemma:vanishings} Let $\phi\in\Delta_n\otimes\Delta_r^{\otimes m}$, $X,Y\in\mathbb{R}^n$, $1\leq a<b<c<d\leq n$ and $1\leq k,l\leq r$. Then
\begin{gather}
\operatorname{Re}\big\langle \kappa_{r*}^m(f_{kl})\cdot \phi,\phi\big\rangle =0,\label{vanishing2}\\
\operatorname{Re}\langle X\wedge Y\cdot\phi,\phi\rangle =0,\label{vanishing3}\\
\operatorname{Im}\langle X\wedge Y\cdot \kappa_{r*}^m(f_{kl})\cdot\phi,\phi\rangle = 0,\label{vanishing4}\\
\operatorname{Re} \langle X\cdot \phi,Y\cdot\phi \rangle = \langle X,Y\rangle |\phi|^2, \label{real-part}\\
\operatorname{Re} \big\langle e_ae_be_ce_d\cdot \kappa_{r*}^m(f_{kl})\cdot\psi,\psi \big\rangle =0.\label{eeeeff}
\end{gather}
\end{Lemma}
\begin{proof}
By using \eqref{clifford-skew-symmetric} repeatedly
\begin{gather*}
 \langle \mu_{r}^a(f_kf_l)\cdot\phi,\phi\rangle = -\overline{\langle \mu_{r}^a(f_kf_l)\phi,\phi\rangle },
\end{gather*}
so that \eqref{vanishing2} follows from \eqref{kappa-en-mus}.

For identity \eqref{vanishing3}, recall that for $X, Y\in \mathbb{R}^n$
\begin{gather*} X\wedge Y = X\cdot Y + \langle X,Y\rangle .\end{gather*}
Thus
\begin{gather*}
 \langle X\wedge Y\cdot \phi,\phi\rangle = -\overline{\langle X\wedge Y \cdot\phi,\phi\rangle }.
\end{gather*}

Identities \eqref{vanishing4}, \eqref{real-part} and \eqref{eeeeff} follow similarly.
\end{proof}

For any $\xi\in\ext^2\big(\mathbb{R}^n\big)^*$ define $\hat\xi\in\operatorname{End}^-\big(\mathbb{R}^n\big)$ as follows
\begin{gather*}
 X \mapsto \hat\xi(X) := (X\lrcorner \xi)^\sharp =\xi(X,\cdot)^\sharp,
\end{gather*}
where $\lrcorner$ denotes contraction and $^\sharp$ denotes metric dualization.

\begin{Definition}[\cite{Espinosa-Herrera}] Let $\phi\in\Delta_n\otimes\Delta_r^{\otimes m}$, $X,Y\in\mathbb{R}^n$ and $1\leq k,l\leq r$.
\begin{itemize}\itemsep=0pt
\item Let
\begin{gather*} \eta_{kl}^{\phi} (X,Y) = \operatorname{Re}\big\langle X\wedge Y\cdot \kappa_{r*}^m(f_{kl})\cdot \phi,\phi\big\rangle \end{gather*}
be the real $2$-forms associated to the spinor $\phi$.
\item Define the antisymmetric endomorphisms $\hat\eta_{kl}^\phi\in\operatorname{End}^-\big(\mathbb{R}^n\big)$ by
\begin{gather*}X\mapsto \hat\eta_{kl}^\phi(X).\end{gather*}
\end{itemize}
\end{Definition}

\begin{Remark}\quad
\begin{itemize}\itemsep=0pt
 \item For $k\not= l$,
\begin{gather*}\eta_{kl}^\phi = - \eta_{lk}^\phi.\end{gather*}
\item By \eqref{vanishing3},
\begin{gather*}\eta_{kk}\equiv 0.\end{gather*}
 \item By \eqref{vanishing4}, if $k\not= l$,
\begin{gather*}\eta_{kl}^{\phi} (X,Y) =\big\langle X\wedge Y\cdot \kappa_{r*}^m(f_{kl})\cdot \phi,\phi\big\rangle .\end{gather*}
\item For $\lambda\in \mathbb{C}$ with $|\lambda|=1$, the spinor $\lambda\phi$ produces the same $2$-forms
\begin{gather*}\eta_{kl}^{\lambda\phi}= \eta_{kl}^\phi.\end{gather*}
\item Note that, depending on the spinor, such 2-forms can actually be identically zero.
\end{itemize}
\end{Remark}

\begin{Lemma}[\cite{Espinosa-Herrera}]\label{lemma:morphisms}
Any spinor $\phi\in\Delta_n\otimes\Delta_r^{\otimes m}$ defines two maps $($extended by linearity$)$
\begin{align*}
\Phi^\phi\colon \ \ext^2 \mathbb{R}^r&\lra \ext^2 \mathbb{R}^n,\\
f_{kl} &\mapsto \Phi^\phi(f_{kl}):=\eta_{kl}^{\phi},
\end{align*}
and
\begin{align*}
\hat\Phi^\phi\colon \ \ext^2 \mathbb{R}^r&\lra \operatorname{End}\big(\mathbb{R}^n\big),\\
f_{kl} &\mapsto \hat\Phi^\phi(f_{kl}):=\hat\eta_{kl}^{\phi}.
\end{align*}
\end{Lemma}

\subsection[Pure spinors: $r\geq 3$]{Pure spinors: $\boldsymbol{r\geq 3}$}

From now on we shall assume that $r\geq 3$.

\begin{Definition} \label{def:pure-spinor}A (non-zero) spinor $\phi\in\Delta_n\otimes \Delta_r^{\otimes m}$ is called a {\em twisted pure spinor} if
\begin{gather}
\big(\eta_{kl}^\phi + 2 \kappa_{r*}^m(f_{kl})\big)\cdot \phi = 0,\label{eq: pure spinor condition 1}\\
\big(\hat\eta_{kl}^{\phi}\big)^2 = - {\rm Id}_{\mathbb{R}^n},\label{eq: pure spinor condition 2}
\end{gather}
for all $1\leq k<l\leq r$.
\end{Definition}

\begin{Remark}\quad
\begin{itemize}\itemsep=0pt
 \item The first condition says that the subalgebra $\operatorname{span}\big(\eta_{kl}^\phi+2f_{kl}\big) \subset\mathfrak{spin}(n)\oplus\mathfrak{spin}(r)$ annihilates the twisted pure spinor.
 \item The second condition ensures that the 2-forms are non-zero and the associated endomorphisms are almost complex structures.
 \item We will show that the two conditions imply
 \begin{gather*}\operatorname{span}\big\{\hat\eta_{kl}^\phi\in \operatorname{End}\big(\mathbb{R}^n\big)\,|\,1\leq k<l\leq r\big\}\cong\mathfrak{spin}(r).\end{gather*}
\end{itemize}
\end{Remark}

\begin{Lemma}\label{commutators} Let $\phi\in\Delta_n\otimes\Delta_r^{\otimes m}$ be a twisted pure spinor.
\begin{enumerate}\itemsep=0pt
 \item[$1.$] If $1\leq i,j,k,l\leq r$ are all different,
 \begin{gather}
 \big[\hat\eta_{kl}^\phi,\hat\eta_{ij}^\phi\big]=0.\label{[kl,ij]=0}
 \end{gather}
\item[$2.$] If $1\leq i,j,k\leq r$ are all different,
 \begin{gather}
 \big[\hat\eta_{ij}^\phi,\hat\eta_{jk}^\phi\big]=-2 \hat\eta_{ik}^\phi.\label{[ij,jk]=-ik}
 \end{gather}
\end{enumerate}
\end{Lemma}
\begin{proof} For identity \eqref{[kl,ij]=0}, suppose $1\leq i,j,k,l\leq r$ are all different. Notice that in $\mathfrak{spin}(r)\subset {\rm Cl}_r^0$,
\begin{gather*}[f_{kl},f_{ij}]=0\end{gather*}
and, since $\kappa_{r*}^m\colon \mathfrak{spin}(r)\subset {\rm Cl}_r^0 \lra \operatorname{End}(\Delta_r^{\otimes m})$ is a Lie algebra homomorphism,
\begin{gather*}
 0 = \kappa_{r*}^m([f_{kl},f_{ij}]) = \big[\kappa_{r*}^m(f_{kl}),\kappa_{r*}^m(f_{ij})\big],
\end{gather*}
i.e.,
\begin{gather*}\kappa_{r*}^m(f_{kl})\kappa_{r*}^m(f_{ij}) =\kappa_{r*}^m(f_{ij})\kappa_{r*}^m(f_{kl}).\end{gather*}
Now recall that, by definition,
\begin{gather*}\eta_{ij}^\phi\cdot\phi = -2 \kappa_{r*}^m(f_{ij})\cdot \phi,\end{gather*}
which implies
\begin{gather*}
\kappa_{r*}^m(f_{kl})\cdot\eta_{ij}^\phi\cdot \phi = -2 \kappa_{r*}^m(f_{kl})\kappa_{r*}^m(f_{ij})\cdot\phi
= -2 \kappa_{r*}^m(f_{ij})\kappa_{r*}^m(f_{kl})\cdot\phi=\kappa_{r*}^m(f_{ij})\cdot\eta_{kl}^\phi\cdot \phi.
\end{gather*}
By Lemma \ref{lemma:vanishings},
\begin{gather*}
\operatorname{Re}\big\langle e_s\wedge e_t \cdot \eta_{ij}^\phi\cdot \kappa_{r*}^m(f_{kl})\cdot\phi,\phi\big\rangle
 = \operatorname{Re}\langle e_s\wedge e_t \cdot \left(\sum_{a<b}\eta_{ij}^\phi(e_a,e_b)e_a\wedge e_b\right)\cdot
\kappa_{r*}^m(f_{kl})\cdot\phi,\phi\rangle \\
\hphantom{\operatorname{Re}\big\langle e_s\wedge e_t \cdot \eta_{ij}^\phi\cdot \kappa_{r*}^m(f_{kl})\cdot\phi,\phi\big\rangle}{}
= \operatorname{Re}\sum_{a<b}\eta_{ij}^\phi(e_a,e_b)\langle e_s\cdot e_t \cdot e_a\cdot e_b\cdot
\kappa_{r*}^m(f_{kl})\cdot\phi,\phi\rangle \\
\hphantom{\operatorname{Re}\big\langle e_s\wedge e_t \cdot \eta_{ij}^\phi\cdot \kappa_{r*}^m(f_{kl})\cdot\phi,\phi\big\rangle}{}
= \sum_{s=a<b}\eta_{ij}^\phi(e_s,e_b)\langle e_t \cdot e_b\cdot
\kappa_{r*}^m(f_{kl})\cdot\phi,\phi\rangle \\
\hphantom{\operatorname{Re}\big\langle e_s\wedge e_t \cdot \eta_{ij}^\phi\cdot \kappa_{r*}^m(f_{kl})\cdot\phi,\phi\big\rangle=}{}
 +\sum_{t=a<b}\eta_{ij}^\phi(e_t,e_b)(-\langle e_s\cdot e_b\cdot
\kappa_{r*}^m(f_{kl})\cdot\phi,\phi\rangle )\\
\hphantom{\operatorname{Re}\big\langle e_s\wedge e_t \cdot \eta_{ij}^\phi\cdot \kappa_{r*}^m(f_{kl})\cdot\phi,\phi\big\rangle=}{}
 +\sum_{a<t=b}\eta_{ij}^\phi(e_a,e_t)\langle e_s\cdot e_a\cdot
\kappa_{r*}^m(f_{kl})\cdot\phi,\phi\rangle \\
\hphantom{\operatorname{Re}\big\langle e_s\wedge e_t \cdot \eta_{ij}^\phi\cdot \kappa_{r*}^m(f_{kl})\cdot\phi,\phi\big\rangle=}{}
 +\sum_{a<b=s}\eta_{ij}^\phi(e_a,e_s)(-\langle e_t \cdot e_a\cdot
\kappa_{r*}^m(f_{kl})\cdot\phi,\phi\rangle )\\
\hphantom{\operatorname{Re}\big\langle e_s\wedge e_t \cdot \eta_{ij}^\phi\cdot \kappa_{r*}^m(f_{kl})\cdot\phi,\phi\big\rangle}{}
= \sum_{s<b}\eta_{ij}^\phi(e_s,e_b)\eta_{kl}^\phi(e_t,e_b) +\sum_{t<b}\eta_{ij}^\phi(e_t,e_b)\big({-}\eta_{kl}^\phi(e_s,e_b)\big)\\
\hphantom{\operatorname{Re}\big\langle e_s\wedge e_t \cdot \eta_{ij}^\phi\cdot \kappa_{r*}^m(f_{kl})\cdot\phi,\phi\big\rangle=}{} +\sum_{b<t}\eta_{ij}^\phi(e_b,e_t)\eta_{kl}^\phi(e_s,e_b) +\sum_{b<s}\eta_{ij}^\phi(e_b,e_s)\big({-}\eta_{kl}^\phi(e_t,e_b)\big)\\
 \hphantom{\operatorname{Re}\big\langle e_s\wedge e_t \cdot \eta_{ij}^\phi\cdot \kappa_{r*}^m(f_{kl})\cdot\phi,\phi\big\rangle }{}
 = -\sum_{b}\eta_{ij}^\phi(e_s,e_b)\eta_{kl}^\phi(e_b,e_t)
+\sum_{b}\eta_{kl}^\phi(e_s,e_b)\eta_{ij}^\phi(e_b,e_t)\\
\hphantom{\operatorname{Re}\big\langle e_s\wedge e_t \cdot \eta_{ij}^\phi\cdot \kappa_{r*}^m(f_{kl})\cdot\phi,\phi\big\rangle}{}
= -\sum_{b}\big(\hat\eta_{kl}^\phi\big)_{tb}\big(\hat\eta_{ij}^\phi\big)_{bs}
+\sum_{b}\big(\hat\eta_{ij}^\phi\big)_{tb}\big(\hat\eta_{kl}^\phi\big)_{bs}\\
\hphantom{\operatorname{Re}\big\langle e_s\wedge e_t \cdot \eta_{ij}^\phi\cdot \kappa_{r*}^m(f_{kl})\cdot\phi,\phi\big\rangle}{}
= \big[\hat\eta_{ij}^\phi,\hat\eta_{kl}^\phi\big]_{ts},
\end{gather*}
the entry in row $t$ and column $s$ of the matrix
\begin{gather*}\big[\hat\eta_{kl}^\phi,\hat\eta_{ij}^\phi\big].\end{gather*}
Analogously,
\begin{gather*}
\operatorname{Re}\big\langle e_s\wedge e_t \cdot \eta_{kl}^\phi\cdot \kappa_{r*}^m(f_{ij})\cdot\phi,\phi\big\rangle
 = \big[\hat\eta_{kl}^\phi,\hat\eta_{ij}^\phi\big]_{ts}.
\end{gather*}
Thus,
\begin{gather*}\big[\hat\eta_{ij}^\phi,\hat\eta_{kl}^\phi\big]=\big[\hat\eta_{kl}^\phi,\hat\eta_{ij}^\phi\big],\end{gather*}
but by definition of the bracket
\begin{gather*}\big[\hat\eta_{kl}^\phi,\hat\eta_{ij}^\phi\big]=-\big[\hat\eta_{ij}^\phi,\hat\eta_{kl}^\phi\big].\end{gather*}
Hence,
\begin{gather*}\big[\hat\eta_{ij}^\phi,\hat\eta_{kl}^\phi\big]=0.\end{gather*}

For identity \eqref{[ij,jk]=-ik}, recall that in $\mathfrak{spin}(r)\subset {\rm Cl}_r^0$,
\begin{gather*}
[f_{ij},f_{jk}] = f_{ij}f_{jk}-f_{jk}f_{ij}= -2f_{ik},
\end{gather*}
so that
\begin{gather*}
 -2\kappa_{r*}^m(f_{ik}) =\kappa_{r*}^m([f_{ij},f_{jk}])= \big[\kappa_{r*}^m(f_{ij}),\kappa_{r*}^m(f_{jk})\big],
\end{gather*}
i.e.,
\begin{gather*}\kappa_{r*}^m(f_{ij})\kappa_{r*}^m(f_{jk})
=\kappa_{r*}^m(f_{jk})\kappa_{r*}^m(f_{ij})-2\kappa_{r*}^m(f_{ik}).\end{gather*}
Now,
\begin{gather*}\eta_{ij}^\phi\cdot\phi = -2 \kappa_{r*}^m(f_{ij})\cdot \phi,\end{gather*}
which implies
\begin{gather*}
\kappa_{r*}^m(f_{jk})\cdot\eta_{ij}^\phi\cdot \phi = -2 \kappa_{r*}^m(f_{jk})\kappa_{r*}^m(f_{ij})\cdot
\phi= -2 \big[\kappa_{r*}^m(f_{ij})\kappa_{r*}^m(f_{jk})+2\kappa_{r*}^m(f_{ik})\big]\cdot \phi\\
\hphantom{\kappa_{r*}^m(f_{jk})\cdot\eta_{ij}^\phi\cdot \phi}{}
= \kappa_{r*}^m(f_{ij})\cdot\eta_{jk}^\phi\cdot \phi - 4 \kappa_{r*}^m(f_{ik})\cdot \phi.
\end{gather*}
Thus, on the one hand,
\begin{gather*}
\operatorname{Re}\big\langle e_s\wedge e_t \cdot \eta_{ij}^\phi\cdot \kappa_{r*}^m(f_{jk})\cdot\phi,\phi\big\rangle\\
\qquad{} = \operatorname{Re}\big\langle e_s\wedge e_t \cdot \eta_{jk}^\phi\cdot
\kappa_{r*}^m(f_{ij})\cdot\phi,\phi\big\rangle - 4\operatorname{Re}\big\langle e_s\wedge e_t\cdot
\kappa_{r*}^m(f_{ik})\cdot\phi \big\rangle \\
 \qquad = \operatorname{Re}\big\langle e_s\wedge e_t \cdot \eta_{jk}^\phi\cdot
\kappa_{r*}^m(f_{ij})\cdot\phi,\phi\big\rangle - 4 \eta_{ik}^\phi(e_s,e_t).
\end{gather*}
By the calculation above
\begin{gather*}
\operatorname{Re}\big\langle e_s\wedge e_t \cdot \eta_{ij}^\phi\cdot \kappa_{r*}^m(f_{jk})\cdot\phi,\phi\big\rangle
 = \big[\hat\eta_{ij}^\phi,\hat\eta_{jk}^\phi\big]_{ts},\\
\operatorname{Re}\big\langle e_s\wedge e_t \cdot \eta_{jk}^\phi\cdot \kappa_{r*}^m(f_{ij})\cdot\phi,\phi\big\rangle= \big[\hat\eta_{jk}^\phi,\hat\eta_{ij}^\phi\big]_{ts},\\
\eta_{ik}^\phi(e_s,e_t) = \big(\hat\eta_{ik}^\phi\big)_{ts},
\end{gather*}
so that
\begin{gather*}
\big[\hat\eta_{ij}^\phi,\hat\eta_{jk}^\phi \big] =\big[\hat\eta_{jk}^\phi,\hat\eta_{ij}^\phi\big] - 4\hat\eta_{ik}^\phi
 = - \big[\hat\eta_{ij}^\phi,\hat\eta_{jk}^\phi \big] - 4\hat\eta_{ik}^\phi,
\end{gather*}
and
\begin{gather*}
2 \big[\hat\eta_{ij}^\phi,\hat\eta_{jk}^\phi\big] = -4\hat\eta_{ik}^\phi.\tag*{\qed}
\end{gather*}\renewcommand{\qed}{}
\end{proof}

\begin{Remark} For a twisted spinor $\phi$ satisfying only condition \eqref{eq: pure spinor condition 1}, the endomorphisms $2\hat\eta_{kl}^\phi$ satisfy the Lie bracket relations of the Lie algebra~$\mathfrak{so}(r)$.
\end{Remark}

\begin{Lemma}\label{lemma:clifford-relations} Let $\phi\in\Delta_n\otimes\Delta_r^{\otimes m}$ be a~twisted pure spinor. Let $1\leq i,j,k,l\leq r$ be all different.
\begin{itemize}\itemsep=0pt
 \item The automorphisms $\hat\eta_{ij}^\phi$ and $\hat\eta_{kl}^\phi$ commute
\begin{gather}
\hat\eta_{ij}^\phi\hat\eta_{kl}^\phi=\hat\eta_{kl}^\phi\hat\eta_{ij}^\phi.\label{ij.kl=kl.ij}
\end{gather}
\item The automorphisms $\hat\eta_{ij}^\phi$
and $\hat\eta_{jk}^\phi$ anticommute
\begin{gather}
\hat\eta_{ij}^\phi\hat\eta_{jk}^\phi =-\hat\eta_{jk}^\phi\hat\eta_{ij}^\phi =- \hat\eta_{ik}^\phi.\label{ij.jk=-jk.ij}
\end{gather}
\item The following identities hold
\begin{gather}
 \hat\eta_{ij}^\phi\hat\eta_{kl}^\phi = -\hat\eta_{ik}^\phi\hat\eta_{jl}^\phi
 = -\hat\eta_{jl}^\phi\hat\eta_{ik}^\phi = \hat\eta_{kl}^\phi\hat\eta_{ij}^\phi = \hat\eta_{jk}^\phi\hat\eta_{il}^\phi
 = \hat\eta_{il}^\phi\hat\eta_{jk}^\phi.\label{eq:ijkl=-ikjl}
\end{gather}
\end{itemize}
\end{Lemma}

\begin{proof} Identity \eqref{ij.kl=kl.ij} is the same as \eqref{[kl,ij]=0} in Lemma~\ref{commutators}.

For identity \eqref{ij.jk=-jk.ij} recall
\begin{gather*}\big(\hat\eta_{ij}^\phi\big)^2 = - {\rm Id}_{\mathbb{R}^n},\end{gather*}
and the identity
\begin{gather*}\hat\eta_{ij}^\phi\hat\eta_{jk}^\phi -\hat\eta_{jk}^\phi\hat\eta_{ij}^\phi= -2\,\hat\eta_{ik}^\phi.\end{gather*}
Compose the last identity on the left and on the right with $\hat\eta_{ij}^\phi$
\begin{gather*}\big(\hat\eta_{ij}^\phi\big)^2\hat\eta_{jk}^\phi\hat\eta_{ij}^\phi
-\hat\eta_{ij}^\phi\hat\eta_{jk}^\phi\big(\hat\eta_{ij}^\phi\big)^2=
-2 \hat\eta_{ij}^\phi\hat\eta_{ik}^\phi\hat\eta_{ij}^\phi,\end{gather*}
so that
\begin{gather*}-\hat\eta_{jk}^\phi\hat\eta_{ij}^\phi +\hat\eta_{ij}^\phi\hat\eta_{jk}^\phi=
-2 \hat\eta_{ij}^\phi\hat\eta_{ik}^\phi\hat\eta_{ij}^\phi.\end{gather*}
Thus,
\begin{gather*}-2\hat\eta_{ik}^\phi= -2 \hat\eta_{ij}^\phi\hat\eta_{ik}^\phi\hat\eta_{ij}^\phi,\end{gather*}
and
\begin{gather*}\hat\eta_{ik}^\phi\hat\eta_{ij}^\phi= \hat\eta_{ij}^\phi\hat\eta_{ik}^\phi\big(\hat\eta_{ij}^\phi\big)^2,\end{gather*}
i.e.,
\begin{gather*}\hat\eta_{ik}^\phi\hat\eta_{ij}^\phi= -\hat\eta_{ij}^\phi\hat\eta_{ik}^\phi.\end{gather*}
Hence,
\begin{gather*}
-2 \hat\eta_{ik}^\phi
 = \big[\hat\eta_{ij}^\phi,\hat\eta_{jk}^\phi\big]
 = \hat\eta_{ij}^\phi\hat\eta_{jk}^\phi -\hat\eta_{jk}^\phi\hat\eta_{ij}^\phi
 = \hat\eta_{ij}^\phi\hat\eta_{jk}^\phi -\big({-}\hat\eta_{ij}^\phi\hat\eta_{jk}^\phi\big)
 = 2\hat\eta_{ij}^\phi\hat\eta_{jk}^\phi.
\end{gather*}

For \eqref{eq:ijkl=-ikjl}, we can see that
\begin{gather*}
 \hat\eta_{ij}^\phi\hat\eta_{kl}^\phi = \hat\eta_{ik}^\phi\hat\eta_{jk}^\phi \hat\eta_{kl}^\phi = -\hat\eta_{ik}^\phi\hat\eta_{jl}^\phi,
\end{gather*}
and similarly for the remaining identities.
\end{proof}

\begin{Lemma}\label{lemma: independent of frame} The definition of twisted pure spinor does not depend on the choice of orthonormal basis $(f_1,\dots,f_r)$ of $\mathbb{R}^r$.
\end{Lemma}
\begin{proof} Suppose $(f_1'\dots,f_r')$ is another orthonormal basis of $\mathbb{R}^r$ so that
\begin{gather*} f_k' = a_{k1}f_1 + \cdots + a_{kr}f_r ,\end{gather*}
for $1\leq k \leq r$, and the matrix $A=(a_{kl})\in {\rm SO}(r)$.
Recall that
\begin{gather*}\eta_{kl}^\phi = \Phi^\phi(f_{kl}).\end{gather*}
If we write the left-hand side of \eqref{eq: pure spinor condition 1}
 with respect to the basis $(f_1',\dots,f_r')$, we have
\begin{gather*}
 \big(\Phi^\phi(f_{kl}') + 2\kappa_{r_*}^m (f_{kl}')\big)\cdot \phi\\
 \qquad{} = \left(\left(\sum_{s< t}(a_{ks}a_{lt} -a_{kt}a_{ls})
\Phi^\phi(f_{st})\right) + 2\kappa_{r_*}^m \left(\sum_{s< t}(a_{ks}a_{lt} -a_{kt}a_{ls})
f_{st}\right) \right) \cdot \phi\\
\qquad{} = \sum_{s< t}(a_{ks}a_{lt} -a_{kt}a_{ls})
(\Phi^\phi(f_{st}) + 2\kappa_{r_*}^m (f_{st})) \cdot \phi = 0.
\end{gather*}
In order to simplify notation, let
\begin{gather*}
 J_{kl} = \hat\Phi^\phi(f_{kl}), \qquad J'_{kl} = \hat\Phi^\phi(f'_{kl}).
\end{gather*}
Now suppose that the second condition of pure spinor is fulfilled for the frame $(f_1,\dots ,f_r)$
\begin{gather*}
J_{kl}^2= -{\rm Id}_{\mathbb{R}^n}.\end{gather*}
With respect to an orthonormal basis $(e_1,\dots,e_n)$ of $\mathbb{R}^n$,
\begin{align}
J'_{kl}(X) &= \sum_{c=1}^n \Phi^\phi (f'_{kl})(X,e_c)e_c= \sum_{c=1}^n \sum_{1\leq s< t\leq r}(a_{ks}a_{lt} -a_{kt}a_{ls})
\Phi^\phi(f_{st})(X,e_c)e_c\nonumber\\
 &= \sum_{1\leq s< t\leq r}(a_{ks}a_{lt} -a_{kt}a_{ls})\sum_{c=1}^n\Phi^\phi(f_{st})(X,e_c)e_c\nonumber\\
 &= \sum_{1\leq s< t\leq r}(a_{ks}a_{lt} -a_{kt}a_{ls}) \label{eq:J' and J}
J_{st}(X).
\end{align}
Since the bases $\{f_{kl}\,|\,1\leq k< l\leq r\}$ and $\{f'_{kl}\,|\,1\leq k< l\leq r\}$ are orthonormal in $\ext^2 \mathbb{R}^r$,
\begin{align}
\delta_{ac}\delta_{bd}=\langle f'_{ab},f'_{cd}\rangle& =\bigg\langle \sum_{1\leq s< t\leq r}(a_{as}a_{bt} -a_{at}a_{bs})f_{st},
\sum_{1\leq u< v\leq r}(a_{cu}a_{dv} -a_{cv}a_{du})f_{uv}\bigg\rangle \nonumber\\
 &=\sum_{1\leq s< t\leq r}\sum_{1\leq u< v\leq r}(a_{as}a_{bt} -a_{at}a_{bs})
(a_{cu}a_{dv} -a_{cv}a_{du})\delta_{su}\delta_{tv}\nonumber\\
 &=\sum_{1\leq s< t\leq r}(a_{as}a_{bt} -a_{at}a_{bs})(a_{cs}a_{dt} -a_{ct}a_{ds}). \label{eq:orthononality-ext2}
\end{align}
By \eqref{eq:J' and J},
\begin{gather*}
 J'_{kl} = \sum_{1\leq s< t\leq r}(a_{ks}a_{lt} -a_{kt}a_{ls})J_{st},
\end{gather*}
we have
\begin{align*}
 J'_{kl} J'_{kl} &= \left(\sum_{1\leq s< t\leq r}(a_{ks}a_{lt} -a_{kt}a_{ls})
J_{st}\right) \left(\sum_{1\leq u< v\leq r}(a_{ku}a_{lv} -a_{kv}a_{lu})J_{uv}\right)\\
 &=\sum_{1\leq s< t\leq r}\sum_{1\leq u< v\leq r}(a_{ks}a_{lt} -a_{kt}a_{ls})(a_{ku}a_{lv} -a_{kv}a_{lu})J_{st}J_{uv}.
\end{align*}
There are three cases:
\begin{itemize}\itemsep=0pt
 \item[(i)] the indices $s,t,u,v$ are all different;
 \item[(ii)] the pairs $(s,t)$ and $(u,v)$ have one, and only one, common entry;
 \item[(iii)] the pairs $(s,t)$ and $(u,v)$ coincide.
\end{itemize}

For (i), note that since $s<t$ and $u<v$, we only have the following six summands with those indices, so that
\begin{gather*}
 (a_{k s} a_{l t} - a_{k t} a_{l s}) (a_{k u} a_{l v} - a_{k v} a_{l u}) J_{s t} J_{u v}
+ (a_{k s} a_{l u} - a_{k u} a_{l s}) (a_{k t} a_{l v} - a_{k v} a_{l t}) J_{s u} J_{t v} \\
\qquad {} + (a_{k s} a_{l u} - a_{k u} a_{l s})
 (a_{k t} a_{l v} - a_{k v} a_{l t}) J_{s u} J_{t v} + (a_{k s} a_{l t} - a_{k t} a_{l s})
 (a_{k u} a_{l v} - a_{k v} a_{l u}) J_{s t} J_{u v} \\
\qquad {} + (a_{k s} a_{l v} - a_{k v} a_{l s})
 (a_{k t} a_{l u} - a_{k u} a_{l t}) J_{s v} J_{t u}
+ (a_{k s} a_{l v} - a_{k v} a_{l s})
 (a_{k t} a_{l u} - a_{k u} a_{l t}) J_{s v} J_{t u}\\
\quad{} = (2(a_{ks}a_{lt}-a_{kt}a_{ls})(a_{ku}a_{lv}-a_{kv}a_{lu}) -2(a_{ks}a_{lu}-a_{ku}a_{ls})(a_{kt}a_{l v }-a_{kv}a_{lt})\\
 \qquad{} +2(a_{ks}a_{lv}-a_{kv}a_{ls})(a_{kt}a_{lu}-a_{ku}a_{lt}))J_{s t} J_{u v} = 0.
\end{gather*}

For (ii), suppose $s=u$ but $t\not=v$. Now we have two summands
\begin{gather*}
 (a_{k s} a_{l t} - a_{k t} a_{l s}) (a_{k s} a_{l v} - a_{k v} a_{l s}) J_{s t} J_{s v}
 + (a_{k s} a_{l v} - a_{k v} a_{l s}) (a_{k s} a_{l t} - a_{k t} a_{l s}) J_{s v} J_{s t}\\
 \qquad{} = (a_{k s} a_{l t} - a_{k t} a_{l s}) (a_{k s} a_{l v} - a_{k v} a_{l s})
 (J_{s t} J_{s v}+J_{s v} J_{s t}) \\
 \qquad{} = (a_{k s} a_{l t} - a_{k t} a_{l s}) (a_{k s} a_{l v} - a_{k v} a_{l s})
 (J_{tv}+J_{vt}) \\
 \qquad{} = (a_{k s} a_{l t} - a_{k t} a_{l s}) (a_{k s} a_{l v} - a_{k v} a_{l s}) (J_{tv}-J_{tv}) =0,
\end{gather*}
and similarly for the other cases.

For (iii), we have
\begin{gather*}
 \sum_{1\leq s< t\leq r}(a_{ks}a_{lt} -a_{kt}a_{ls})^2 J_{st}^2 = -{\rm Id}_{\mathbb{R}^n},
\end{gather*}
where we have used \eqref{eq:orthononality-ext2}.
\end{proof}

\begin{Proposition}\label{prop: Clifford transfer} A twisted pure spinor $\phi\in\Delta_n\otimes\Delta_r^{\otimes m}$ induces a linear even-Clifford Hermitian structure
of rank~$r$ on~$\mathbb{R}^n$
\begin{align*}
{\rm Cl}_r^0&\lra \operatorname{End}\big(\mathbb{R}^n\big), \\
f_{ij}&\mapsto \hat\eta_{ij}^\phi,
 \end{align*}
so that
\begin{gather*}\mathbb{R}^n \cong
\begin{cases}
\mathbb{R}^m\otimes \tilde\Delta_r & \text{if $r\not\equiv 0$ {\rm (mod 4)}},\\
\mathbb{R}^{m_1}\otimes\tilde\Delta_r^+ \oplus \mathbb{R}^{m_2}\otimes\tilde\Delta_r^- &
\text{if $r\equiv 0$ {\rm (mod 4)}},\end{cases}
\end{gather*}
as a representation of ${\rm Cl}_r^0$, for some $m, m_1, m_2\in\mathbb{N}$, where $\tilde\Delta_r$ denotes the $($unique$)$ non-trivial real representation of ${\rm Cl}_r^0$ if $r\not\equiv 0$ {\rm (mod} $4)$, and $\tilde\Delta_r^+$ and $\tilde\Delta_r^-$ denote the two non-trivial real representations of ${\rm Cl}_r^0$ if $r\equiv 0$ {\rm (mod} $4)$.
\end{Proposition}
\begin{proof} By Lemma \ref{lemma:clifford-relations}, the map
\begin{align*}
\big({\rm Cl}_r^0\big)^2&\lra \operatorname{End}^-\big(\mathbb{R}^n\big),\\
f_{ij} &\mapsto \hat\eta_{ij}^\phi
\end{align*}
extends to an algebra homomorphism
\begin{gather*}
{\rm Cl}_r^0 \lra \operatorname{End}\big(\mathbb{R}^n\big).
\end{gather*}
Since the matrices $\hat\eta_{ij}^\phi$ square to $-{\rm Id}_{\mathbb{R}^n}$, this representation of~${\rm Cl}_r^0$ contains no trivial summands. By \cite[Theorem~5.6]{Lawson}, we know that the algebra ${\rm Cl}_r^0$ has (up to isomorphism) only one or two non-trivial irreducible representations depending on whether $r\not\equiv 0$ ({\rm mod} 4) or $r\equiv 0$ ({\rm mod} 4) respectively.
\end{proof}

\begin{Remark} Note that, unlike \cite{Moroianu-Semmelmann}, in our case $\hat\eta_{kk}^\phi=0$.\end{Remark}

\begin{Lemma}\label{lemma: orbit of pure spinor} Let $\phi\in\Delta_n\otimes\Delta_r^{\otimes m}$ be a twisted pure spinor and $[g,h]\in {\rm Spin}^{r}(n)$. The spinor $\kappa_{n,r}^m([g,h])(\phi)$ is also a twisted pure spinor.
\end{Lemma}
\begin{proof} Consider the orthonormal bases
\begin{gather*}
(e'_1,\dots,e'_n)=(\lambda_n(g)(e_1),\dots,\lambda_n(g)(e_n)),\\
(f'_1,\dots,f'_r)=(\lambda_r(h)(f_1),\dots,\lambda_r(h)(f_r)),
\end{gather*}
of $\mathbb{R}^n$ and $\mathbb{R}^r$ respectively.
We will verify the pure spinor identities for $\varphi:=\kappa_{n,r}^m([g,h])(\phi)$ using these bases.
Indeed,{\samepage
\begin{gather*}
 \Phi^\varphi(f'_kf'_l) \cdot \varphi = \sum_{a<b} \langle e'_ae'_b\cdot \kappa_{r*}^m(f'_kf'_l)\cdot \varphi,\varphi\rangle e'_ae'_b\cdot\varphi\\
 \hphantom{\Phi^\varphi(f'_kf'_l) \cdot \varphi}{} = \sum_{a<b} \big\langle \lambda_n(g)(e_a)\lambda_n(g)(e_b)\cdot
\kappa_{r*}^m(\lambda_r(h)(f_k)\lambda_r(h)(f_l))\cdot \varphi,\varphi\big\rangle\\
\hphantom{\Phi^\varphi(f'_kf'_l) \cdot \varphi=}{}\times \lambda_n(g)(e_a)\lambda_n(g)(e_b)\cdot\varphi\\
\hphantom{\Phi^\varphi(f'_kf'_l) \cdot \varphi}{}
 = \sum_{a<b} \big\langle \lambda_n(g)(e_ae_b)\cdot
\kappa_{r*}^m(\lambda_r(h)(f_kf_l))\cdot \varphi,\varphi\big\rangle
\lambda_n(g)(e_ae_b)\cdot\varphi\\
\hphantom{\Phi^\varphi(f'_kf'_l) \cdot \varphi}{}
 = \sum_{a<b} \big\langle \lambda_n\times\lambda_r([g,h])\big(e_ae_b\cdot
\kappa_{r*}^m(f_kf_l)\big)\cdot \kappa_{n,r}^m([g,h])(\phi),\kappa_{n,r}^m([g,h])(\phi)\big\rangle\\
\hphantom{\Phi^\varphi(f'_kf'_l) \cdot \varphi=}{}\times
\lambda_n(g)(e_ae_b)\cdot\kappa_{n,r}^m([g,h])(\phi)\\
\hphantom{\Phi^\varphi(f'_kf'_l) \cdot \varphi}{}
 = \sum_{a<b} \big\langle \kappa_{n,r}^m([g,h])(e_ae_b\cdot
\kappa_{r*}^m(f_kf_l)\cdot \phi),\kappa_{n,r}^m([g,h])(\phi)\big\rangle\\
\hphantom{\Phi^\varphi(f'_kf'_l) \cdot \varphi=}{}\times
\lambda_n(g)(e_ae_b)\cdot\kappa_{n,r}^m([g,h])(\phi)\\
 \hphantom{\Phi^\varphi(f'_kf'_l) \cdot \varphi}{}
 = \sum_{a<b} \big\langle e_ae_b\cdot
\kappa_{r*}^m(f_kf_l)\cdot \phi,\phi\big\rangle
\kappa_{n,r}^m([g,h])(e_ae_b\cdot\phi)\\
\hphantom{\Phi^\varphi(f'_kf'_l) \cdot \varphi}{}
 = \kappa_{n,r}^m([g,h])\left(\sum_{a<b} \big\langle e_ae_b\cdot
\kappa_{r*}^m(f_kf_l)\cdot \phi,\phi\big\rangle
e_ae_b\cdot\phi\right)\\
\hphantom{\Phi^\varphi(f'_kf'_l) \cdot \varphi}{}
 = \kappa_{n,r}^m([g,h])\big(\Phi^\phi(f_kf_l)\cdot\phi\big)\\
 \hphantom{\Phi^\varphi(f'_kf'_l) \cdot \varphi}{}
 = \kappa_{n,r}^m([g,h])\big({-}2\kappa_{r*}^m(f_kf_l)\cdot\phi\big)\\
 \hphantom{\Phi^\varphi(f'_kf'_l) \cdot \varphi}{}
 = -2\kappa_{r*}^m(\lambda_r(h)(f_kf_l))\cdot\kappa_{n,r}^m([g,h])(\phi)\\
 \hphantom{\Phi^\varphi(f'_kf'_l) \cdot \varphi}{}
 = -2\kappa_{r*}^m(f'_kf'_l)\cdot\varphi,
\end{gather*}
which proves condition \eqref{eq: pure spinor condition 1} for $\varphi$.}

For condition \eqref{eq: pure spinor condition 2}, consider
\begin{gather*}
 \Phi^\varphi(f'_kf'_l) = \sum_{a<b} \big\langle e'_ae'_b\cdot \kappa_{r*}^m(f'_kf'_l)\cdot \varphi,\varphi\big\rangle e'_ae'_b\\
 \hphantom{\Phi^\varphi(f'_kf'_l)}{}= \sum_{a<b} \big\langle \lambda_n(g)(e_a)\lambda_n(g)(e_b)\cdot
\kappa_{r*}^m(\lambda_r(h)(f_k)\lambda_r(h)(f_l))\cdot \varphi,\varphi\big\rangle
e'_ae'_b\\
\hphantom{\Phi^\varphi(f'_kf'_l)}{} = \sum_{a<b} \big\langle \lambda_n(g)(e_ae_b)\cdot
\kappa_{r*}^m(\lambda_r(h)(f_kf_l))\cdot \varphi,\varphi\big\rangle e'_ae'_b\\
\hphantom{\Phi^\varphi(f'_kf'_l)}{} = \sum_{a<b} \big\langle \lambda_n\times\lambda_r([g,h])(e_ae_b\cdot
\kappa_{r*}^m(f_kf_l))\cdot \kappa_{n,r}^m([g,h])(\phi),\kappa_{n,r}^m([g,h])(\phi)\big\rangle e'_ae'_b\\
\hphantom{\Phi^\varphi(f'_kf'_l)}{} = \sum_{a<b} \big\langle \kappa_{n,r}^m([g,h])(e_ae_b\cdot
\kappa_{r*}^m(f_kf_l)\cdot \phi),\kappa_{n,r}^m([g,h])(\phi)\big\rangle e'_ae'_b\\
\hphantom{\Phi^\varphi(f'_kf'_l)}{} = \sum_{a<b} \big\langle e_ae_b\cdot \kappa_{r*}^m(f_kf_l)\cdot \phi,\phi\big\rangle e'_ae'_b\\
\hphantom{\Phi^\varphi(f'_kf'_l)}{} = \sum_{a<b} \Phi^\phi(f_kf_l)(e_a,e_b) e'_ae'_b,
\end{gather*}
which means that the matrix representing $\Phi^\varphi(f'_kf'_l)$ with respect to the basis $(e'_1,\dots,e'_n)$ has the same coefficients as the matrix representing $\Phi^\phi(f_kf_l)$ with respect to the basis $(e_1,\dots,e_n)$. Hence,
\begin{gather*}\big[\Phi^\varphi(f'_kf'_l)\big]^2=-{\rm Id}_{\mathbb{R}^n}.\tag*{\qed}\end{gather*}\renewcommand{\qed}{}
\end{proof}

\begin{Lemma}\label{lemma: finding stabilizer 1} Let $\phi\in\Delta_n\otimes\Delta_r^{\otimes m}$ be a twisted pure spinor. If $\xi\in\mathfrak{spin}(n)$ is such that
\begin{gather*}\xi\cdot \phi=0,\end{gather*}
then
\begin{gather*}\big[\hat{\xi},\hat{\eta}_{kl}\big]=0\quad\quad\mbox{and}\quad\quad \tr\big(\hat{\xi}\hat{\eta}_{kl}^\phi\big)=0\end{gather*}
for all $1\leq k<l \leq r$.
\end{Lemma}
\begin{proof}
Let
\begin{gather*}\xi=\sum_{1\leq a<b\leq n} \xi(e_a,e_b)e_ae_b.\end{gather*}
Note that, if $1\leq s<t\leq n$ and $1\leq k<l\leq r $,
\begin{gather*}
 0 = \operatorname{Re}\big\langle e_se_t\cdot\kappa_{r^*}^m(f_kf_l)\cdot\xi\cdot\phi,\phi\big\rangle = \operatorname{Re}\big\langle e_se_t\cdot\xi\cdot\kappa_{r^*}^m(f_kf_l)\cdot\phi,\phi\big\rangle \\
 \hphantom{0}{} =
 \operatorname{Re}\bigg\langle e_se_t\cdot\left(\sum_{1\leq a<b\leq n} \xi(e_a,e_b)e_ae_b\right)\cdot\kappa_{r^*}^m(f_kf_l)\cdot\phi,\phi\bigg\rangle \\
\hphantom{0}{} =
 \sum_{1\leq a<b\leq n} \xi(e_a,e_b)\big\langle e_se_t\cdot e_ae_b\cdot\kappa_{r^*}^m(f_kf_l)\cdot\phi,\phi\big\rangle \\
\hphantom{0}{} =
 -\sum_{b} \big(\hat\eta_{kl}^\phi\big)_{tb}\hat\xi_{bs}
 +\sum_{b} \hat\xi_{tb}\big(\hat\eta_{kl}^\phi\big)_{bs}
=
 \big[\hat{\xi},\hat{\eta}_{kl}^\phi\big]_{ts},
\end{gather*}
i.e.,
\begin{gather*}\big[\hat{\xi},\hat{\eta}_{kl}^\phi\big]=0.\end{gather*}

Since $r\geq 3$, for $q\not= k,l$,
\begin{gather*}
 \tr\big(\hat{\xi}\hat{\eta}_{kl}^\phi\big) =
 \tr\big({-}\hat{\eta}_{kq}^\phi\hat{\xi}\hat{\eta}_{kl}^\phi\hat{\eta}_{kq}^\phi\big)
 =
 -\tr\big(\hat{\eta}_{kq}^\phi\hat{\xi}\hat{\eta}_{lq}^\phi\big)
 =
 -\tr\big(\hat{\xi}\hat{\eta}_{kq}^\phi\hat{\eta}_{lq}^\phi\big)
 =
 -\tr\big(\hat{\xi}\hat{\eta}_{kl}^\phi\big).\tag*{\qed}
\end{gather*}\renewcommand{\qed}{}
\end{proof}

\begin{Lemma}\label{lemma:stabilizer} Let $\phi\in\Delta_n\otimes\Delta_r^{\otimes m}$ be a twisted pure spinor. The annihilator subalgebra of $\phi$ in $\mathfrak{spin}(n)\oplus \mathfrak{spin}(r)$ is contained in one of the subalgebras $\tilde{\mathfrak{g}}$ listed in Table~{\rm \ref{table: annihilating algebras}},
where
\begin{gather*}\widetilde{\mathfrak{spin}}(r)=\operatorname{span}\big(\eta_{kl}^\phi+2f_kf_l\big)\subset \mathfrak{spin}(n)\oplus \mathfrak{spin}(r).\end{gather*}
\end{Lemma}

\begin{table}[h!]\centering
$
\begin{array}{|c|c|}
 \hline
 r \mbox{\ {\rm (mod 8)} } & \tilde{\mathfrak{g}}
 \rule{0pt}{3ex}
 \\
 \hline
 0 & \mathfrak{so}(m_1)\oplus\mathfrak{so}(m_2)\oplus \widetilde{\mathfrak{spin}}(r)\rule{0pt}{3ex}\\
 \hline
 1,7 & \mathfrak{so}(m)\oplus \widetilde{\mathfrak{spin}}(r)\rule{0pt}{3ex}\\
\hline
 2,6 & \mathfrak{u}(m)\oplus \widetilde{\mathfrak{spin}}(r)\rule{0pt}{3ex}\\
\hline
 3,5 & \mathfrak{sp}(m)\oplus \widetilde{\mathfrak{spin}}(r)\rule{0pt}{3ex}\\
\hline
 4 & \mathfrak{sp}(m_1)\oplus \mathfrak{sp}(m_2)\oplus \widetilde{\mathfrak{spin}}(r)\rule{0pt}{3ex}\\
\hline
 \end{array}
$
\caption{}\label{table: annihilating algebras}
\end{table}

\begin{proof} Let
\begin{gather*}
 \xi \in \mathfrak{spin}(n),\qquad
 \sigma = \sum_{1\leq k<l\leq r} \sigma_{kl}f_kf_l \in\mathfrak{spin}(r),
\end{gather*}
be such that
\begin{gather*}\big(\xi+\kappa_{r^*}^{m} (\sigma)\big)\cdot \phi =0.\end{gather*}
Expanding this identity,
\begin{gather*}
 0 = \big(\xi+\kappa_{r^*}^{m} (\sigma)\big)\cdot \phi =
 \xi\cdot\phi + \sum_{1\leq k<l\leq r} \sigma_{kl}\kappa_{r^*}^{m}(f_kf_l)\cdot\phi\\
 \hphantom{0}{} =
 \xi\cdot\phi + {1\over 2}\sum_{1\leq k<l\leq r} \sigma_{kl}2\kappa_{r^*}^{m}(f_kf_l)\cdot\phi =
 \xi\cdot\phi - {1\over 2}\sum_{1\leq k<l\leq r} \sigma_{kl}\eta_{kl}^\phi\cdot\phi\\
\hphantom{0}{} =
 \left(\xi - {1\over 2}\sum_{1\leq k<l\leq r} \sigma_{kl}\eta_{kl}^\phi\right)\cdot\phi.
\end{gather*}
By Lemma \ref{lemma: finding stabilizer 1},
\begin{gather*}\hat{\xi} - {1\over 2}\sum_{1\leq k<l\leq r} \sigma_{kl}\hat{\eta}_{kl}^\phi\in C_{\mathfrak{so}(n)}(\mathfrak{spin}(r))\end{gather*}
and is orthogonal to $\operatorname{span}(\hat{\eta}_{kl})$.
Thus,
\begin{gather*}\xi+\sigma = \left(\xi - {1\over 2}\sum_{1\leq k<l\leq r} \sigma_{kl}\eta_{kl}^\phi\right) + \left({1\over 2}\sum_{1\leq k<l\leq r} \sigma_{kl}(\eta_{kl}^\phi + 2\,f_kf_l)\right) \in \mathfrak{spin}(n)\oplus\mathfrak{spin}(r).\end{gather*}
The table follows from \cite[Theorems 3.1 and 3.2]{AH}.
\end{proof}

We refer the reader to \cite[Theorem 3.2]{AGH} for a description of the (connected components of the identity element of the) corresponding Lie groups.

\begin{Lemma}\label{lemma: Spin(r) orbit pure spinor}Let $\phi\in\Delta_n\otimes\Delta_r^{\otimes m}$ be a twisted pure spinor. Every element in the orbit $\widehat{{\rm Spin}(r)}\cdot \phi$ induces the same linear even-Clifford Hermitian structure of rank~$r$ on $\mathbb{R}^n$, where $\widehat{{\rm Spin}(r)}$ denotes the canonical copy of ${\rm Spin}(r)$ in ${\rm Spin}^r(n)$ given by the elements $[(1,g)]$.
\end{Lemma}
\begin{proof} Let $g\in \widehat{{\rm Spin}(r)}\subset {\rm Spin}^r(n)$, i.e., $\lambda_n^r(g) = (1,g_2)\in {\rm SO}(n)\times {\rm SO}(r)$. Then
\begin{gather*}
\hat\Phi^{g(\phi)}(f_kf_l)(X) = \sum_{b=1}^n \big\langle X\wedge e_b\cdot \kappa_{r^*}^m(f_kf_l)\cdot g(\phi),g(\phi)\big\rangle e_b \\
\hphantom{\hat\Phi^{g(\phi)}(f_kf_l)(X)}{} =
 \sum_{b=1}^n\big\langle X\wedge e_b\cdot \kappa_{r^*}^m(g_2(f_k')g_2(f_l'))\cdot g(\phi),g(\phi)\big\rangle e_b\\
\hphantom{\hat\Phi^{g(\phi)}(f_kf_l)(X)}{} =
 \sum_{b=1}^n\big\langle g(X\wedge e_b\cdot \kappa_{r^*}^m(f_k'f_l')\cdot \phi),g(\phi)\big\rangle e_b\\
\hphantom{\hat\Phi^{g(\phi)}(f_kf_l)(X)}{} =
 \sum_{b=1}^n\big\langle X\wedge e_b\cdot \kappa_{r^*}^m(f_k'f_l')\cdot \phi,\phi\big\rangle e_b =
\hat\Phi^\phi(f_k'f_l')(X),
\end{gather*}
where $g_2(f_k')=f_k$.
\end{proof}

\begin{Lemma}\label{lemma: producing more pure spinors} Let $\phi\in\Delta_n\otimes\Delta_r^{\otimes m}$ be a twisted pure spinor. Let $X\in\mathbb{R}^n$ such that $|X|=1$. Then, $X\cdot\phi$ is also a twisted pure spinor.
\end{Lemma}
\begin{proof}\allowdisplaybreaks Without loss of generality we can assume that $X=e_1$. First, notice that
\begin{gather*}
\eta_{kl}^{e_1\cdot\phi}
 = \sum_{a<b} \big\langle e_ae_b\cdot\kappa_{r*}^m(f_{kl})\cdot e_1\cdot\phi,e_1\cdot\phi\big\rangle e_ae_b\\
 \hphantom{\eta_{kl}^{e_1\cdot\phi}}{} = \sum_{1<b} \big\langle e_1e_b\cdot\kappa_{r*}^m(f_{kl})\cdot e_1\cdot\phi,e_1\cdot\phi\big\rangle e_1e_b
 +\sum_{2\leq a<b} \big\langle e_ae_b\cdot\kappa_{r*}^m(f_{kl})\cdot e_1\cdot\phi,e_1\cdot\phi\big\rangle e_ae_b\\
 \hphantom{\eta_{kl}^{e_1\cdot\phi}}{}= -\sum_{1<b} \big\langle e_1e_1e_b\cdot\kappa_{r*}^m(f_{kl})\cdot\phi,e_1\cdot\phi\big\rangle e_1e_b
 +\sum_{2\leq a<b} \big\langle e_1e_ae_b\cdot\kappa_{r*}^m(f_{kl})\cdot\phi,e_1\cdot\phi\big\rangle e_ae_b\\
 \hphantom{\eta_{kl}^{e_1\cdot\phi}}{} = -\sum_{1<b} \big\langle e_1e_b\cdot\kappa_{r*}^m(f_{kl})\cdot\phi,\cdot\phi\big\rangle e_1e_b
 +\sum_{2\leq a<b} \big\langle e_ae_b\cdot\kappa_{r*}^m(f_{kl})\cdot\phi,\cdot\phi\big\rangle e_ae_b
\end{gather*}
has the same coefficients as those in $\eta_{kl}^{e_1\cdot\phi}$ but with some signs changed. More precisely
\begin{alignat*}{3}
&\eta_{kl}^{e_1\cdot\phi}(e_1,e_b) = -\eta_{kl}^{\phi}(e_1,e_b)\qquad&&\text{for} \quad 1\leq b,&\\
& \eta_{kl}^{e_1\cdot\phi}(e_a,e_b) = \eta_{kl}^{\phi}(e_a,e_b)\qquad&& \text{for} \quad 2\leq a<b.&
\end{alignat*}
Thus
\begin{gather*}
\big(\eta_{kl}^{e_1\cdot\phi}+2\kappa_{r*}^m(f_{kl})\big)\cdot (e_1\cdot\phi)
 = -\left(\sum_{1<b} \big\langle e_1e_b\cdot\kappa_{r*}^m(f_{kl})\cdot\phi,\cdot\phi\big\rangle e_1e_b\right)\cdot e_1\cdot\phi\\
 \qquad{} +\left(\sum_{2\leq a<b} \big\langle e_ae_b\cdot\kappa_{r*}^m(f_{kl})\cdot\phi,\cdot\phi\big\rangle e_ae_b\right)\cdot e_1\cdot\phi
 + 2\kappa_{r*}^m(f_{kl})\cdot (e_1\cdot\phi)\\
 \quad{} = \left(\sum_{1<b} \big\langle e_1e_b\cdot\kappa_{r*}^m(f_{kl})\cdot\phi,\cdot\phi\big\rangle e_1e_1e_b\right)\cdot\phi\\
 \qquad{} +\left(\sum_{2\leq a<b} \big\langle e_ae_b\cdot\kappa_{r*}^m(f_{kl})\cdot\phi,\cdot\phi\big\rangle e_1e_ae_b\right)\cdot\phi
 + 2e_1\cdot\kappa_{r*}^m(f_{kl})\cdot\phi\\
 \quad{} = e_1\cdot\left(\sum_{1<b} \big\langle e_1e_b\cdot\kappa_{r*}^m(f_{kl})\cdot\phi,\cdot\phi\big\rangle e_1e_b\right)\cdot\phi\\
 \qquad{} +e_1\cdot\left(\sum_{2\leq a<b} \big\langle e_ae_b\cdot\kappa_{r*}^m(f_{kl})\cdot\phi,\cdot\phi\big\rangle e_ae_b\right)\cdot\phi
 + e_1\cdot 2\kappa_{r*}^m(f_{kl})\cdot\phi\\
 \quad{}= e_1\cdot\left(\sum_{a<b} \big\langle e_ae_b\cdot\kappa_{r*}^m(f_{kl})\cdot\phi,\cdot\phi\big\rangle e_ae_b\right)\cdot\phi
 +e_1\cdot 2\kappa_{r*}^m(f_{kl})\cdot\phi\\
\quad{} = e_1\cdot\big(\eta_{kl}^{\phi} +2\kappa_{r*}^m(f_{kl})\big)\cdot\phi = 0.
\end{gather*}
Regarding the endomorphism $\hat\eta_{kl}^{e_1\cdot\phi}$ we have
\begin{gather*}
\big(\big(\hat\eta_{kl}^{e_1\cdot\phi}\big)^2\big)_{ts} =\sum_{b}\big(\hat\eta_{kl}^{e_1\cdot\phi}\big)_{tb}\big(\hat\eta_{kl}^{e_1\cdot\phi}\big)_{bs}\\
\hphantom{\big(\big(\hat\eta_{kl}^{e_1\cdot\phi}\big)^2\big)_{ts}}{}
=\big(\hat\eta_{kl}^{e_1\cdot\phi}\big)_{t1}\big(\hat\eta_{kl}^{e_1\cdot\phi}\big)_{1s}
 +\sum_{b\geq2}\big(\hat\eta_{kl}^{e_1\cdot\phi}\big)_{tb}\big(\hat\eta_{kl}^{e_1\cdot\phi}\big)_{bs}\\
\hphantom{\big(\big(\hat\eta_{kl}^{e_1\cdot\phi}\big)^2\big)_{ts}}{} =\eta_{kl}^{e_1\cdot\phi}(e_1,e_t)\eta_{kl}^{e_1\cdot\phi}(e_s,e_1)
 +\sum_{b\geq2}\eta_{kl}^{e_1\cdot\phi}(e_b,e_t)\eta_{kl}^{e_1\cdot\phi}(e_s,e_b)\\
\hphantom{\big(\big(\hat\eta_{kl}^{e_1\cdot\phi}\big)^2\big)_{ts}}{} =\big({-}\eta_{kl}^{\phi}(e_1,e_t)\big)\big({-}\eta_{kl}^{\phi}(e_s,e_1)\big)
 +\sum_{b\geq2}\eta_{kl}^{\phi}(e_b,e_t)\eta_{kl}^{\phi}(e_s,e_b)\\
\hphantom{\big(\big(\hat\eta_{kl}^{e_1\cdot\phi}\big)^2\big)_{ts}}{} =\sum_{b}\eta_{kl}^{\phi}(e_b,e_t)\eta_{kl}^{\phi}(e_s,e_b) =-\delta_{ts},
\end{gather*}
i.e.,
\begin{gather*}\big(\hat\eta_{kl}^{e_1\cdot\phi}\big)^2=-{\rm Id}_{\mathbb{R}^n}.\tag*{\qed}\end{gather*}\renewcommand{\qed}{}
\end{proof}

\subsection[Pure spinors: $r=2$]{Pure spinors: $\boldsymbol{r=2}$}

We have left out of our discussion the case $r=2$ due to the following two reasons:
\begin{enumerate}\itemsep=0pt
\item The prototypical pure ${\rm Spin}^c$ spinor is given by $\varphi=u_{1,\dots,1}\in\Delta_{2n}$. It satisfies the equation
\begin{gather*}e_{2j-1}\cdot \varphi = \sqrt{-1} e_{2j}\cdot\varphi\end{gather*}
for $1\leq j\leq n$.
This means that, as described in the introduction, the orthogonal complex structure determined by $\varphi$ is the standard one on $\mathbb{R}^{2n}$,
\begin{gather*}J_0=\left(
\begin{matrix}
 & -1 & & & \\
1 & & & & \\
 & & \ddots & & \\
 & & & & -1\\
 & & & 1 &
\end{matrix}
\right).\end{gather*}
Furthermore,
\begin{gather*} e_{2j-1}e_{2j}\cdot\varphi = \sqrt{-1} \varphi ,\end{gather*}
so that the associated real 2-form
\begin{gather*}
\eta^\varphi := \sum_{1\leq a<b\leq 2n} \sqrt{-1}\langle e_a e_b \cdot \varphi,\varphi\rangle e_ae_b= -\sum_{a=1}^n e_{2a-1}e_{2a}
\end{gather*}
is such that
\begin{gather*}
 \eta^\varphi\cdot\varphi = -\sum_{a=1}^n e_{2a-1}e_{2a}\cdot\varphi= -\sum_{a=1}^n {\rm i}\varphi= - n \sqrt{-1} \varphi,
\end{gather*}
i.e., the associated 2-form $\eta^\varphi$ and the spinor $\varphi$ satisfy
\begin{gather*}\big(\eta^\varphi + n \sqrt{-1}\big)\cdot \varphi =0,\end{gather*}
which has the coefficient $n$ instead of $2$, and describes the $\big(n\sqrt{-1}\big)$-eigenspace of the corresponding K\"ahler form (see~\cite{Friedrich}).

\item Recall that ${\rm Spin}(2)$ is very different from all other spin groups ${\rm Spin}(r)$, $r\geq 3$, since it is abelian, non-simple and non-simply-connected. All of these differences are somehow reflected by the fact that there are no pure ${\rm Spin}^2(2n)$-spinors according to Definition~\ref{def:pure-spinor}, except for $n=2$. Instead, there are spinors satisfying the equations
\begin{gather}
 \big(\eta_{12}^\phi + n \kappa_2^1(f_{12})\big) \cdot\phi = 0,\label{eq: twisted pure spinor rank 2 eq 1}\\
 \big(\hat\eta_{12}^\phi\big)^2 = -{\rm Id}_{\mathbb{R}^{2n}},\label{eq: twisted pure spinor rank 2 eq 2}
\end{gather}
with coefficient $n$ instead of $2$, just as in the ${\rm Spin}^c$ description above. However, in this rank, we only need the twisted pure spinor to induce a complex structure, which is dictated by~\eqref{eq: twisted pure spinor rank 2 eq 2}. Thus~\eqref{eq: twisted pure spinor rank 2 eq 1} becomes redundant.
\end{enumerate}

\subsection{Existence of pure spinors}\label{subsec: existence pure spinors}

In this subsection we present explicit pure spinors for the ranks $r= 3,7$.
Let us define the following maps:
\begin{align*}
 G\colon \ \{\pm 1\}^{\times m} &\lra \{\pm 1\}^{\times 2m},\\
 (\varepsilon_1,\dots,\varepsilon_m) &\longmapsto (\varepsilon_1,\varepsilon_1,\dots,\varepsilon_m,\varepsilon_m),\\
 H\colon \ \{\pm 1\}^{\times m} &\lra \{0,1,\dots,m\},\\
 (\varepsilon_1,\dots,\varepsilon_m) &\longmapsto \sum_{j=1}^m {1-\varepsilon_j\over 2}.
\end{align*}
Define
\begin{gather*}\{\pm1\}_j^{\times m} := H^{-1}(j),\end{gather*}
which is the set of elements in $\{\pm1\}^{\times m}$ with exactly $j$ entries equal to $(-1)$. Note that
\begin{gather*}\big|\{\pm1\}_j^{\times m}\big|={m\choose j}.\end{gather*}

\subsubsection[Dimension $n=4m$, rank $r=3$]{Dimension $\boldsymbol{n=4m}$, rank $\boldsymbol{r=3}$}\label{subsec:QK-pure-spinor}

Consider the following spinors
\begin{gather*}
\psi_j = \sum_{(\varepsilon_1,\dots,\varepsilon_m)\in\{\pm1\}_j^{\times m}} u_{G(\varepsilon_1,\dots,\varepsilon_m)}, \qquad
\varphi_j =\sum_{(\varepsilon_1,\dots,\varepsilon_m)\in\{\pm1\}_j^{\times m}}v_{(\varepsilon_1,\dots,\varepsilon_m)}.
\end{gather*}
The twisted spinor
\begin{gather}
\phi=\sqrt{3\over m+2}{1\over \sqrt{m+1}}\sum_{j=0}^m {1\over {m\choose j}} \psi_j\otimes\varphi_{m-j} \in \Delta_{4m}\otimes\Delta_3^{\otimes m} \label{eq: rank 3 pure spinor}
\end{gather}
is pure. The 2-forms associated to $\phi$ are
\begin{gather*}
 \eta_{12}^\phi = \sum_{j=1}^{m} (e_{4j-3}e_{4j-2}+e_{4j-1}e_{4j}),\qquad
 \eta_{13}^\phi = \sum_{j=1}^{m} (-e_{4j-3}e_{4j-1}+e_{4j-2}e_{4j}),\\
 \eta_{23}^\phi = \sum_{j=1}^{m} (-e_{4j-3}e_{4j}-e_{4j-2}e_{4j-1}),
\end{gather*}
which span a copy of $\mathfrak{spin}(3)\in \mathfrak{so}(4m)$. For instance, let us compute
\begin{gather*}
\eta_{13}^{\phi}(e_r,e_s) = \operatorname{Re}\big\langle e_r \wedge e_s \cdot \kappa_{3\ast}^{m}(f_{13})\cdot
\phi,\phi\big\rangle \\
{} = \frac{3}{(m+2)(m+1)}\operatorname{Re}\bigg \langle \sum_{j=0}^{m}\frac{1}{{m\choose j}}
\big(e_re_s\cdot\psi_{j}\otimes\kappa_{3\ast}^{m}(f_{13})\cdot\varphi_{m-j}\big), \sum_{j=0}^{m}\frac{1}{{m\choose
j}}\psi_{j}\otimes\varphi_{m-j}\bigg\rangle .
\end{gather*}
Consider
\begin{gather*}
e_re_s\cdot\psi_{j} = e_{4k_1-j_1}e_{4k_2-j_2}\cdot\psi_{j}= \sum_{(\varepsilon_{1},\dots, \varepsilon_{m})\in\{\pm
1\}_{j}^{m}}e_{4k_1-j_1}e_{4k_2-j_2}\cdot u_{G(\varepsilon_1,\dots,\varepsilon_m)},
\end{gather*}
where $4k_1-j_1< 4k_2-j_2$, $0\leq j_1$, $j_2\leq 4$. Define $\big(\varepsilon_1^{0},\varepsilon_1^{1},\dots,
\varepsilon_m^{0},\varepsilon_m^{1} \big):=(\varepsilon_{1},\varepsilon_1 ,\dots, \varepsilon_{m},\varepsilon_m )$. Note that
\begin{align*}
e_{4k-j}\cdot u_{G(\varepsilon_{1},\dots,\varepsilon_m)} &= e_{4k-j}\cdot u_{(\varepsilon_{1}^0,\varepsilon_{1}^{1},\dots,\varepsilon_{m}^{0},\varepsilon_{m}^{1})}\\
&= -(i)^j(\varepsilon_{m-k+1})^{[\frac{j+1}{2}]} u_{(\varepsilon_{1}^0,\varepsilon_{1}^1,\dots,(-\varepsilon_{m-k+1})^{[\frac{j}{2}]},\dots,\varepsilon_{m}^0, \varepsilon_{m}^1)}.
\end{align*}
Thus,
\begin{gather*}
 e_re_s\cdot\psi_{j}=
\sum_{(\varepsilon_{1},\dots ,\varepsilon_{m})\in\{\pm
1\}_{j}^{m}}(i)^{j_1+j_2}(\varepsilon_{m-k_1+1})^{[\frac{j_1+1}{2}]}(\varepsilon_{m-k_2+1})^{[\frac{j_2+1}{2}]}\\
\hphantom{e_re_s\cdot\psi_{j}=}{} \times u_{(\varepsilon_{1}^0, \varepsilon_{1}^1 ,\dots, (-\varepsilon_{m-k_1+1})^{[\frac{j_1}{2}]}, \dots,
(-\varepsilon_{m-k_2+1})^{[\frac{j_2}{2}]}, \dots,\varepsilon_{m}^0,\varepsilon_{m}^1)}.
\end{gather*}
On the other hand,
\begin{gather*}
\kappa_{3\ast}^{m}(f_{13})\cdot \varphi_{m-j} = \sum_{(\varepsilon_{1},\dots ,\varepsilon_{m})\in\{\pm1\}_{m-j}^{m}}\kappa_{3\ast}^{m}(f_{13})\cdot v_{(\varepsilon_{1},\dots,\varepsilon_m)}\\
\hphantom{\kappa_{3\ast}^{m}(f_{13})\cdot \varphi_{m-j}}{} =
 \sum_{(\varepsilon_{1},\dots,
\varepsilon_{m})\in\{\pm 1\}_{m-j}^{m}}\left( \sum_{l=1}^{m}\varepsilon_{l}\right)
v_{(\varepsilon_{1},\dots,-\varepsilon_{l},\dots,\varepsilon_m)}\\
\hphantom{\kappa_{3\ast}^{m}(f_{13})\cdot \varphi_{m-j}}{}
= \sum_{(\varepsilon_{1},\dots ,\varepsilon_{m})\in\{\pm
1\}_{m-j+1}^{m}}(m-(j-1))v_{(\varepsilon_{1},\dots,\varepsilon_m)}\\
\hphantom{\kappa_{3\ast}^{m}(f_{13})\cdot \varphi_{m-j}=}{}
-\sum_{(\varepsilon_{1},\dots, \varepsilon_{m})\in\{\pm 1\}_{m-j-1}^{m}}(j+1)v_{(\varepsilon_{1},\dots,\varepsilon_m)}.
\end{gather*}
For $k_1<k_2$
\begin{gather*}
\big\langle u_{(\varepsilon_{1}^0, \varepsilon_{1}^1, \dots, (-\varepsilon_{m-k_1+1})^{[\frac{j_1}{2}]} ,\dots,
(-\varepsilon_{m-k_2+1})^{[\frac{j_2}{2}]}, \dots, \varepsilon_{m}^0,\varepsilon_{m}^1)} \otimes v_{(\varepsilon_{1},\dots,\varepsilon_m)},\\
\qquad{} u_{(\tilde{\varepsilon}_{1}^0, \tilde{\varepsilon}_{1}^1, \dots ,\tilde{\varepsilon}_{m}^0,\tilde{\varepsilon}_{m}^1)} \otimes
v_{(\tilde{\varepsilon}_{1},\dots,\tilde{\varepsilon}_m)}\big\rangle =0.
\end{gather*}
Thus, for $k_1<k_2$
\begin{gather*}\eta_{13}^\phi(e_{4k_1-j_1},e_{4k_2-j_2})=0.\end{gather*}
Now consider $k_1=k_2=k$. In this case
\begin{gather*}
\eta_{13}^\phi(e_{4k-j_1},e_{4k-j_2})=
 \frac{3}{(m+2)(m+1)}\operatorname{Re}\bigg\langle \sum_{j=0}^{m}\frac{1}{{m\choose j}}\\
 \hphantom{\eta_{13}^\phi(e_{4k-j_1},e_{4k-j_2})=}{}\times
 \sum_{(\varepsilon_{1},\dots, \varepsilon_{m})\in\{\pm
1\}_{j}^{m}}(i)^{j_1+j_2}(\varepsilon_{m-k+1})^{[\frac{j_1+1}{2}]}(\varepsilon_{m-k+1})^{[\frac{j_2+1}{2}]} \\
 \hphantom{\eta_{13}^\phi(e_{4k-j_1},e_{4k-j_2})=}{}\times
 u_{(\varepsilon_{1}^0, \varepsilon_{1}^1, \dots, (-\varepsilon_{m-k+1})^{[\frac{j_2}{2}]} ,\dots, (-\varepsilon_{m-k+1})^{[\frac{j_1}{2}]},\dots, \varepsilon_{m}^0,\varepsilon_{m}^1)}\\
\hphantom{\eta_{13}^\phi(e_{4k-j_1},e_{4k-j_2})=}{} \otimes
 \bigg[ \sum_{(\varepsilon_{1},\dots, \varepsilon_{m})\in\{\pm
1\}_{m-j+1}^{m}}(m-(j-1))v_{(\varepsilon_{1},\dots,\varepsilon_m)}\\
\hphantom{\eta_{13}^\phi(e_{4k-j_1},e_{4k-j_2})=}{}
-\sum_{(\varepsilon_{1},\dots, \varepsilon_{m})\in\{\pm
1\}_{m-j-1}^{m}}(j+1)v_{(\varepsilon_{1},\dots,\varepsilon_m)}\bigg],\\
\hphantom{\eta_{13}^\phi(e_{4k-j_1},e_{4k-j_2})=}{}
 \sum_{j=0}^{m}\frac{1}{{m\choose j}}\sum_{(\varepsilon_{1},\dots, \varepsilon_{m})\in\{\pm
1\}_{j}^{m}}u_{G(\varepsilon_1,\dots,\varepsilon_m)}\otimes
\sum_{(\varepsilon_{1}\dots \varepsilon_{m})\in\{\pm 1\}_{m-j}^{m}}v_{(\varepsilon_1,\dots,\varepsilon_m)}\bigg\rangle.
\end{gather*}
We have the following cases:
\begin{enumerate}\itemsep=0pt\allowdisplaybreaks
\item If $[\frac{j_1}{2}]=[\frac{j_2}{2}]$ so that $j_1+j_2=1$ (mod 4):
\begin{gather*}
 \eta_{13}^\phi(e_{4k-j_1},e_{4k-j_2}) = \frac{3}{(m+2)(m+1)}\\
\hphantom{\eta_{13}^\phi(e_{4k-j_1},e_{4k-j_2}) =}{}
\times \operatorname{Re}\bigg\{ i\bigg\langle \sum_{j=0}^{m}\frac{1}{{m\choose j}}\sum_{(\varepsilon_{1},\dots,
\varepsilon_{m})\in\{\pm 1\}_{j}^{m}}\varepsilon_{m-k+1} u_{G(\varepsilon_{1}, \dots, \varepsilon_{m})}\\
\hphantom{\eta_{13}^\phi(e_{4k-j_1},e_{4k-j_2}) =}{} \otimes \bigg[ \sum_{(\varepsilon_{1},\dots ,\varepsilon_{m})\in\{\pm
1\}_{m-j+1}^{m}}(m-(j-1))v_{(\varepsilon_{1},\dots,\varepsilon_m)}\\
\hphantom{\eta_{13}^\phi(e_{4k-j_1},e_{4k-j_2}) =}{}
-\sum_{(\varepsilon_{1},\dots, \varepsilon_{m})\in\{\pm
1\}_{m-j-1}^{m}}(j+1)v_{(\varepsilon_{1},\dots,\varepsilon_m)}\bigg],\\
\hphantom{\eta_{13}^\phi(e_{4k-j_1},e_{4k-j_2}) =}{}
 \sum_{j=0}^{m}\frac{1}{{m\choose j}}\sum_{(\varepsilon_{1},\dots, \varepsilon_{m})\in\{\pm
1\}_{j}^{m}}u_{G(\varepsilon_1\dots\varepsilon_m)}\\
\hphantom{\eta_{13}^\phi(e_{4k-j_1},e_{4k-j_2}) =}{}
\otimes
\sum_{(\varepsilon_{1},\dots, \varepsilon_{m})\in\{\pm 1\}_{m-j}^{m}}v_{(\varepsilon_1,\dots,\varepsilon_m)}\bigg\rangle\bigg\} =0.
\end{gather*}
\item If $[\frac{j_1}{2}]\neq[\frac{j_2}{2}]$:
\begin{itemize}\allowdisplaybreaks\itemsep=0pt
\item[i)] $j_1=2$, $j_2=0$ and $j_1+j_2=2$
\begin{gather*}
 \eta_{13}^\phi(e_{4k-j_1},e_{4k-j_2}) =
 \frac{-3}{(m+2)(m+1)}\bigg\langle \sum_{j=0}^{m}\frac{1}{{m\choose j}}\\
\hphantom{\eta_{13}^\phi(e_{4k-j_1},e_{4k-j_2}) =}{}\times
\bigg[{} -\sum_{(\varepsilon_{1},\dots,\hat{\varepsilon}_{m-k+1},\dots, \varepsilon_{m})\in\{\pm
1\}_{j-1}^{m-1}}u_{G(\varepsilon_{1},\dots, \varepsilon_{m-k} 1,\dots,\varepsilon_{m})}\\
\hphantom{\eta_{13}^\phi(e_{4k-j_1},e_{4k-j_2}) =}{} + \sum_{(\varepsilon_{1},\dots,\hat{\varepsilon}_{m-k+1},\dots, \varepsilon_{m})\in\{\pm 1\}_{j}^{m-1}}u_{G(\varepsilon_{1},\dots,
\varepsilon_{m-k -1}, \dots,\varepsilon_{m})} \bigg]\\
\hphantom{\eta_{13}^\phi(e_{4k-j_1},e_{4k-j_2}) =}{}
 \otimes\bigg[ \sum_{(\varepsilon_{1},\dots, \varepsilon_{m})\in\{\pm
1\}_{m-j+1}^{m}}(m-(j-1))v_{(\varepsilon_{1},\dots,\varepsilon_m)}\\
\hphantom{\eta_{13}^\phi(e_{4k-j_1},e_{4k-j_2}) =}{}
-\sum_{(\varepsilon_{1},\dots, \varepsilon_{m})\in\{\pm
1\}_{m-j-1}^{m}}(j+1)v_{(\varepsilon_{1},\dots,\varepsilon_m)}\bigg] ,\\
\hphantom{\eta_{13}^\phi(e_{4k-j_1},e_{4k-j_2}) =}{}
\sum_{j=0}^{m}\frac{1}{{m\choose j}}\sum_{(\varepsilon_{1},\dots, \varepsilon_{m})\in\{\pm
1\}_{j}^{m}}u_{G(\varepsilon_1,\dots,\varepsilon_m)}\\
\hphantom{\eta_{13}^\phi(e_{4k-j_1},e_{4k-j_2}) =}{}
\otimes
\sum_{(\varepsilon_{1},\dots, \varepsilon_{m})\in\{\pm 1\}_{m-j}^{m}}v_{(\varepsilon_1,\dots,\varepsilon_m)}\bigg\rangle\\
\hphantom{\eta_{13}^\phi(e_{4k-j_1},e_{4k-j_2})}{}
 = \frac{-3}{(m+2)(m+1)}\bigg\rangle \sum_{j=0}^{m}\frac{1}{{m\choose j}}\\
\hphantom{\eta_{13}^\phi(e_{4k-j_1},e_{4k-j_2})=}{} \times
 \bigg\{ {-}\bigg[ \sum_{(\varepsilon_{1},\dots,\hat{\varepsilon}_{m-k+1},\dots, \varepsilon_{m})\in\{\pm
1\}_{j-1}^{m-1}}u_{G(\varepsilon_{1},\dots, \varepsilon_{m-k} 1,\dots,\varepsilon_{m})}\\
\hphantom{\eta_{13}^\phi(e_{4k-j_1},e_{4k-j_2})=}{}
\otimes \sum_{(\varepsilon_{1},\dots,
\varepsilon_{m})\in\{\pm 1\}_{m-j+1}^{m}}(m-(j-1))v_{(\varepsilon_{1},\dots,\varepsilon_m)}\bigg]\\
\hphantom{\eta_{13}^\phi(e_{4k-j_1},e_{4k-j_2})=}{}
 - \bigg[ \sum_{(\varepsilon_{1},\dots,\hat{\varepsilon}_{m-k+1},\dots, \varepsilon_{m})\in\{\pm
1\}_{j}^{m-1}}u_{G(\varepsilon_{1},\dots, \varepsilon_{m-k -1} ,\dots,\varepsilon_{m})}\\
\hphantom{\eta_{13}^\phi(e_{4k-j_1},e_{4k-j_2})=}{}
\otimes \sum_{(\varepsilon_{1},\dots, \varepsilon_{m})\in\{\pm
1\}_{m-j-1}^{m}}(j+1)v_{(\varepsilon_{1},\dots,\varepsilon_m)}\bigg] \bigg\},\\
\hphantom{\eta_{13}^\phi(e_{4k-j_1},e_{4k-j_2})=}{}
 \sum_{j=0}^{m}\frac{1}{{m\choose j}}\sum_{(\varepsilon_{1},\dots, \varepsilon_{m})\in\{\pm
1\}_{j}^{m}}u_{G(\varepsilon_1,\dots,\varepsilon_m)}\\
\hphantom{\eta_{13}^\phi(e_{4k-j_1},e_{4k-j_2})=}{}
\otimes \sum_{(\varepsilon_{1},\dots, \varepsilon_{m})\in\{\pm 1\}_{m-j}^{m}}v_{(\varepsilon_1,\dots,\varepsilon_m)}\bigg\rangle \\
\hphantom{\eta_{13}^\phi(e_{4k-j_1},e_{4k-j_2})}{}
 = \frac{-3}{(m+2)(m+1)}\sum_{j=0}^{m}\frac{1}{{m\choose j}}\bigg[ {-}{{{m-1}\choose{j-1}}}{{m\choose
{m-j+1}}}\\
\hphantom{\eta_{13}^\phi(e_{4k-j_1},e_{4k-j_2})=}{}\times
\frac{1}{{m\choose{j-1}}}{(m-j+1)} {{- {{m-1}\choose j}}{{m}\choose{m-j-1}}}
 \frac{1}{{m\choose
{j+1}}}(j+1)\bigg] \\
\hphantom{\eta_{13}^\phi(e_{4k-j_1},e_{4k-j_2})}{}
 = \frac{-3}{m(m+2)(m+1)}\sum_{j=0}^{m}\big(2j^2-2mj-m\big)=1.
\end{gather*}
\item[ii)] If $j_1=2$, $j_2=1$ and $j_1+j_2=3$
\begin{gather*}
 \eta_{13}^\phi(e_{4k-j_1},e_{4k-j_2}) =
 \frac{3}{(m+2)(m+1)}\operatorname{Re}\bigg\{ {-}{\rm i}\bigg\langle \sum_{j=0}^{m}\frac{1}{{m\choose
j}}\\
\hphantom{\eta_{13}^\phi(e_{4k-j_1},e_{4k-j_2}) =}{}\times
\sum_{(\varepsilon_{1},\dots ,\varepsilon_{m})\in\{\pm 1\}_{j}^{m}}u_{G(\varepsilon_{1},\dots,
-\varepsilon_{m-k+1},\dots,\varepsilon_{m})}\\
\hphantom{\eta_{13}^\phi(e_{4k-j_1},e_{4k-j_2}) =}{}
\otimes\bigg[ \sum_{(\varepsilon_{1},\dots, \varepsilon_{m})\in\{\pm
1\}_{m-j+1}^{m}}(m-(j-1))v_{(\varepsilon_{1},\dots,\varepsilon_m)}\\
\hphantom{\eta_{13}^\phi(e_{4k-j_1},e_{4k-j_2}) =}{}
-\sum_{(\varepsilon_{1},\dots, \varepsilon_{m})\in\{\pm
1\}_{m-j-1}^{m}}(j+1)v_{(\varepsilon_{1},\dots,\varepsilon_m)}\bigg] ,\\
\hphantom{\eta_{13}^\phi(e_{4k-j_1},e_{4k-j_2}) =}{}
 \sum_{j=0}^{m}\frac{1}{{m\choose j}}\sum_{(\varepsilon_{1},\dots, \varepsilon_{m})\in\{\pm
1\}_{j}^{m}}u_{G(\varepsilon_1,\dots,\varepsilon_m)}\\
\hphantom{\eta_{13}^\phi(e_{4k-j_1},e_{4k-j_2}) =}{}
\otimes
\sum_{(\varepsilon_{1},\dots ,\varepsilon_{m})\in\{\pm 1\}_{m-j}^{m}}v_{(\varepsilon_1,\dots,\varepsilon_m)}\bigg\rangle\bigg\} =0.
\end{gather*}
\item[iii)] If $j_1=3$, $j_2=0$ and $j_1+j_2=3$
\begin{gather*}
 \eta_{13}^\phi(e_{4k-j_1},e_{4k-j_2}) =0.
\end{gather*}
\item[iv)] If $j_1=3$, $j_2=1$ and $j_1+j_2=4$
\begin{gather*}
 \eta_{13}^\phi(e_{4k-j_1},e_{4k-j_2})=-1.
\end{gather*}
\end{itemize}
\end{enumerate}

The spinor $\phi$ is annihilated by
\begin{gather*}
 \eta_{12}^\phi + 2 \kappa_{3*}^m(f_{12}), \qquad
 \eta_{13}^\phi + 2 \kappa_{3*}^m(f_{13}), \qquad
 \eta_{23}^\phi + 2 \kappa_{3*}^m(f_{23}),
\end{gather*}
and the forms
\begin{gather*}
 \beta_{ij}^1 = e_{4i-3}e_{4j-3} + e_{4i-2}e_{4j-2} + e_{4i-1}e_{4j-1} + e_{4i}e_{4j}, \\
 \beta_{ij}^2 = e_{4i-3}e_{4j-2} - e_{4i-1}e_{4j}, \\
 \beta_{ij}^3 = e_{4i-3}e_{4j-1} + e_{4i-2}e_{4j}, \\
 \beta_{ij}^4 = e_{4i-3}e_{4j} - e_{4i-2}e_{4j-1},
\end{gather*}
where $1\leq i\leq j\leq m$. Thus,
\begin{gather*}\operatorname{span}\big( \big\{\beta_{ij}^s \,|\, 1\leq i\leq j\leq m, \,
1\leq s\leq 4\big\}\cup \big\{\eta_{kl}^\phi + 2\,f_{kl}\,|\, 1\leq k<l\leq 3\big\} \big) \\
\qquad {} = \mathfrak{sp}(m) \oplus \widetilde{\mathfrak{spin}}(3) \subset \mathfrak{spin}(4m)\oplus\mathfrak{spin}(3)\end{gather*}
annihilates $\phi$, which is consistent with Lemma~\ref{lemma:stabilizer}. Its isotropy group is
\begin{gather*}{\rm Sp}(m){\rm Sp}(1)\equiv {\rm Sp}(m)\widetilde{\rm Spin}(3)\subset {\rm Spin}^3(4m).\end{gather*}
Let us check, for instance,
\begin{gather*}
\big(\eta_{13}^{\phi}+2\kappa_{3\ast}^{m}(f_{13})\big)\phi
 = \sqrt{\frac{3}{m+2}}\frac{1}{\sqrt{m+1}}\!\sum_{j=0}^{m}\frac{1}{{m\choose
j}}\big \{ \eta_{13}^{\phi}\cdot\psi_j\otimes\varphi_{m-j} + \psi_{j}\otimes
2\kappa_{3\ast}^{m}(f_{13})\cdot\varphi_{m-j}\big \} .
\end{gather*}
Observe that
\begin{gather*}
\eta_{13}^{\phi}\cdot \psi_j = -2\big[(j+1)\psi_{j+1}+(j-1-m)\psi_{j-1}\big],
\end{gather*}
so that
\begin{gather*}
\big(\eta_{13}^{\phi}+2\kappa_{3\ast}^{m}(f_{13})\big)\phi =
 \frac{2\sqrt{3}}{\sqrt{(m+2)(m+1)}}\sum_{j=0}^{m}\frac{1}{{m\choose j}}\\
 \hphantom{\big(\eta_{13}^{\phi}+2\kappa_{3\ast}^{m}(f_{13})\big)\phi =}{}\times
 \big\{
 \big[ (-j-1)\psi_{j+1}-(j-1-m)\psi_{j-1}\big] \otimes\varphi_{m-j}\\
\hphantom{\big(\eta_{13}^{\phi}+2\kappa_{3\ast}^{m}(f_{13})\big)\phi =}{}
 + \psi_{j}\otimes \big[ (m-j+1)\varphi_{m-j+1} -(j+1)\varphi_{m-j-1}\big] \big \} \\
\hphantom{\big(\eta_{13}^{\phi}+2\kappa_{3\ast}^{m}(f_{13})\big)\phi}{}
 =
 \frac{2\sqrt{3}}{\sqrt{(m+2)(m+1)}}\sum_{j=0}^{m}\left\{ \left[\frac{-j}{{m\choose
j-1}}+\frac{m-j+1}{{m\choose j}}\right] \psi_{j}\otimes \varphi_{m-j+1}\right.\\
\left.\hphantom{\big(\eta_{13}^{\phi}+2\kappa_{3\ast}^{m}(f_{13})\big)\phi=}{}
 +\left[ \frac{-j+m}{{m\choose j+1}}-\frac{j+1}{{m\choose j}}\right] \psi_{j}\otimes \varphi_{m-j-1}\right\} = 0.
\end{gather*}

\subsubsection[Dimension $n=8$, rank $r=7$]{Dimension $\boldsymbol{n=8}$, rank $\boldsymbol{r=7}$} \label{subsec: Spin(7)}

{\allowdisplaybreaks The spinor $\phi_1\in\Delta_8^+\otimes\Delta_7$ given as follows
\begin{gather}
{\phi_1} = {1\over 2}\big[u_{(-1,-1,-1,-1)}\otimes v_{(1,1,1)} -u_{(1,-1,-1,1)}\otimes v_{(1,1,-1)}+u_{(1,-1,1,-1)}\otimes v_{(1,-1,1)} \nonumber\\
 \hphantom{{\phi_1} =}{} -u_{(1,1,-1,-1)}\otimes v_{(1,-1,-1)} -u_{(-1,-1,1,1)}\otimes v_{(-1,1,1)} +u_{(-1,1,-1,1)}\otimes v_{(-1,1,-1)} \nonumber\\
 \hphantom{{\phi_1} =}{}-u_{(-1,1,1,-1)}\otimes v_{(-1,-1,1)} +u_{(1,1,1,1)}\otimes v_{(-1,-1,-1)} \big]\label{eq: rank 7 pure spinor}
\end{gather}
is pure. The $2$-forms associated to $\phi_1$ are
\begin{alignat*}{3}
& \eta^{\phi_1}_{1 2} = e_{1} e_{2} - e_{3} e_{4} + e_{5} e_{6} + e_{7} e_{8},\qquad &&
 \eta^{\phi_1}_{1 3} = e_{1} e_{3} + e_{2} e_{4} + e_{5} e_{7} - e_{6} e_{8},&\\
& \eta^{\phi_1}_{1 4} = e_{1} e_{4} - e_{2} e_{3} + e_{5} e_{8} + e_{6} e_{7},\qquad &&
 \eta^{\phi_1}_{1 5} = e_{1} e_{5} - e_{2} e_{6} - e_{3} e_{7} - e_{4} e_{8},&\\
& \eta^{\phi_1}_{1 6} = e_{1} e_{6} + e_{2} e_{5} + e_{3} e_{8} - e_{4} e_{7},\qquad &&
 \eta^{\phi_1}_{1 7} = e_{1} e_{7} - e_{2} e_{8} + e_{3} e_{5} + e_{4} e_{6},&\\
& \eta^{\phi_1}_{2 3} = -e_{1} e_{4} + e_{2} e_{3} + e_{5} e_{8} + e_{6} e_{7},\qquad &&
 \eta^{\phi_1}_{2 4} = e_{1} e_{3} + e_{2} e_{4} - e_{5} e_{7} + e_{6} e_{8},&\\
& \eta^{\phi_1}_{2 5} = e_{1} e_{6} + e_{2} e_{5} - e_{3} e_{8} + e_{4} e_{7},\qquad &&
 \eta^{\phi_1}_{2 6} = -e_{1} e_{5} + e_{2} e_{6} - e_{3} e_{7} - e_{4} e_{8},&\\
& \eta^{\phi_1}_{2 7} = e_{1} e_{8} + e_{2} e_{7} + e_{3} e_{6} - e_{4} e_{5},\qquad &&
 \eta^{\phi_1}_{3 4} = -e_{1} e_{2} + e_{3} e_{4} + e_{5} e_{6} + e_{7} e_{8},&\\
& \eta^{\phi_1}_{3 5} = e_{1} e_{7} + e_{2} e_{8} + e_{3} e_{5} - e_{4} e_{6},\qquad &&
 \eta^{\phi_1}_{3 6} = -e_{1} e_{8} + e_{2} e_{7} + e_{3} e_{6} + e_{4} e_{5},&\\
& \eta^{\phi_1}_{3 7} = -e_{1} e_{5} - e_{2} e_{6} + e_{3} e_{7} - e_{4} e_{8},\qquad &&
 \eta^{\phi_1}_{4 5} = e_{1} e_{8} - e_{2} e_{7} + e_{3} e_{6} + e_{4} e_{5},&\\
& \eta^{\phi_1}_{4 6} = e_{1} e_{7} + e_{2} e_{8} - e_{3} e_{5} + e_{4} e_{6},\qquad &&
 \eta^{\phi_1}_{4 7} = -e_{1} e_{6} + e_{2} e_{5} + e_{3} e_{8} + e_{4} e_{7},&\\
& \eta^{\phi_1}_{5 6} = e_{1} e_{2} + e_{3} e_{4} + e_{5} e_{6} - e_{7} e_{8},\qquad &&
 \eta^{\phi_1}_{5 7} = e_{1} e_{3} - e_{2} e_{4} + e_{5} e_{7} + e_{6} e_{8},&\\
& \eta^{\phi_1}_{6 7} = e_{1} e_{4} + e_{2} e_{3} - e_{5} e_{8} + e_{6} e_{7},\qquad &&&
\end{alignat*}
and
\begin{gather*}\operatorname{span}\big\{\hat\eta_{kl}^{\phi_1} \,|\, 1\leq k<l\leq 7 \big\} = \mathfrak{spin}(7) \subset \mathfrak{so}(8).\end{gather*}
The subalgebra
\begin{gather*}\operatorname{span}\big\{\eta_{kl}^{\phi_1} + 2 f_{kl}\,|\, 1\leq k<l\leq 7 \big\} = \widetilde{\mathfrak{spin}}(7) \subset
\mathfrak{spin}(8)\oplus \mathfrak{spin}(7)\end{gather*}
annihilates ${\phi_1}$, which is consistent with Lemma \ref{lemma:stabilizer}. Its isotropy group is
\begin{gather*}\widetilde{\rm Spin}(7)\subset {\rm Spin}^7(8).\end{gather*}

}

{\allowdisplaybreaks Here, we present another twisted pure spinor $\phi_2=-{\rm i}e_1\cdot \phi_1 \in \Delta_8^-\otimes \Delta_7$ (by Lemma \ref{lemma: producing more pure spinors})
\begin{gather*}
\phi_2 ={1\over 2}\big[ u_{(-1,-1,-1,1)}\otimes v_{(1,1,1)} -u_{(1,-1,-1,-1)}\otimes v_{(1,1,-1)}
 + u_{(1,-1,1,1)}\otimes v_{(1,-1,1)} \\
 \hphantom{\phi_2 =}{} - u_{(1,1,-1,1)}\otimes v_{(1,-1,-1)} - u_{(-1,-1,1,-1)}\otimes v_{(-1,1,1)} + u_{(-1,1,-1,-1)}\otimes v_{(-1,1,-1)}\\
\hphantom{\phi_2 =}{}
 - u_{(-1,1,1,1)}\otimes v_{(-1,-1,1)} + u_{(1,1,1,-1)}\otimes v_{(-1,-1,-1)}\big].
\end{gather*}
The $2$-forms associated to $\phi_2$ are
\begin{alignat*}{3}
& \eta^{\phi_2}_{1 2} = -e_{1} e_{2} - e_{3} e_{4} + e_{5} e_{6} + e_{7} e_{8},\qquad &&
 \eta^{\phi_2}_{1 3} = -e_{1} e_{3} + e_{2} e_{4} + e_{5} e_{7} - e_{6} e_{8},&\\
& \eta^{\phi_2}_{1 4} = -e_{1} e_{4} - e_{2} e_{3} + e_{5} e_{8} + e_{6} e_{7},\qquad &&
 \eta^{\phi_2}_{1 5} = -e_{1} e_{5} - e_{2} e_{6} - e_{3} e_{7} - e_{4} e_{8},&\\
& \eta^{\phi_2}_{1 6} = -e_{1} e_{6} + e_{2} e_{5} + e_{3} e_{8} - e_{4} e_{7},\qquad &&
 \eta^{\phi_2}_{1 7} = -e_{1} e_{7} - e_{2} e_{8} + e_{3} e_{5} + e_{4} e_{6},&\\
& \eta^{\phi_2}_{2 3} = e_{1} e_{4} + e_{2} e_{3} + e_{5} e_{8} + e_{6} e_{7},\qquad &&
 \eta^{\phi_2}_{2 4} = -e_{1} e_{3} + e_{2} e_{4} - e_{5} e_{7} + e_{6} e_{8},&\\
& \eta^{\phi_2}_{2 5} = -e_{1} e_{6} + e_{2} e_{5} - e_{3} e_{8} + e_{4} e_{7},\qquad &&
 \eta^{\phi_2}_{2 6} = e_{1} e_{5} + e_{2} e_{6} - e_{3} e_{7} - e_{4} e_{8},&\\
& \eta^{\phi_2}_{2 7} = -e_{1} e_{8} + e_{2} e_{7} + e_{3} e_{6} - e_{4} e_{5},\qquad &&
 \eta^{\phi_2}_{3 4} = e_{1} e_{2} + e_{3} e_{4} + e_{5} e_{6} + e_{7} e_{8},&\\
& \eta^{\phi_2}_{3 5} = -e_{1} e_{7} + e_{2} e_{8} + e_{3} e_{5} - e_{4} e_{6},\qquad &&
 \eta^{\phi_2}_{3 6} = e_{1} e_{8} + e_{2} e_{7} + e_{3} e_{6} + e_{4} e_{5},&\\
& \eta^{\phi_2}_{3 7} = e_{1} e_{5} - e_{2} e_{6} + e_{3} e_{7} - e_{4} e_{8},\qquad &&
 \eta^{\phi_2}_{4 5} = -e_{1} e_{8} - e_{2} e_{7} + e_{3} e_{6} + e_{4} e_{5},&\\
& \eta^{\phi_2}_{4 6} = -e_{1} e_{7} + e_{2} e_{8} - e_{3} e_{5} + e_{4} e_{6},\qquad &&
 \eta^{\phi_2}_{4 7} = e_{1} e_{6} + e_{2} e_{5} + e_{3} e_{8} + e_{4} e_{7},&\\
& \eta^{\phi_2}_{5 6} = -e_{1} e_{2} + e_{3} e_{4} + e_{5} e_{6} - e_{7} e_{8},\qquad &&
 \eta^{\phi_2}_{5 7} = -e_{1} e_{3} - e_{2} e_{4} + e_{5} e_{7} + e_{6} e_{8},&\\
& \eta^{\phi_2}_{6 7} = -e_{1} e_{4} + e_{2} e_{3} - e_{5} e_{8} + e_{6} e_{7},\qquad &&&
\end{alignat*}
and we have another copy
\begin{gather*}\operatorname{span}\big\{\hat\eta_{kl}^{\phi_2} \,|\, 1\leq k<l\leq 7 \big\} = \mathfrak{spin}(7) \subset \mathfrak{so}(8).\end{gather*}
The subalgebra
\begin{gather*}\operatorname{span}\big\{\eta_{kl}^{\phi_2} + 2 f_{kl}\,|\, 1\leq k<l\leq 7 \big\} = \widetilde{\mathfrak{spin}}(7) \subset
\mathfrak{spin}(8)\oplus \mathfrak{spin}(7)\end{gather*}
annihilates ${\phi_2}$. Its isotropy group is another copy
\begin{gather*}\widetilde{\rm Spin}(7)\subset {\rm Spin}^7(8).\end{gather*}

}

The common annihilator of the twisted pure spinors $\phi_1$ and $\phi_2$ is generated by the following elements of $\mathfrak{spin}(8)\oplus \mathfrak{spin}(7)$:
\begin{alignat*}{3}
& e_3 e_4 - e_7 e_8 - f_1 f_2 + f_5 f_6,\qquad && e_2 e_4 - e_6 e_8 + f_1 f_3 - f_5 f_7,&\\
& e_2 e_3 - e_5 e_8 - f_1 f_4 + f_6 f_7,\qquad && e_2 e_6 + e_4 e_8 - f_1 f_5 - f_3 f_7,&\\
& e_2 e_5 + e_3 e_8 + f_1 f_6 + f_4 f_7,\qquad && e_2 e_3 + e_6 e_7 + f_2 f_3 + f_6 f_7,& \\
& e_2 e_4 - e_5 e_7 + f_2 f_4 - f_5 f_7,\qquad && e_2 e_5 + e_4 e_7 + f_2 f_5 + f_4 f_5,& \\
& e_2 e_6 - e_3 e_7 - f_3 f_7 + f_2 f_6,\qquad && e_4 e_6 - e_2 e_8 - f_3 f_5 + f_1 f_7,& \\
& e_3 e_6 + e_2 e_7 + f_3 f_6 + f_2 f_7,\qquad && e_2 e_7 - e_4 e_5 - f_4 f_5 + f_2 f_7,& \\
&e_2 e_8 - e_3 e_5 - f_1 f_7 + f_4 f_6,\qquad && e_3 e_4 + e_5 e_6 + f_3 f_4 + f_5 f_6,&
\end{alignat*}
as can be seen easily by taking appropriate sums of the generators. This elements span a copy of~$\mathfrak{g}_2$.

Let us remark that the two copies of $\mathfrak{spin}(7)$ within $ \mathfrak{so}(8)$ provided by the sets of 2-forms associated to the spinors $\phi_1$ and $\phi_2$ are related by triality, and their intersection is a copy of $\mathfrak{g}_2$.

\section{Special Riemannian holonomy} \label{sec:holonomy}

In this section, we present the geometrical consequences of the existence of parallel twisted pure spinors. In particular, we establish a correspondence between special Riemannnian holonomies and parallel twisted pure spinors.

\subsection[Generic holonomy ${\rm SO}(n)$]{Generic holonomy $\boldsymbol{{\rm SO}(n)}$}

\begin{Proposition}[\cite{Espinosa-Herrera}] Every oriented Riemannian manifold admits a spinorially twisted spin structure such that an associated spinor bundle admits a parallel spinor field.
\end{Proposition}

Indeed, there exists a lift of the horizontal (diagonal) map
of the following diagram
\begin{gather*}
\xymatrix{
 & {\rm Spin}(n)\times_{\mathbb{Z}_2} {\rm Spin}(n) \ar[d]\\
{\rm SO}(n) \ar[ur] \ar[r] & {\rm SO}(n)\times {\rm SO}(n).
}
\end{gather*}
Let $\mathcal{B}$ be the unitary basis of $\Delta_n$ described in Section~\ref{preliminaries} and~$\gamma_n$ be the corresponding real or quaternionic structure of~$\Delta_n$. The twisted spinor $\phi_0\in\Delta_n\otimes \Delta_n$,
\begin{gather*}
\phi_0 := \sum_{\psi\in\mathcal{B}} \psi\otimes \gamma_n(\psi)
\end{gather*}
is ${\rm SO}(n)$ invariant.

\begin{Proposition}[\cite{Espinosa-Herrera}] The $2$-forms associated to $\phi_0$ are multiples of the basic $2$-forms $e_p\wedge e_q$ of $\mathfrak{so}(n)$,
\begin{gather*}\eta_{pq}^{\phi_0} = 2^{[n/2]} e_p\wedge e_q.\end{gather*}
\end{Proposition}

Note that $\phi_0$ is not pure. However, it satisfies the equations
\begin{gather*}
e_pe_q\cdot\phi_0 + \kappa_{n*}^1(f_pf_q)\cdot\phi_0 =0,
\end{gather*}
for $1\leq p<q\leq n$.

\subsection{Holonomy reduction due to parallel twisted pure spinors}

\begin{Definition}Let $M$ be a ${\rm Spin}^r$ Riemannian manifold and $F$ its (locally defined) auxiliary rank $r$ Riemannian vector bundle.
\begin{itemize}\itemsep=0pt
 \item A twisted spinor field $\phi\in\Gamma(S(M,F,m))$ is called pure if $\phi_p$ is pure for every $p\in M$.
 \item Given a connection $\theta$ on the auxiliary bundle $P_{{\rm SO}(r)}$, a twisted spinor field $\psi\in\Gamma(S(M,F,\allowbreak m))$ is parallel if
 \begin{gather*}\nabla_X^\theta\psi=0\end{gather*}
 for all $X\in\Gamma(TM)$.
\end{itemize}
\end{Definition}

\begin{Theorem} Let $M$ be a ${\rm Spin}^{r}$ manifold carrying a twisted pure spinor field $\phi\in\Gamma(S(M,F,\allowbreak m))$ for some $m\in\mathbb{N}$, where $r\geq 3$. Then, in terms of local orthonormal frames $(f_1,\dots,f_r)$ of~$F$,
 \begin{enumerate}\itemsep=0pt
 \item[$1)$] there is a well-defined subbundle $Q\subset \ext^2T^*M$ locally generated by $\{\eta_{kl}^\phi\, |\, 1\leq k<l\leq r\}$;
 \item[$2)$] there is a well-defined subbundle $\hat Q$ of $\operatorname{End}^-(TM)$ locally generated by $\{\hat\eta_{kl}^\phi\, |\, 1\leq k<l\leq r\}$ whose fibre is isomorphic to $\mathfrak{spin}(r)$;
 \item[$3)$] there is a rank $r$ almost even-Clifford Hermitian structure induced by the local maps
\begin{align*}
\big({\rm Cl}^0(F)\big)^2 &\longrightarrow \operatorname{End}^-(TM),\\
 f_{ij}&\mapsto \hat\eta_{kl}^\phi.
\end{align*}
 \end{enumerate}
\end{Theorem}
\begin{proof}The proof follows from Proposition \ref{prop: Clifford transfer} and Lemma \ref{lemma: orbit of pure spinor}.\end{proof}

\begin{Theorem}\label{theo: parallelness} Let $M$ be a ${\rm Spin}^r$ Riemannian manifold whose auxiliary bundle $P_{{\rm SO}(r)}$ is endowed with a connection~$\theta$. If $M$ carries a parallel twisted pure spinor field $\phi\in\Gamma(S(M,F,m))$ for some $m\in\mathbb{N}$, $r\geq 3$, then the manifold~$M$ admits a~rank~$r$ parallel even-Clifford Hermitian structure and its holonomy algebra is contained in one of the algebras $\mathfrak{g}$ of Table~{\rm \ref{table: normalizer subalgebras}}.
\begin{table}[h!]\centering
$
\begin{array}{|c|c|}
 \hline
 r \mbox{\ {\rm (mod 8)} } & \mathfrak{g}
 \rule{0pt}{3ex}
 \\
 \hline
 0 & \mathfrak{so}(m_1)\oplus\mathfrak{so}(m_2)\oplus {\mathfrak{spin}}(r)\rule{0pt}{3ex}\\
 \hline
 1,7 & \mathfrak{so}(m)\oplus {\mathfrak{spin}}(r)\rule{0pt}{3ex}\\
\hline
 2,6 & \mathfrak{u}(m)\oplus {\mathfrak{spin}}(r)\rule{0pt}{3ex}\\
\hline
 3,5 & \mathfrak{sp}(m)\oplus {\mathfrak{spin}}(r)\rule{0pt}{3ex}\\
\hline
 4 & \mathfrak{sp}(m_1)\oplus \mathfrak{sp}(m_2)\oplus {\mathfrak{spin}}(r)\rule{0pt}{3ex}\\
\hline
 \end{array}
$
\caption{}\label{table: normalizer subalgebras}
\end{table}
\end{Theorem}

\begin{proof}\allowdisplaybreaks Suppose $\nabla^\theta \phi=0$. Let $(e_1,\dots,e_n)$ and $(f_1,\dots,f_r)$ be local orthonormal frames for $TM$ and $F$ respectively,
and $X\in\Gamma(TM)$. Recall that
\begin{gather*}
\nabla_Xe_j = \omega_{1j}(X) e_1 +\dots +\omega_{nj}(X) e_n ,\\
\nabla_Xf_j = \theta_{1j}(X) f_1 +\dots +\theta_{rj}(X) f_r.
\end{gather*}
On the one hand
\begin{gather*}
 \nabla_X \big(\eta_{kl}^\phi(e_s,e_t)\big) = \big(\nabla_X \eta_{kl}^\phi\big)(e_s,e_t) + \eta_{kl}^\phi(\nabla_Xe_s,e_t)+\eta_{kl}^\phi(e_s,\nabla_Xe_t)\\
\hphantom{\nabla_X \big(\eta_{kl}^\phi(e_s,e_t)\big)}{} = \big(\nabla_X \eta_{kl}^\phi\big)(e_s,e_t) +
\sum_{a=1}^n\omega_{as}(X)\eta_{kl}^\phi(e_a,e_t)+\sum_{a=1}^n\omega_{at}(X)\eta_{kl}^\phi(e_s,e_a)
\end{gather*}
and, on the other,
\begin{gather*}
 \nabla_X \big(\eta_{kl}^\phi(e_s,e_t)\big) = \nabla_X \big\langle e_s e_t\cdot \kappa_{r*}^m(f_{kl}) \cdot \phi,\phi \big\rangle \\
\hphantom{\nabla_X \big(\eta_{kl}^\phi(e_s,e_t)\big)}{}
 = \big\langle \nabla_X (e_s e_t)\cdot \kappa_{r*}^m(f_{kl}) \cdot \phi,\phi \big\rangle
+\big\langle e_s e_t\cdot \nabla_X (\kappa_{r*}^m(f_{kl})) \cdot \phi,\phi \big\rangle \\
\hphantom{\nabla_X \big(\eta_{kl}^\phi(e_s,e_t)\big)=}{}
 +\big\langle e_s e_t\cdot \kappa_{r*}^m(f_{kl}) \cdot \nabla_X^\theta \phi,\phi \big\rangle
+\big\langle e_s e_t\cdot \kappa_{r*}^m(f_{kl}) \cdot \phi,\nabla_X^\theta\phi \big\rangle \\
\hphantom{\nabla_X \big(\eta_{kl}^\phi(e_s,e_t)\big)}{}
 = \big\langle \nabla_X (e_s e_t)\cdot \kappa_{r*}^m(f_{kl}) \cdot \phi,\phi \big\rangle
+\big\langle e_s e_t\cdot \nabla_X (\kappa_{r*}^m(f_{kl})) \cdot \phi,\phi \big\rangle \\
\hphantom{\nabla_X \big(\eta_{kl}^\phi(e_s,e_t)\big)}{}
 = \bigg\langle \sum_{a=1}^n\omega_{as}(X)e_a e_t\cdot \kappa_{r*}^m(f_{kl}) \cdot \phi,\phi \bigg\rangle
+\bigg\langle \sum_{a=1}^n\omega_{at}(X)e_s e_a\cdot \kappa_{r*}^m(f_{kl}) \cdot \phi,\phi \bigg\rangle \\
\hphantom{\nabla_X \big(\eta_{kl}^\phi(e_s,e_t)\big)=}{}
 +\bigg\langle\! \sum_{a=1}^r\theta_{ak}(X)e_s e_t\cdot \kappa_{r*}^m(f_{al}) \cdot \phi,\phi \bigg\rangle
+\bigg\langle\! \sum_{a=1}^r\theta_{al}(X)e_s e_t\cdot \kappa_{r*}^m(f_{ka})) \cdot \phi,\phi \bigg\rangle \\
\hphantom{\nabla_X \big(\eta_{kl}^\phi(e_s,e_t)\big)}{}
 = \sum_{a=1}^n\omega_{as}(X)\big\langle e_a e_t\cdot \kappa_{r*}^m(f_{kl}) \cdot \phi,\phi \big\rangle
+\sum_{a=1}^n\omega_{at}(X)\big\langle e_s e_a\cdot \kappa_{r*}^m(f_{kl}) \cdot \phi,\phi \big\rangle \\
\hphantom{\nabla_X \big(\eta_{kl}^\phi(e_s,e_t)\big)=}{}
 +\sum_{a=1}^r\theta_{ak}(X)\big\langle e_s e_t\cdot \kappa_{r*}^m(f_{al}) \cdot \phi,\phi \big\rangle
+\sum_{a=1}^r\theta_{al}(X)\big\langle e_s e_t\cdot \kappa_{r*}^m(f_{ka}) \cdot \phi,\phi \big\rangle \\
\hphantom{\nabla_X \big(\eta_{kl}^\phi(e_s,e_t)\big)}{}
 = \sum_{a=1}^n\omega_{as}(X)\eta_{kl}^\phi(e_a,e_t)
+\sum_{a=1}^n\omega_{at}(X)\eta_{kl}^\phi(e_s,e_a)\\
\hphantom{\nabla_X \big(\eta_{kl}^\phi(e_s,e_t)\big)=}{}
 +\sum_{a=1}^r\theta_{ak}(X)\eta_{al}^\phi(e_s, e_t)
+\sum_{a=1}^r\theta_{al}(X)\eta_{ka}^\phi( e_s, e_t).
\end{gather*}
Thus,
\begin{gather*}
 \big(\nabla_X \big(\eta_{kl}^\phi\big)\big)(e_s,e_t)
= \sum_{a=1}^r\theta_{ak}(X)\eta_{al}^\phi(e_s, e_t)
+\sum_{a=1}^r\theta_{al}(X)\eta_{ka}^\phi( e_s, e_t),
\end{gather*}
i.e.,
\begin{gather*}
 \nabla \eta_{kl}^\phi = \sum_{a=1}^r\theta_{ak}\otimes\eta_{al}^\phi +\theta_{al}\otimes\eta_{ka}^\phi.
\end{gather*}
Furthermore, for $X,Y\in\Gamma(TM)$,
\begin{gather}
R^M(X,Y)\big(\hat\eta_{kl}^\phi\big)= \sum_{a=1}^r \Theta_{ak}(X,Y)\hat\eta_{al}^\phi+\Theta_{al}(X,Y)\hat\eta_{ka}^\phi, \label{eq: curvature identity endomorphisms}
\end{gather}
which means
\begin{gather*}R^M(X,Y)\in N_{\mathfrak{so}(n)}(\mathfrak{spin}(r)).\end{gather*}
Such normalizer subalgebras were computed in \cite{AH} and the corresponding Lie groups were computed in \cite{AGH}.
\end{proof}

\subsection{Curvature calculations}

In this section we carry out various curvature calculations (as in \cite{Moroianu-Semmelmann}) which lead us to the ge\-ne\-ralization of the second Seiberg--Witten equation.

Recall that if $\xi$ and $\zeta$ are antisymmetic two forms, and $\hat\xi$ and $\hat\zeta$ are the corresponding antisymmetric endomorphisms with respect to a positive definite inner product, then
\begin{gather*}\tr\big(\hat\xi\hat\zeta\big) = -2 \xi\bullet\zeta,\end{gather*}
where $\bullet$ denotes the induced inner product on 2-forms.

If $X,Y,Z,W\in\Gamma(TM)$, let
\begin{gather*}R^M(X,Y,Z,W)=\big\langle R^M(X,Y),Z,W\big\rangle ,\end{gather*}
denote the curvature 4-form of $M$,
\begin{gather*}
 \Theta = \sum_{1\leq k<l \leq r} \Theta_{kl}\otimes f_{kl}
\end{gather*}
the curvature of $\theta$
and
\begin{gather*}
 \eta^\phi = \sum_{1\leq k<l \leq r} \eta_{kl}^\phi\otimes f_{kl},
\end{gather*}
the 2-form with values in $\ext^2F$ associated to $\phi$.

\begin{Lemma}\label{Lemma4.6} Let $M$ be a ${\rm Spin}^r$ Riemannian manifold whose auxiliary bundle $P_{{\rm SO}(r)}$ is endowed with a connection $\theta$. If $M$ carries a parallel twisted pure spinor field $\phi\in\Gamma(S(M,F,m))$ for some $m\in\mathbb{N}$, $r\geq 3$, then
\begin{gather}
\Theta=-{4\over \dim(M)}R^M\bullet \eta^\phi,\label{eq: first version second SW equation}
\end{gather}
where $R^M\bullet \eta^\phi$ denotes the image of $R^M\otimes \eta^\phi$ under the map $($contraction of underlined factors$)$
\begin{gather*}
\big(\ext^2 T^*M\otimes \underline{\ext^2 T^*M}\big)\otimes\big(\underline{\ext^2 T^*M}\otimes\ext^2 F\big)
 \xrightarrow{\,\,\,\,\bullet\,\,\,\,} \ext^2 T^*M\otimes\ext^2 F.
\end{gather*}
\end{Lemma}
\begin{proof} Let $X,Y\in\Gamma(TM)$. Since $r\geq 3$, multiply~\eqref{eq: curvature identity endomorphisms} by $\hat{\eta}^\phi_{kp}$ on the left, $p\not = k,l$,
\begin{gather*}
\hat{\eta}^\phi_{kp}{R(X,Y)}\hat{\eta}^\phi_{kl}-\hat{\eta}^\phi_{pl}{R(X,Y)} \\
\qquad {} = \Theta_{kp}(X,Y)\eta_{kl}^\phi + \Theta_{pl}(X,Y)(-{\rm Id}_{TM}) + \sum_{a\not= k,l,p} \Theta_{ak}(X,Y) \eta_{kp}^\phi\hat{\eta}^\phi_{al}+ \Theta_{al}(X,Y)\hat{\eta}^\phi_{pa},
\end{gather*}
where the last summation is absent if $r=3$. Since
\begin{gather*}
\tr\big(\eta_{kp}^\phi\hat{\eta}^\phi_{lq}\big) = 0, \qquad \text{if} \quad r\not = 4,\\
\tr\big(\eta_{kl}^\phi\big) = 0,
\end{gather*}
it makes sense to take the trace on the right-hand side in order to isolate $\Theta_{pl}(X,Y)$. The left-hand side gives
\begin{gather*}
\tr\big(\hat{\eta}^\phi_{kp}{R(X,Y)}\hat{\eta}^\phi_{kl}-\hat{\eta}^\phi_{pl}{R(X,Y)}\big)\\
\qquad {} =
\tr\big(\hat{\eta}^\phi_{kp}{R(X,Y)}\hat{\eta}^\phi_{kl}\big)-\tr\big(\hat{\eta}^\phi_{pl}{R(X,Y)}\big)
 = 2\tr\big(\hat{\eta}^\phi_{lp}{R(X,Y)}\big),
\end{gather*}
and the right-hand side
\begin{gather*}
\tr\left(\Theta_{kp}(X,Y)\hat\eta_{kl}^\phi + \Theta_{pl}(X,Y)(-{\rm Id}_{TM})
 + \sum_{q\not= k,l,p} \Theta_{kq}(X,Y) \eta_{kp}^\phi\hat{\eta}^\phi_{lq}+ \Theta_{ql}(X,Y)\hat{\eta}^\phi_{pq}\right)\\
 \qquad{} = -\dim(M)\Theta_{pl}(X,Y).
\end{gather*}
Thus
\begin{gather}
2\tr\big(\hat{\eta}^\phi_{pl}{R(X,Y)}\big)=\dim(M)\Theta_{pl}(X,Y), \label{eq: isolation Theta}
\end{gather}
which proves Lemma~\ref{Lemma4.6}.
\end{proof}

\begin{Remark}\eqref{eq: first version second SW equation} is, in fact, the prototype of the second Clifford monopole equation,
where $-R^M\bullet$ can be substituted by a~symmetric endomorphisms of $\ext^2 T^*M$.
\end{Remark}

\begin{Lemma}\label{lemma: identity curvatures}Let $M$ be a ${\rm Spin}^r$ $n$-dimensional Riemannian manifold whose auxiliary bundle $P_{{\rm SO}(r)}$ is endowed with a connection $\theta$. If $M$ carries a parallel twisted pure spinor field $\phi\in\Gamma(S(M,F,m))$ for some $m\in\mathbb{N}$, $r\geq 3$, $r\not= 4$, $n\not=8$, $n+4r-16\not=0$, $n+8r-16\not=0$,
then
\begin{gather*}
\hat\Theta_{kl} = {{\rm R}\over n\left({n\over 4}+2r-4\right)}\hat\eta_{kl}^\phi
\end{gather*}
for $1\leq k<l\leq r$, and $M$ is Einstein.
\end{Lemma}
\begin{proof}\allowdisplaybreaks Let $X,Y\in\Gamma(TM)$. In order to simplify notation and carry out further calculations, let
\begin{gather} \vartheta_{pl}(X,Y):= \tr\big(\hat{\eta}^\phi_{pl}{R(X,Y)}\big)=\tr\big({R(X,Y)}\hat{\eta}^\phi_{pl}\big)\nonumber\\
 \hphantom{\vartheta_{pl}(X,Y)}{} =
 \sum_i \big\langle {R(X,Y)}\hat{\eta}^\phi_{pl}(e_i),e_i\big\rangle =
 -\sum_i \big\langle {R(X,Y)}(e_i),\hat{\eta}^\phi_{pl}(e_i)\big\rangle \nonumber\\
\hphantom{\vartheta_{pl}(X,Y)}{} =
 \sum_i \big\langle {R(X,e_i)}\big(\hat{\eta}^\phi_{pl}(e_i)\big), Y\big\rangle +\sum_i \big\langle {R\big(X,\hat{\eta}^\phi_{pl}(e_i)\big) }(Y),e_i\big\rangle \nonumber\\
\hphantom{\vartheta_{pl}(X,Y)}{}
 =
 2\sum_i \big\langle {R(X,e_i)}\big(\hat{\eta}^\phi_{pl}(e_i)\big), Y\big\rangle .\label{eq: definition vartheta}
\end{gather}
Now multiply \eqref{eq: curvature identity endomorphisms}
on the left by $\hat\eta_{kl}^\phi$
\begin{gather}
\hat{\eta}^\phi_{kl}{R(X,Y)}\hat{\eta}^\phi_{kl}+{R(X,Y)} =
\hat{\eta}^\phi_{kl} \sum_{q\not= k,l} \Theta_{kq}(X,Y) \hat{\eta}^\phi_{lq}+ \Theta_{ql}(X,Y)\hat{\eta}^\phi_{kq}\nonumber\\
\hphantom{\hat{\eta}^\phi_{kl}{R(X,Y)}\hat{\eta}^\phi_{kl}+{R(X,Y)}}{} =
- \sum_{q\not= k,l} \Theta_{kq}(X,Y) \hat{\eta}^\phi_{kq}+ \Theta_{lq}(X,Y)\hat{\eta}^\phi_{lq}\nonumber\\
\hphantom{\hat{\eta}^\phi_{kl}{R(X,Y)}\hat{\eta}^\phi_{kl}+{R(X,Y)}}{} =
 2\Theta_{kl}(X,Y)\hat\eta_{kl}^\phi -\sum_{q} \Theta_{kq}(X,Y) \hat{\eta}^\phi_{kq}+ \Theta_{lq}(X,Y)\hat{\eta}^\phi_{lq}. \label{eq: curvature identity endomorphisms 2}
\end{gather}
Consider the sum
\begin{gather*}
\sum_{i,j}\big\langle \hat{\eta}^\phi_{kl}{R(X,e_i)}\hat{\eta}^\phi_{kl}(e_i),e_j\big\rangle e_j +\sum_{i,j}\langle {R(X,e_i)}(e_i),e_j\rangle e_j\\
{}= 2\sum_{i,j}\big\langle \Theta_{kl}(X,e_i)\hat\eta_{kl}^\phi(e_i),e_j\big\rangle e_j -\sum_{i,j}\bigg\langle \sum_{q} \Theta_{kq}(X,e_i) \hat{\eta}^\phi_{kq}(e_i)
 + \Theta_{lq}(X,e_i)\hat{\eta}^\phi_{lq}(e_i),e_j\bigg\rangle e_j.
\end{gather*}
The left-hand side is
\begin{gather*}
\sum_{i,j}\big\langle \hat{\eta}^\phi_{kl}{R(X,e_i)}\hat{\eta}^\phi_{kl}(e_i),e_j\big\rangle e_j
+\sum_{i,j}\langle {R(X,e_i)}(e_i),e_j\rangle e_j\\
\qquad{} =
-\sum_{i,j}\big\langle {R(X,e_i)}\hat{\eta}^\phi_{kl}(e_i),\hat{\eta}^\phi_{kl}e_j\big\rangle e_j +\sum_{i,j}\langle {R(X,e_i)}(e_i),e_j\rangle e_j\\
 \qquad{} =
\hat\eta_{kl}^\phi\sum_{i,j}\big\langle {R(X,e_i)}\hat{\eta}^\phi_{kl}(e_i),\hat{\eta}^\phi_{kl}e_j\big\rangle \hat\eta_{kl}^\phi (e_j)+{\rm Ric}(X)\\
\qquad{} = \hat\eta_{kl}^\phi\sum_{j}{1\over 2}\vartheta_{kl}(X,e_j')e_j'+{\rm Ric}(X)
 =
{1\over 2}\hat\eta_{kl}^\phi\big(\hat\vartheta_{kl}(X)\big)+{\rm Ric}(X) ,
\end{gather*}
and the right-hand side is
\begin{gather*}
 2\sum_{i,j}\big\langle \Theta_{kl}(X,e_i)\hat\eta_{kl}^\phi(e_i),e_j\big\rangle e_j
 -\sum_{i,j}\bigg\langle \sum_{q} \Theta_{kq}(X,e_i) \hat{\eta}^\phi_{kq}(e_i)
 + \Theta_{lq}(X,e_i)\hat{\eta}^\phi_{lq}(e_i),e_j\bigg\rangle e_j\\
{} =
- 2\sum_{j}\big\langle \hat\Theta_{kl}(X),\hat\eta_{kl}^\phi(e_j)\big\rangle e_j
 +\sum_{j}\sum_q\big\langle \hat\Theta_{kq}(X),\hat{\eta}^\phi_{kq}(e_j)\big\rangle e_j
 + \sum_{j}\sum_q\big\langle \hat\Theta_{lq}(X),\hat{\eta}^\phi_{lq}(e_j)\big\rangle e_j\\
{} = 2\hat\eta_{kl}^\phi\big(\hat\Theta_{kl}(X)\big) -\sum_q\hat{\eta}^\phi_{kq}\big(\hat\Theta_{kq}(X)\big)
 - \sum_q\hat{\eta}^\phi_{lq}\big(\hat\Theta_{lq}(X)\big).
\end{gather*}
Thus
\begin{gather*}
{\rm Ric}(X) +{1\over 2}\hat\eta_{kl}^\phi\big(\hat\vartheta_{kl}(X)\big) =
 2\hat\eta_{kl}^\phi\big(\hat\Theta_{kl}(X)\big) -\sum_q\hat{\eta}^\phi_{kq}\big(\hat\Theta_{kq}(X)\big)
 - \sum_q\hat{\eta}^\phi_{lq}\big(\hat\Theta_{lq}(X)\big).
\end{gather*}
By \eqref{eq: isolation Theta} and \eqref{eq: definition vartheta}
\begin{gather*}\hat\vartheta_{kl}(X)={n\over 2}\hat\Theta_{kl}(X),\end{gather*}
so that
\begin{gather*}
{\rm Ric}(X) +{n\over 4}\hat\eta_{kl}^\phi(\hat\Theta_{kl}(X)) =
 2\hat\eta_{kl}^\phi\big(\hat\Theta_{kl}(X)\big) -\sum_q\hat{\eta}^\phi_{kq}\big(\hat\Theta_{kq}(X)\big) - \sum_q\hat{\eta}^\phi_{lq}\big(\hat\Theta_{lq}(X)\big),
\end{gather*}
i.e.,
\begin{gather*}
 0={\rm Ric}^M + \left({n\over 4}-2\right)\hat\eta_{kl}^\phi\circ\hat\Theta_{kl} +\sum_q\hat{\eta}^\phi_{kq}\circ\hat\Theta_{kq} + \sum_q\hat{\eta}^\phi_{lq}\circ\hat\Theta_{lq}.
\end{gather*}
Let
\begin{gather*}T_k:=\sum_q \hat{\eta}^\phi_{kq}\circ\hat\Theta_{kq}.\end{gather*}
Now
\begin{gather}
0={\rm Ric} + \left({n\over 4}-2\right)\hat\eta_{kl}^\phi\circ\hat\Theta_{kl} + T_k + T_l . \label{eq: Ricci curvature identity 2}
\end{gather}
Recall that $k\not = l$, and consider
\begin{gather*}
0=\sum_{1\leq l\leq r,\,l\not= k} {\rm Ric} +\left({n\over 4}-2\right)\sum_{1\leq l\leq r,\,l\not= k}\hat\eta_{kl}^\phi\circ\hat\Theta_{kl}
 +\sum_{1\leq l\leq r,\,l\not= k} T_k + \sum_{1\leq l\leq r,\,l\not= k} T_l,
\end{gather*}
i.e.,
\begin{gather*}
0 = (r-1) {\rm Ric}+ \left({n\over 4}-2\right)T_k +(r-1) T_k + \sum_{l} T_l - T_k\\
\hphantom{0}{} = (r-1) {\rm Ric}+ \left({n\over 4}-2\right)T_k +(r-2) T_k + \sum_l T_l\\
\hphantom{0}{} = (r-1) {\rm Ric}+ \left(r+{n\over 4}-4\right)T_k + \sum_l T_l.
\end{gather*}
Thus, if $p\not =k$,
\begin{gather*}
0 =(r-1) {\rm Ric}+ \left(r+{n\over 4}-4\right)T_k + \sum_l T_l,\\
0 =(r-1) {\rm Ric}+ \left(r+{n\over 4}-4\right)T_p + \sum_l T_l,
\end{gather*}
and, if $(r+n/4-4)\not =0$, we have
\begin{gather*}T_k=T_p=:\mathbb{T}.\end{gather*}
Thus,
\begin{gather*}
0=(r-1) {\rm Ric}+\left(2r+{n\over 4}-4\right)\mathbb{T} ,
\end{gather*}
and if $4-{n\over 4}-2r\not=0$,
\begin{gather*}
 \mathbb{T}= {r-1\over \left(4-{n\over 4}-2r\right)} {\rm Ric} .
\end{gather*}
Going back to \eqref{eq: Ricci curvature identity 2}
\begin{gather*}
0 ={\rm Ric}+ \left({n\over 4}-2\right)\hat\eta_{kl}^\phi\circ\hat\Theta_{kl} +T_k + T_l =
{\rm Ric}+ \left({n\over 4}-2\right)\hat\eta_{kl}^\phi\circ\hat\Theta_{kl} +{2r-2\over \left(4-{n\over 4}-2r\right)} {\rm Ric} ,
\end{gather*}
i.e.,
\begin{gather*}
\left({n\over 4}-2\right)\hat\eta_{kl}^\phi\circ\hat\Theta_{kl} =-{\rm Ric} -{2r-2\over \left(4-{n\over 4}-2r\right)} {\rm Ric}
={\left({n\over 4}-2\right)\over \left(4-{n\over 4}-2r\right)} {\rm Ric}.
\end{gather*}
If $n\not = 8$,
\begin{gather}
\hat\eta_{kl}^\phi\circ\hat\Theta_{kl} = {1\over \left(4-{n\over 4}-2r\right)} {\rm Ric}^M. \label{eq: eta theta Ric}
\end{gather}
Since ${\rm Ric}$ is symmetric, $\hat\Theta_{kl}$ and $\hat\eta_{kl}^\phi$ commute, and they commute with ${\rm Ric}$ for all $1\leq k<l\leq r$. If $k$, $l$, $p$, $q$ are all different,
\begin{gather*}
\tr\big(\hat\Theta_{pq}\hat\eta_{kl}^\phi\big) =
-\tr\big(\hat\Theta_{pq}\hat\eta_{pq}^\phi\hat\eta_{pq}^\phi\hat\eta_{sk}^\phi\hat\eta_{sl}^\phi\big) =
-{1\over \left(4-{n\over 4}-2r\right)}\tr\big({\rm Ric}^M\hat\eta_{pq}^\phi\hat\eta_{sk}^\phi\hat\eta_{sl}^\phi\big) \\
\hphantom{\tr\big(\hat\Theta_{pq}\hat\eta_{kl}^\phi\big)}{} =
-{1\over \left(4-{n\over 4}-2r\right)}\tr\big({\rm Ric}^M\hat\eta_{sk}^\phi\hat\eta_{pq}^\phi\hat\eta_{sl}^\phi\big) =
-{1\over \left(4-{n\over 4}-2r\right)}\tr\big(\hat\eta_{sk}^\phi{\rm Ric}^M\hat\eta_{pq}^\phi\hat\eta_{sl}^\phi\big) \\
\hphantom{\tr\big(\hat\Theta_{pq}\hat\eta_{kl}^\phi\big)}{} =
-{1\over \left(4-{n\over 4}-2r\right)}\tr\big({\rm Ric}^M\hat\eta_{pq}^\phi\hat\eta_{sl}^\phi\hat\eta_{sk}^\phi\big) =
{1\over \left(4-{n\over 4}-2r\right)}\tr\big({\rm Ric}^M\hat\eta_{pq}^\phi\hat\eta_{kl}^\phi\big) \\
\hphantom{\tr\big(\hat\Theta_{pq}\hat\eta_{kl}^\phi\big)}{}
 = -\tr\big(\hat\Theta_{pq}\hat\eta_{kl}^\phi\big) ,
\end{gather*}
so that
\begin{gather}
\tr\big(\hat\Theta_{pq}\hat\eta_{kl}^\phi\big)=0.\label{eq: vanishing traces 3}
\end{gather}
If $p$, $k$, $l$ are all different
\begin{gather*}
\tr\big(\hat\Theta_{pk}\hat\eta_{kl}^\phi\big) = -\tr\big(\hat\Theta_{pk}\hat\eta_{pk}^\phi\hat\eta_{pk}^\phi\hat\eta_{kl}^\phi\big)
= {1\over \left(4-{n\over 4}-2r\right)}\tr\big({\rm Ric}^M\hat\eta_{kl}^\phi\hat\eta_{pk}^\phi\big) \\
\hphantom{\tr\big(\hat\Theta_{pk}\hat\eta_{kl}^\phi\big)}{}
 = {1\over \left(4-{n\over 4}-2r\right)}\tr\big(\hat\eta_{kl}^\phi{\rm Ric}^M\hat\eta_{pk}^\phi\big) =
 {1\over \left(4-{n\over 4}-2r\right)}\tr\big({\rm Ric}^M\hat\eta_{pk}^\phi\hat\eta_{kl}^\phi\big) \\
\hphantom{\tr\big(\hat\Theta_{pk}\hat\eta_{kl}^\phi\big)}{}=
 -\tr\big(\hat\Theta_{pk}\hat\eta_{kl}^\phi\big) ,
\end{gather*}
i.e.,
\begin{gather}
\tr\big(\hat\Theta_{pk}\hat\eta_{kl}^\phi\big) =0.\label{eq: vanishing traces 4}
\end{gather}
Now set $X=\hat\eta_{st}^\phi (e_a)$ and $Y=e_a$ in \eqref{eq: curvature identity endomorphisms 2} and sum over $a$
\begin{gather*}
\sum_a\hat{\eta}^\phi_{kl}{R\big(\hat\eta_{st}^\phi(e_a),e_a\big)}\hat{\eta}^\phi_{kl}+{R\big(\hat\eta_{st}^\phi(e_a),e_a\big)}\\
\qquad{} = 2\sum_a\Theta_{kl}\big(\hat\eta_{st}^\phi(e_a),e_a\big)\hat\eta_{kl}^\phi
 -\sum_a\sum_{q} \Theta_{kq}\big(\hat\eta_{st}^\phi(e_a),e_a\big) \hat{\eta}^\phi_{kq}+ \Theta_{lq}\big(\hat\eta_{st}^\phi(e_a),e_a\big)\hat{\eta}^\phi_{lq}.
\end{gather*}
The right-hand side is equal to
\begin{gather*}
 2\bigg(\sum_a\big\langle \hat\Theta_{kl}\hat\eta_{st}^\phi(e_a),e_a\big\rangle \bigg)\hat\eta_{kl}^\phi
 -\sum_{q} \bigg[\bigg(\!\sum_a\big\langle \hat\Theta_{kq}\hat\eta_{st}^\phi(e_a),e_a\big\rangle \bigg) \hat{\eta}^\phi_{kq}
 + \bigg(\!\sum_a\big\langle \hat\Theta_{lq}\hat\eta_{st}^\phi(e_a),e_a\big\rangle \bigg)\hat{\eta}^\phi_{lq}\bigg] \\
\qquad{} = 2\tr\big(\hat\Theta_{kl}\hat\eta_{st}^\phi\big)\hat\eta_{kl}^\phi
 -\sum_{q} \big[\tr\big(\hat\Theta_{kq}\hat\eta_{st}^\phi\big) \hat{\eta}^\phi_{kq}
 + \tr\big(\hat\Theta_{lq}\hat\eta_{st}^\phi\big)\hat{\eta}^\phi_{lq}\big] .
\end{gather*}
Now let us analyze the terms of the left-hand side by evaluating it on a vector field $Z\in\Gamma(TM)$
\begin{gather*}
\sum_a\hat{\eta}^\phi_{kl}{R\big(\hat\eta_{st}^\phi(e_a),e_a\big)}\hat{\eta}^\phi_{kl} (Z)
 =\hat{\eta}^\phi_{kl}\bigg(\sum_a{R\big(\hat\eta_{st}^\phi(e_a),e_a\big)}\hat{\eta}^\phi_{kl}(Z)\bigg)\\
 \qquad {} =\hat{\eta}^\phi_{kl}\bigg(\sum_{a,b}\big\langle {R\big(\hat\eta_{st}^\phi(e_a),e_a\big)}\hat{\eta}^\phi_{kl}(Z),e_b\big\rangle e_b\bigg)
 = -\sum_{a,b}\big\langle {R\big(\hat\eta_{st}^\phi(e_a),e_a\big)}e_b,\hat{\eta}^\phi_{kl}(Z)\big\rangle \hat{\eta}^\phi_{kl}(e_b)\\
 \qquad {}= \sum_{a,b}\big\langle \hat{\eta}^\phi_{kl}{R\big(\hat\eta_{st}^\phi(e_a),e_a\big)}e_b,Z\big\rangle \hat{\eta}^\phi_{kl}(e_b)
 = -\sum_{a,b}\big\langle \hat{\eta}^\phi_{kl}{R\big(\hat\eta_{st}^\phi(e_a),e_a\big)}Z,e_b\big\rangle \hat{\eta}^\phi_{kl}(e_b)\\
 \qquad{} = \sum_{a,b}\big\langle {R\big(\hat\eta_{st}^\phi(e_a),e_a\big)}Z,\hat{\eta}^\phi_{kl}(e_b)\big\rangle \hat{\eta}^\phi_{kl}(e_b)
 = \sum_{a,b}\big\langle {R\big(\hat\eta_{st}^\phi(e_a),e_a\big)}Z,e_b'\big\rangle e_b'\\
 \qquad{} = -\sum_{a,b}\big\langle {R\big(e_a,\hat\eta_{st}^\phi(e_a)\big)}Z,e_b'\big\rangle e_b'
 = -\sum_{a,b}\big\langle {R(Z,e_b')}e_a,\hat\eta_{st}^\phi(e_a)\big\rangle e_b'\\
 \qquad{} = \sum_{a,b}\big\langle \hat\eta_{st}^\phi{R(Z,e_b')}e_a,e_a\big\rangle e_b' =
\sum_{b}\bigg(\sum_{a}\big\langle \hat\eta_{st}^\phi{R(Z,e_b')}e_a,e_a\big\rangle \bigg)e_b'\\
 \qquad{} = \sum_{b}\tr\big(\hat\eta_{st}^\phi{R(Z,e_b')}\big)e_b' =\sum_{b}\vartheta_{st}(Z,e_b')e_b' = \hat\vartheta_{st}(Z),
\end{gather*}
i.e.,
\begin{gather*}
\sum_a\hat{\eta}^\phi_{kl}{R\big(\hat\eta_{st}^\phi(e_a),e_a\big)}\hat{\eta}^\phi_{kl} = \hat\vartheta_{st},
\end{gather*}
and, as can also be seen from the middle of the previous calculation,
\begin{gather*}
\sum_a{R\big(\hat\eta_{st}^\phi(e_a),e_a\big)} =\hat\vartheta_{st}.
\end{gather*}
Thus, we have
\begin{gather*}
2\hat\vartheta_{st} =2\tr\big(\hat\Theta_{kl}\hat\eta_{st}^\phi\big)\hat\eta_{kl}^\phi
 -\sum_{q} \big[\tr\big(\hat\Theta_{kq}\hat\eta_{st}^\phi\big) \hat{\eta}^\phi_{kq}
 + \tr\big(\hat\Theta_{lq}\hat\eta_{st}^\phi\big)\hat{\eta}^\phi_{lq}\big],
\end{gather*}
which, by \eqref{eq: first version second SW equation} is equivalent to
\begin{gather*}
n\hat\Theta_{st} =2\tr\big(\hat\Theta_{kl}\hat\eta_{st}^\phi\big)\hat\eta_{kl}^\phi
 -\sum_{q} \big[\tr\big(\hat\Theta_{kq}\hat\eta_{st}^\phi\big) \hat{\eta}^\phi_{kq}
 + \tr\big(\hat\Theta_{lq}\hat\eta_{st}^\phi\big)\hat{\eta}^\phi_{lq}\big] .
\end{gather*}
If $s=k<l\not =t$,
\begin{gather*}
n\hat\Theta_{st} =
2\tr\big(\hat\Theta_{sl}\hat\eta_{st}^\phi\big)\hat\eta_{sl}^\phi
 -\sum_{q} \big[\tr\big(\hat\Theta_{sq}\hat\eta_{st}^\phi\big) \hat{\eta}^\phi_{sq}
 + \tr\big(\hat\Theta_{lq}\hat\eta_{st}^\phi\big)\hat{\eta}^\phi_{lq}\big] .
\end{gather*}
By \eqref{eq: vanishing traces 3} and \eqref{eq: vanishing traces 4},
\begin{gather*}
n\hat\Theta_{st} =-\tr\big(\hat\Theta_{st}\hat\eta_{st}^\phi\big)\hat\eta_{st}^\phi.
\end{gather*}
Recalling \eqref{eq: eta theta Ric}
\begin{gather*}
n\hat\Theta_{st} = -{1\over \left(4-{n\over 4}-2r\right)}\tr\big({\rm Ric}^M\big)\hat\eta_{st}^\phi=-{{\rm R}\over \left(4-{n\over 4}-2r\right)}\hat\eta_{st}^\phi,
\end{gather*}
i.e.,
\begin{gather*}
\hat\Theta_{st} = {{\rm R}\over n\left({n\over 4}+2r-4\right)}\hat\eta_{st}^\phi.
\end{gather*}
Finally, observe that
\begin{gather*}
{\rm Ric}^M = \left(4-{n\over 4}-2r\right) \hat\eta_{kl}^\phi\hat\Theta_{kl}=
-{{\rm R}\over n\left(4-{n\over 4}-2r\right)} \left(4-{n\over 4}-2r\right) \hat\eta_{kl}^\phi \hat\eta_{kl}^\phi =
{{\rm R}\over n} {\rm Id}_{TM}.\tag*{\qed}
\end{gather*}\renewcommand{\qed}{}
\end{proof}

\begin{Remark} The previous lemma implies the identity
\begin{gather}
\Theta ={{\rm R}\over n\left({n\over 4}+2r-4\right)}\eta^\phi,\label{eq: second equation}
\end{gather}
which shows that the pair formed by the parallel twisted pure spinor $\phi$ and the connec\-tion~$\theta$ satisfies, up to a factor, an equation analogous to the second Seiberg--Witten equation in dimension~4.
\end{Remark}

\subsection{Spinorial characterization of special Riemannian holonomies}
\subsubsection[K\"ahlerian homonomies ${\rm U}(n)$ and ${\rm SU}(n)$]{K\"ahlerian homonomies $\boldsymbol{{\rm U}(n)}$ and $\boldsymbol{{\rm SU}(n)}$}

The K\"ahler and hyper-K\"ahler cases have been treated spinorially by various authors \cite{Hitchin,Kirchberg, Lawson, Moroianu, Wang}. For the sake of completeness, we collect and use some of their ideas to prove the following two corollaries.

\begin{Corollary} An oriented Riemannian manifold $M$ is K\"ahler if and only if it admits a ${\rm Spin}^c$ structure endowed with a connection and carrying a parallel $($classical$)$ pure spinor field.
\end{Corollary}
\begin{proof} Let us assume $M$ is a $2m$-dimensional K\"ahler manifold, $J$ its complex structure, $\ext^{p,q}$~denote the vector bundle of exterior differential forms of type $(p,q)$ and
\begin{gather*}\kappa_M=\ext^{m,0}= \det\big(\ext^{1,0}\big).\end{gather*}
By \cite{Hitchin}, the locally defined Spin bundle decomposes as follows
\begin{gather*}S(TM) = \big(\ext^{0,0}\oplus \cdots \oplus \ext^{0,m}\big)\otimes \kappa_M^{1/2},\end{gather*}
so that the anti-canonical ${\rm Spin}^c$ bundle
\begin{gather*}S(TM)\otimes \kappa_M^{-1/2} =\ext^{0,0}\oplus \cdots \oplus \ext^{0,m}\end{gather*}
contains a trivial summand. Thus, the manifold $M$ admits a parallel spinor field $\psi\in \Gamma\big(\ext^{0,0}\big)$ such that
\begin{gather*}(X+{\rm i}J(X))\cdot\psi=0\end{gather*}
for all $X\in \Gamma(TM)$ (see~\cite{Friedrich}).

Conversely, suppose $M$ admits a ${\rm Spin}^c$ structure carrying a parallel pure spinor field $\psi\in \Gamma \big(S^c(TM)\big)$. If $X\in \Gamma(TM)$, there exists $Y\in \Gamma(TM)$ such that
\begin{gather*}X\cdot \psi= {\rm i} Y\cdot \psi.\end{gather*}
By defining $Y=J(X)$, we see that $J$ is an orthogonal complex structure, and by differentiating
\begin{gather*}\nabla_ZX\cdot \psi = {\rm i} \nabla_Z(J(X))\cdot \psi.\end{gather*}
Note that the vector $\nabla_ZX$ satisfies
\begin{gather*}\nabla_ZX\cdot \psi= {\rm i} J(\nabla_ZX)\cdot \psi,\end{gather*}
so that
\begin{gather*}\big(\nabla_Z(J(X))-J(\nabla_ZX)\big)\cdot \psi=0.\end{gather*}
Since real tangent vectors do not annihilate spinors,
\begin{gather*}\nabla J=0.\tag*{\qed}\end{gather*}\renewcommand{\qed}{}
\end{proof}

\begin{Corollary}Let $M$ be a $2m$-dimensional irreducible oriented Riemannian manifold. The mani\-fold~$M$ is Calabi--Yau if and only if it admits a ${\rm Spin}^c$ structure endowed with a connection carrying two parallel classical pure spinor fields which are complex-linearly independent at one point.
\end{Corollary}
\begin{proof} Let us assume $M$ is Calabi--Yau and $J$ is its complex structure. Since $M$ is Spin and $\kappa_M$ is trivial, we can consider a ${\rm Spin}^c$ structure with trivial auxiliary complex line bundle $L=\kappa_M$ and flat connection. The ${\rm Spin}^c$ spinor bundle
\begin{gather*}S(TM)\otimes \kappa_M^{-1/2} =
 \ext^{0,0}\oplus \cdots \oplus \ext^{0,m}
\end{gather*}
contains two trivial summands generated by parallel spinor fields $\psi_1\in\Gamma\big(\ext^{0,0}\big)$ and $\psi_2\in\Gamma\big(\ext^{0,m}\big)$ such that
\begin{gather*}(X+{\rm i} J(X))\cdot\psi_1=0 \qquad\text{and}\qquad (X-{\rm i} J(X))\cdot\psi_2=0\end{gather*}
for all $X\in \Gamma(TM)$ (see \cite{Friedrich}).

Conversely, suppose $M$ admits a ${\rm Spin}^c$ bundle carrying two parallel classical pure spinor fields $\psi_1$ and $\psi_2$ such that they are complex-linearly independent at one point. We claim that they must be complex-linearly independent everywhere. Suppose there is $z\in\mathbb{C}$ such that $(\psi_1)_y=z(\psi_2)_y$ for some $y\in M$. Since $\psi_1$ and $\psi_2$ are parallel, the spinor field $\psi_1-z\psi_2$ is parallel and its length is constant and equal to zero. Therefore, $\psi_1=z\psi_2$ everywhere.

Thus, the projective classes $[(\psi_1)_x]\not=[(\psi_2)_x]$ for every $x\in M$. If $X\in \Gamma(TM)$, there exist $Y_1, Y_2\in \Gamma(TM)$ such that
\begin{gather*}X\cdot \psi_1= {\rm i} Y_1\cdot \psi_1\qquad\text{and}\qquad X\cdot \psi_2= {\rm i} Y_2\cdot \psi_2.
\end{gather*}
By defining
\begin{gather*}Y_1=J_1(X) \qquad\text{and}\qquad Y_2=J_2(X),
\end{gather*}
we obtain two parallel complex structures. Since orthogonal complex structures are in one to one correspondence with projective classes of classical ${\rm Spin}^c$ pure spinors, $J_1\not=J_2$. If $\Theta$ denotes the curvature 2-form of the connection on the auxiliary ${\rm Spin}^c$ line bundle, by~\cite{Moroianu}
\begin{gather*}{\rm Ric}(X)\cdot \psi_1=i\hat\Theta(X)\cdot \psi_1\qquad\text{and}\qquad {\rm Ric}(X)\cdot \psi_2={\rm i}\hat\Theta(X)\cdot \psi_2 ,\end{gather*}
and
\begin{gather*}J_1\circ {\rm Ric}(X) =\hat\Theta(X)=J_2\circ {\rm Ric}(X),\end{gather*}
i.e., $J_1$ and $J_2$ coincide in the image of the Ricci tensor. The distribution
\begin{gather*}D=\{X\in TM\,|\, J_1(X)=J_2(X)\}\end{gather*}
is parallel and, by irreducibility, it is either equal to $TM$ or trivial. Since $J_1\not = J_2$, $D$ must be trivial and, therefore, ${\rm Ric}=0$ and $\theta$ is flat.
\end{proof}

\subsubsection[Quaternion-K\"ahlerian holonomies ${\rm Sp}(n){\rm Sp}(1)$ and ${\rm Sp}(n)$]{Quaternion-K\"ahlerian holonomies $\boldsymbol{{\rm Sp}(n){\rm Sp}(1)}$ and $\boldsymbol{{\rm Sp}(n)}$}

\begin{Corollary} A Riemannian manifold is quaternion-K\"ahler if and only if it admits a ${\rm Spin}^3$ structure endowed with a~connection and a twisted spinor bundle carrying a parallel twisted pure spinor field.
\end{Corollary}
\begin{proof} Let us assume $M$ is quaternion-K\"ahler so that its orthonormal frame bundle has a~parallel reduction to a principal bundle with fiber ${\rm Sp}(m){\rm Sp}(1)$. We have the following diagram
\begin{gather*}
\xymatrix{
 & {\rm Spin}^3(4m) \ar[d]\\
{\rm Sp}(m){\rm Sp}(1) \ar[ur] \ar[r] & {\rm SO}(4m)\times {\rm SO}(3),
}
\end{gather*}
so that the manifold admits a ${\rm Spin}^3$ structure with an induced connection. We can associate a twisted spinor bundle with fibre $\Delta_{4m}\otimes \Delta_3^{m}$ which contains a trivial ${\rm Sp}(m){\rm Sp}(1)$ summand generated by a pure spinor, such as the spinor given in~\eqref{eq: rank 3 pure spinor} in Section~\ref{subsec:QK-pure-spinor}.

Conversely, if $M$ admits a ${\rm Spin}^3$ structure with a connection and carrying a parallel pure spinor, by Theorem~\ref{theo: parallelness}, we have a parallel quaternion-K\"ahler structure.
\end{proof}

\begin{Corollary} A Riemannian manifold is hyper-K\"ahler if and only if it admits a ${\rm Spin}^3$ structure endowed with a connection and a twisted spinor bundle carrying two parallel twisted pure spinor fields complex-lnearly independent at one point.
\end{Corollary}
\begin{proof} Let us assume $M$ is $4m$-dimensional hyper-K\"ahler, $m\geq 2$. Its structure group reduces further to ${\rm Sp}(m)$. The associated bundle $\Delta_{4m}\otimes \Delta_3^{m}$ contains the $\widehat{{\rm Spin}(3)}$ orbit of the pure spinor in Section~\ref{subsec:QK-pure-spinor}, which consists of pure spinors inducing the same quaternionic structure (see Lemma~\ref{lemma: Spin(r) orbit pure spinor}) and fixed by ${\rm Sp}(m)$.

Conversely, suppose $M$ admits a ${\rm Spin}^3$ structure with a connection and carrying two parallel pure spinors $\psi_1$ and $\psi_2$ complex-linearly independent at one point. Due to the parallelism, they must be complex-linearly independent everywhere. By Theorem~\ref{theo: parallelness}, $M$ is quaternion-K\"ahler and Einstein~\cite{Berger}. By Lemma~\ref{lemma: identity curvatures},
\begin{gather*}{{\rm R}\over 4m (m+2)}\eta^{\psi_1}_{kl}=\Theta_{kl}={{\rm R}\over 4m (m+2)}\eta^{\psi_2}_{kl},\end{gather*}
which also hold when $\dim(M)=8$.
If ${\rm R}\not=0$,
\begin{gather*}\eta^{\psi_1}_{kl}=\eta^{\psi_2}_{kl},\end{gather*}
which means $\psi_1$ and $\psi_2$ have the same annihilator $\mathfrak{sp}(m)\oplus \widetilde{\mathfrak{spin}}(3)\subset\mathfrak{spin}(4m)\oplus \mathfrak{spin}(3)$.
However, restricted to this subalgebra, the representation $\Delta_{4m}\otimes \Delta_3^{\otimes m}$ has only one trivial 1-dimensional summand and
$\psi_2$ must be a multiple of $\psi_1$, which is a contradiction. Hence, $M$ is Ricci-flat.
\end{proof}

\subsubsection[Exceptional holonomies ${\rm Spin}(7)$ and $G_2$]{Exceptional holonomies $\boldsymbol{{\rm Spin}(7)}$ and $\boldsymbol{G_2}$}

\begin{Corollary} A Riemannian $8$-dimensional manifold has holonomy contained in ${\rm Spin}(7)$ if and only if it admits a~${\rm Spin}^7$ structure endowed with a connection and carrying a parallel pure spinor field.
\end{Corollary}
\begin{proof} Let us assume $M$ is an $8$-dimensional Riemannian manifold with holonomy contained in ${\rm Spin}(7)$. Its orthonormal frame bundle has a parallel reduction to a principal bundle with fiber ${\rm Spin}(7)$. We have the following diagram
\begin{gather*}
\xymatrix{
 & {\rm Spin}^7(8) \ar[d]\\
{\rm Spin}(7) \ar[ur] \ar[r] & {\rm SO}(8)\times {\rm SO}(7),
}
\end{gather*}
so that the manifold admits a ${\rm Spin}^7$ structure with an induced connection. We can associate a~twisted spinor bundle with fibre $\Delta_{8}\otimes \Delta_7$ which contains a~trivial ${\rm Spin}(7)$ summand generated by a pure spinor, such as the spinor given in~\eqref{eq: rank 7 pure spinor} in Section~\ref{subsec: Spin(7)}.

Conversely, if $M$ admits a ${\rm Spin}^7$ structure with a connection and carrying a parallel pure spinor, by Theorem~\ref{theo: parallelness}, it admits a parallel rank $7$ even-Clifford structure.
\end{proof}

Unlike the complex and quaternionic cases, the $G_2$ holonomy reduction does not arise by the existence of two linearly independent parallel twisted pure spinors belonging to the same ${\rm Spin}^7(8)$-orbit. As it could be expected, the holonomy reduction to this exceptional Lie group is due to triality, which in our context is
expressed by the interaction of twisted pure spinors whose 2-forms are related by triality.

\begin{Corollary} A simply connected $8$-dimensional Riemannian manifold $M$ is a product of a~flat $1$-dimensional manifold $\mathcal{S}$ and a $7$-dimensional manifold~$N$ with holonomy contained in~$G_2$ if and only if~$M$ admits a~${\rm Spin}^7$ structure with a connection $($which includes the lift of the Levi-Civita connection$)$, carrying a parallel twisted pure spinor~$\phi_1$ and a parallel spinor~$\phi_2$ such that the vector field
\begin{gather*}X^{\phi_1,\phi_2}:=\sum_{i=1}^8 \operatorname{Re}\langle e_i\cdot \phi_1,\phi_2\rangle e_i\end{gather*}
is nonzero at some point of $M$.
\end{Corollary}
\begin{proof} Since $G_2\subset {\rm SO}(7)$, we can embed it as a block ${\rm SO}(7)\subset {\rm SO}(8)$, in such a~way the~$G_2$ leaves invariant the first canonical vector $e_1$ of $\mathbb{R}^8$. Now, for a 7-dimensional Riemannian manifold $N$ with holonomy contained in $G_2$, this corresponds to considering the product $\mathcal{S}\times N $ where $\mathcal{S}$ is a flat 1-dimensional manifold. Now that $G_2\subset {\rm SO}(8)$, there are many copies of ${\rm Spin}(7)\subset {\rm SO}(8)$ that contain it. Choose one of such copies. We have the diagram
\begin{gather*}
\xymatrix{
 & {\rm Spin}^7(8) \ar[d]\\
G_2\subset {\rm Spin}(7) \ar[ur] \ar[r] & {\rm SO}(8)\times {\rm SO}(7).
}
\end{gather*}
Now, $\Delta_8\otimes\Delta_7$ decomposes under the image of such a ${\rm Spin}(7)$
\begin{gather*}
\Delta_8\otimes\Delta_7 =\Delta_8^+\otimes\Delta_7\oplus\Delta_8^-\otimes\Delta_7\cong\Delta_7\otimes\Delta_7\oplus\Delta_7\otimes\Delta_7
\end{gather*}
and one of the two summands contains an invariant element which is a pure spinor, say $\phi_1\in\Delta_8^+\otimes\Delta_7$. All this translates into a parallel spinor field for $N\times\mathcal{S}$. Since we have a globally defined vector field $e_1$, we can consider the parallel spinor field $\phi_2=e_1\cdot \phi_1\in\Delta_8^-\otimes \Delta_7$. By Lemma~\ref{lemma: producing more pure spinors}, the spinor $\phi_2$ is also pure and its stabilizer is also a~copy of ${\rm Spin}(7)$, but is not in the same orbit of ${\rm Spin}^7(8)$. The corresponding copies of ${\rm Spin}(7)$ in ${\rm SO}(8)$ intersect in the original $G_2$. Note that
\begin{gather*}
X^{\phi_1,\phi_2} =\sum_{i=1}^8 \operatorname{Re}\langle e_i\cdot \phi_1,\phi_2\rangle e_i
= \operatorname{Re}\langle e_1\cdot \phi_1,e_1\cdot\phi_1\rangle e_1 + \sum_{i=2}^8 \operatorname{Re}\langle e_i\cdot \phi_1,e_1\cdot\phi_1\rangle e_i\\
\hphantom{X^{\phi_1,\phi_2}}{} = \operatorname{Re} |\phi_1|^2e_1\not = 0.
\end{gather*}

Conversely, suppose we have a simply connected 8-dimensional Riemannian manifold $M$ with a ${\rm Spin}^7(8)$ structure and a connection that carries a parallel twisted pure spinor $\phi_1$ and another parallel spinor $\phi_2$ such that the vector field
\begin{gather*}X^{\phi_1,\phi_2}=\sum_{i=1}^8 \operatorname{Re}\langle e_i\cdot \phi_1,\phi_2\rangle e_i\end{gather*}
is such that $X_p\not=0$ at some $p\in M$. Note that for any $Z\in\Gamma(TM)$,
\begin{gather*}
\nabla_ZX^{\phi_1,\phi_2} = \nabla_Z\sum_{i=1}^8 \operatorname{Re}\langle e_i\cdot \phi_1,\phi_2\rangle e_i = \sum_{i=1}^8 \nabla_Z\big(\operatorname{Re}\langle e_i\cdot \phi_1,\phi_2\rangle e_i\big)\\
\hphantom{\nabla_ZX^{\phi_1,\phi_2}}{}
 = \sum_{i=1}^8 \operatorname{Re}\big\langle \nabla_Ze_i\cdot \phi_1,\phi_2\big\rangle e_i
 + \sum_{i=1}^8 \operatorname{Re}\big\langle e_i\cdot \nabla_Z^S\phi_1,\phi_2\big\rangle e_i\\
\hphantom{\nabla_ZX^{\phi_1,\phi_2}=}{} + \sum_{i=1}^8 \operatorname{Re}\big\langle e_i\cdot \phi_1,\nabla_Z^S\phi_2\big\rangle e_i
 + \sum_{i=1}^8 \operatorname{Re}\langle e_i\cdot \phi_1,\phi_2\rangle \nabla_Ze_i\\
 \hphantom{\nabla_ZX^{\phi_1,\phi_2}}{}
 = \sum_{i=1}^8 \operatorname{Re}\langle \nabla_Ze_i\cdot \phi_1,\phi_2\rangle e_i
 + \sum_{i=1}^8 \operatorname{Re}\langle e_i\cdot \phi_1,\phi_2\rangle \nabla_Ze_i\\
\hphantom{\nabla_ZX^{\phi_1,\phi_2}}{}
 = \sum_{i=1}^8 \operatorname{Re}\bigg\langle \sum_{j=1}^n\langle \nabla_Ze_i,e_j\rangle e_j\cdot \phi_1,\phi_2\bigg\rangle e_i
 + \sum_{i=1}^8 \operatorname{Re}\langle e_i\cdot \phi_1,\phi_2\rangle \sum_{j=1}^8\langle \nabla_Ze_i,e_j\rangle e_j\\
\hphantom{\nabla_ZX^{\phi_1,\phi_2}}{}
 = \sum_{i=1}^8 \operatorname{Re}\bigg\langle \sum_{j=1}^8\omega_{ji}(Z)e_j\cdot \phi_1,\phi_2\bigg\rangle e_i
 + \sum_{i=1}^8 \operatorname{Re}\langle e_i\cdot \phi_1,\phi_2\rangle \sum_{j=1}^8\omega_{ji}(Z)e_j\\
\hphantom{\nabla_ZX^{\phi_1,\phi_2}}{}
 = \sum_{i=1}^8 \sum_{j=1}^8\omega_{ji}(Z)\operatorname{Re}\langle e_j\cdot \phi_1,\phi_2\rangle e_i
 + \sum_{i=1}^8 \sum_{j=1}^8\omega_{ji}(Z)\operatorname{Re}\langle e_i\cdot \phi_1,\phi_2\rangle e_j\\
 \hphantom{\nabla_ZX^{\phi_1,\phi_2}}{}
 = \sum_{i=1}^8 \sum_{j=1}^8\omega_{ji}(Z)\operatorname{Re}\langle e_j\cdot \phi_1,\phi_2\rangle e_i
 + \sum_{j=1}^8 \sum_{i=1}^8\omega_{ij}(Z)\operatorname{Re}\langle e_j\cdot \phi_1,\phi_2\rangle e_i
 = 0,
\end{gather*}
i.e., $X^{\phi_1,\phi_2}$ is parallel. As a consequence, it has constant length which is non-zero, and therefore~$M$ decomposes as a product of a flat 1-dimensional manifold $\mathcal{S}$ and a $7$-dimensional Riemannian manifold~$N$. Since $M$ carries a parallel twisted pure spinor, it has holonomy contained in ${\rm Spin}(7)\subset {\rm SO}(8)$. But now that we have proved that $M=\mathcal{S}\times N$, the holonomy of $M$ must reduce to the subgroup
\begin{gather*}{\rm Spin}(7)\cap (\{1\}\times {\rm SO}(7)) =G_2, \end{gather*}
i.e., the $7$-dimensional Riemannian manifold $N$ has holonomy contained in $G_2$.
\end{proof}

\section{Clifford monopole equations}\label{sec:Clifford-monopoles}

Let $M$ be a ${\rm Spin}^r$ manifold with auxiliary bundle $P_{{\rm SO}(r)}$ endowed with a connection $\theta$, $F$ the associated Riemannian rank $r$ vector bundle and $m\in\mathbb{N}$ be such that the twisted Dirac operator
\begin{align*}
\dirac^\theta\colon \ \Gamma(S(M,F,m))&\lra \Gamma(S(M,F,m)),\\
\dirac^\theta\phi &= \sum_{i=1}^{\dim(M)} e_i\cdot \nabla_{e_i}^\theta \phi
\end{align*}
is well-defined, where the vectors $e_i$ form a local orthonormal frame of the tangent bundle. The {\em Clifford monopole equations} are
\begin{gather}
 \dirac^\theta\phi = 0,\qquad \Theta = E\big(\eta^\phi\big),\label{eq:clifford-monopole-equations}
\end{gather}
where
\begin{gather*}
 \Theta = \sum_{1\leq k<l \leq r} \Theta_{kl}\otimes f_{kl} \in \Gamma\big(\ext^2T^*M\otimes \ext^2F\big)
\end{gather*}
is the curvature of $\theta$,
\begin{gather*}
 \eta^\phi = \sum_{1\leq k<l \leq r} \eta_{kl}^\phi\otimes f_{kl} \in \Gamma\big(\ext^2T^*M\otimes \ext^2F\big)
\end{gather*}
is the 2-form with values in $\ext^2 F$ associated to $\phi$, and $E$ is a suitable endomorphism of 2-forms. A pair $(\phi,\theta)$ satisfying~\eqref{eq:clifford-monopole-equations} will be called a {\em Clifford monopole}.

Here, we will show that the Clifford monopole equations restrict to the Seiberg--Witten equations on 4-manifolds, and will also exhibit Clifford monopoles on manifolds with special Riemannian holonomy.

\looseness=-1 As mentioned in the introduction, preliminary work with A.~Quintero indicates the existence of a smooth compact moduli space which, according to the Mathai--Quillen--Atiyah--Jeffrey formalism, will give raise to a topological quantum field theory. Such topological field theory, at least in dimensions $8$, might turn out to be a topological twist of an $N=2$ supersymmetric theory.

\subsection{The Clifford monopole equations on 4-manifolds}

In this subsection, we will show that the Clifford monopole equations restrict to the Seiberg--Witten monopole equations for appropriate choices of the parameters.

Let us start by recalling that every 4-manifold admits a ${\rm Spin}^c$ structure~\cite{Morgan}, i.e., a ${\rm Spin}^2$ structure in our notation. Here, we choose $m=1$ so that $\phi\in\Gamma(S(\Delta_4\otimes \Delta_2))$. The orthogonal decomposition
\begin{gather*}
\Delta_4\otimes \Delta_2 = \big(\Delta_4^+\oplus \Delta_4^-\big)\otimes \big(\Delta_2^+\oplus \Delta_2^-\big) =
 \Delta_4^+\otimes \Delta_2^+ \oplus \Delta_4^+\otimes \Delta_2^- \oplus \Delta_4^-\otimes \Delta_2^+ \oplus \Delta_4^-\otimes \Delta_2^- .
\end{gather*}
implies $\phi=\phi_1+\phi_2+\phi_3+\phi_4$ where
\begin{gather*}
\phi_1 \in \Delta_4^+\otimes\Delta_2^+,\qquad
\phi_2 \in \Delta_4^+\otimes\Delta_2^-,\qquad
\phi_3 \in \Delta_4^-\otimes\Delta_2^+,\qquad
\phi_4 \in \Delta_4^-\otimes\Delta_2^-.
\end{gather*}
Thus, we have four Dirac equations
\begin{gather*}
\dirac^\theta\phi_s =0,\qquad s=1,2,3,4.
\end{gather*}

{\allowdisplaybreaks
Note that, in this case, $f_1f_2$ represents a globally defined section of $\ext^2 F$ and acts as multiplication by $\pm {\rm i}$ on $\Delta_4\otimes\Delta_2^\pm$ respectively. Now consider
\begin{gather*}
\eta^\phi = \sum_{1\leq a<b\leq 4}\langle e_ae_b\cdot f_{12}\cdot (\phi_1+\phi_2+\phi_3+\phi_4),
 (\phi_1+\phi_2+\phi_3+\phi_4)\rangle e_ae_b\otimes f_{12}\\
\hphantom{\eta^\phi}{}= \sum_{s,t=1}^4 \sum_{1\leq a<b\leq 4}\langle e_ae_b\cdot f_{12}\cdot \phi_s,\phi_t\rangle e_ae_b\otimes f_{12}\\
\hphantom{\eta^\phi}{} = \sum_{s=1}^4 \sum_{1\leq a<b\leq 4}\langle e_ae_b\cdot f_{12}\cdot \phi_s,\phi_s\rangle e_ae_b\otimes f_{12}\\
\hphantom{\eta^\phi}{}= \sum_{s=1}^4(-1)^{s-1}{\rm i}\sum_{1\leq a<b\leq 4}\langle e_ae_b\cdot (\phi_s),\phi_s\rangle e_ae_b\otimes f_{12}\\
\hphantom{\eta^\phi}{} = \big({\rm i}\tau^{\phi_1} - {\rm i}\tau^{\phi_2} + {\rm i}\tau^{\phi_3} - {\rm i}\tau^{\phi_4}\big) \otimes f_{12} ,
\end{gather*}
where
\begin{gather*}\tau^{\psi}(X,Y):=\langle X\wedge Y\cdot\psi,\psi\rangle \end{gather*}
for $\psi\in\Delta_4\otimes\Delta_2$. The 2-forms $\tau^{\phi_s}$ are purely imaginary, and
\begin{gather*}{\rm i}\tau^{\phi_1},{\rm i}\tau^{\phi_2}\in \ext^2_+T^*M,\qquad {\rm i}\tau^{\phi_3},{\rm i}\tau^{\phi_4}\in \ext^2_-T^*M.\end{gather*}
Indeed, $\phi_1=(\alpha u_+\otimes u_+ + \beta u_-\otimes u_- )\otimes v_+ $ for some $\alpha,\beta\in\mathbb{C}$, so that
\begin{gather*}{\rm i}\tau^{\phi_1}=
\big({-}|\alpha|^2+|\beta|^2\big)(e_1e_2+e_3e_4) -2\operatorname{Im}\big(\alpha\overline{\beta}\big)(e_1e_3-e_2e_4)\\
\hphantom{{\rm i}\tau^{\phi_1}=}{}
+2\operatorname{Re}\big(\alpha\overline{\beta}\big)(e_1e_4+e_2e_3)\in\ext^2_+T^*M.
\end{gather*}
Similarly for the other 2-forms.

}

At this point we choose $E$ to be such that $\ext^2_+T^*M$ and $\ext^2_-T^*M$ are invariant subspaces. Since $f_{12}$ trivializes $\ext^2 F$, the second Clifford monopole equation
\begin{gather*}\Theta_{12}\otimes f_{12} = E\big(\eta^{\phi}\big)\otimes f_{12}\end{gather*}
becomes
\begin{gather*}\Theta_{12} = E\big(\eta_{12}^\phi\big),\end{gather*}
which splits as follows
\begin{gather*}
{\rm i}\Theta_{12}^+ = -E\big(\tau^{\phi_1}-\tau^{\phi_2}\big),\qquad
{\rm i}\Theta_{12}^- = -E\big(\tau^{\phi_3}-\tau^{\phi_4}\big).
\end{gather*}

If we limit ourselves to work with spinor fields in $\Gamma\big(S\big(\Delta_4\otimes\Delta_2^+\big)\big)$, i.e., $\phi_2=0$ and $\phi_4=0$, we have
\begin{gather*}
{\rm i}\Theta_{12}^+ = -E\big(\tau^{\phi_1}\big),\qquad
{\rm i}\Theta_{12}^- = -E\big(\tau^{\phi_3}\big),
\end{gather*}
and setting $E={1\over 4}{\rm Id}_{\bigwedge^2T^*M}$, we have
\begin{gather*}
{\rm i}\Theta_{12}^+ = -{1\over 4}\tau^{\phi_1},\qquad
{\rm i}\Theta_{12}^- = -{1\over 4}\tau^{\phi_3}.
\end{gather*}
In 4 dimensions \cite{Morgan}, however, if there is a non-zero solution $\phi_1$ for
\begin{gather*}\dirac^\theta\phi_1=0,\end{gather*}
then if $\phi_3$ is such that
\begin{gather*}\dirac^\theta\phi_3=0,\end{gather*}
it must vanish identically, i.e., $\phi_3\equiv0$.

Thus, we are left with the Seiberg--Witten equations for pairs $(\phi_1,\theta)$~\cite{Friedrich}
\begin{gather*}
\dirac^\theta \phi_1 = 0,\qquad
{\rm i}\Theta = -{1\over 4}\tau^{\phi_1}.
\end{gather*}

If, on the other hand, we consider only spinor fields in $\Gamma(S(\Delta_4\otimes\Delta_2^-))$, i.e., $\phi_1=0$ and $\phi_3=0$, and set $E=-{1\over 4}{\rm Id}_{\bigwedge^2T^*M}$, we end up with the Seiberg--Witten equations $(\phi_2,\theta)$
\begin{gather*}
\dirac^\theta \phi_2 = 0,\qquad
{\rm i}\Theta = -{1\over 4}\tau^{\phi_2}.
\end{gather*}

\subsection[8-manifolds with ${\rm Spin}(7)$ holonomy]{8-manifolds with $\boldsymbol{{\rm Spin}(7)}$ holonomy}

Let $M$ be an $8$-dimensional manifold with holonomy contained in ${\rm Spin}(7)$ and $\theta$ denote the Levi-Civita connection restricted to the holonomy principal bundle. We have seen that there exists a parallel pure spinor field $\phi_1$ which characterizes such holonomy reduction (see~\eqref{eq: rank 7 pure spinor}), where $m=1$. By~\eqref{eq: first version second SW equation}
\begin{gather*}
\Theta = -{1\over 4}R^M\bullet\eta^{\phi_1}.
\end{gather*}
Thus $(\phi_1,\theta)$ gives a solution to~\eqref{eq:clifford-monopole-equations} with $E=-{1\over 4}R^M\bullet$.

\begin{Remark}Note that spinor can be scaled in order to remove positive constants so that $(\phi_1/2,\theta)$ gives a solution to~\eqref{eq:clifford-monopole-equations} with $E=-R^M\bullet$.
\end{Remark}

\subsection{Quaternion-K\"ahler manifolds}
Let $M$ be a $4m$-dimensional quaternion-K\"ahler manifold and denote by $\theta$ the connection induced by the Levi-Civita connection on the corresponding ${\rm SO}(3)$-bundle.
We have seen that there exists a parallel pure spinor field $\phi$ which characterizes the holonomy (see \eqref{eq: rank 3 pure spinor}).
By \eqref{eq: first version second SW equation} and \eqref{eq: second equation},
\begin{gather*}
\hat\Theta = -{1\over m}R^M\bullet \eta^{\phi}= {{\rm R}\over 4m\left(m+2\right)} \eta^{\phi} .
\end{gather*}
Thus, we can say that either
\begin{itemize}\itemsep=0pt
\item the pair $\big( \phi/\sqrt{m},\theta\big)$ satisfies~\eqref{eq:clifford-monopole-equations} with $E=-R^M\bullet$;
\item or if ${\rm R}$ is non-negative, the pair $\Big( \sqrt{{\rm R}\over 4m(m+2)}\phi,\theta\Big)$ satisfies~\eqref{eq:clifford-monopole-equations} with $E={\rm Id}_{\ext^2 T^*M}$;
\item or if ${\rm R}$ is non-positive, the pair $\Big( \sqrt{-{\rm R}\over 4m(m+2)}\phi,\theta\Big)$ satisfies~\eqref{eq:clifford-monopole-equations} with $E=-{\rm Id}_{\ext^2 T^*M}$.
\end{itemize}

\begin{Remark}The choices $E=\pm{\rm Id}_{\ext^2 T^*M}$ are the ones that remind us of the Seiberg--Witten equations in dimension 4.
\end{Remark}

\subsection{K\"ahler manifolds}
Let $(M,\langle \cdot,\cdot\rangle , J)$ be a $2n$-dimensional K\"ahler manifold with scalar curvature ${\rm R}$ and canonical ${\rm Spin}^c$ structure with canonical connection $\theta$. Up to a choice of local frames, the complex structure is in correspondence with the parallel the twisted pure spinor
\begin{gather*}\phi=\underbrace{u_{-1}\otimes\cdots \otimes u_{-1}}_{n\,\, {\rm times}}\otimes v_{1}\in \Delta_{2n}^-\otimes \Delta_2^+\end{gather*}
and
\begin{gather*}
 \hat\eta_{12}^{\phi} = \left(\begin{matrix}
 0 & -1 & & & \\
 1 & 0 & & & \\
 & & \ddots & & \\
 & & & 0 & -1 \\
 & & & 1 & 0
 \end{matrix} \right).
\end{gather*}
Note that, in this rank, the sections $f_{12}$, $\eta_{12}^\phi$, $\hat\eta_{12}^\phi$ are actually globally defined.

Since
\begin{gather*}\nabla^\theta \phi =0,\end{gather*}
we have
\begin{gather*}
 \nabla \eta_{12}^\phi= 0
\qquad \text{and}\qquad
\big[{R^M(X,Y)},\hat\eta_{12}^\phi\big] =0,\end{gather*}
where $X,Y\in\Gamma(TM)$. Multiplying by $\hat\eta_{12}^\phi$ on the left
\begin{gather*}
\hat{\eta}^\phi_{12}{R^M(X,Y)}\hat{\eta}^\phi_{12}+{R^M(X,Y)} =0.
\end{gather*}
Setting $Y=e_j$ and summing over $i$, $j$,
\begin{gather*}
0=\sum_{i,j}\big\langle \hat{\eta}^\phi_{12}{R^M(X,e_i)}\hat{\eta}^\phi_{12}(e_i),e_j\big\rangle e_j+\sum_{i,j}\big\langle {R^M(X,e_i)}(e_i),e_j\big\rangle e_j\\
\hphantom{0}{} = -\sum_{i,j}\big\langle {R^M(X,e_i)}\hat{\eta}^\phi_{12}(e_i),\hat{\eta}^\phi_{12}e_j\big\rangle e_j+\sum_{i,j}\big\langle {R^M(X,e_i)}(e_i),e_j\big\rangle e_j\\
\hphantom{0}{} = \hat\eta_{12}^\phi\sum_{i,j}\big\langle {R^M(X,e_i)}\hat{\eta}^\phi_{12}(e_i),\hat{\eta}^\phi_{12}e_j\big\rangle \hat\eta_{12}^\phi (e_j)+{\rm Ric}(X)\\
\hphantom{0}{} = \hat\eta_{12}^\phi\sum_{i,j}\big\langle {R(X,e_i)}\hat{\eta}^\phi_{12}(e_i),e_j'\big\rangle e_j' +{\rm Ric}(X)\\
\hphantom{0}{} = \hat\eta_{12}^\phi\big({-}{R^M\bullet\eta_{12}^\phi}(X)\big)+{\rm Ric}(X) ,
\end{gather*}
i.e.,
\begin{gather*}
{\rm Ric} =\hat\eta_{12}^\phi{R^M\bullet\eta_{12}^\phi},
\end{gather*}
so that
\begin{gather*}
{R^M\bullet\eta_{12}^\phi} =-{\rm Ric} \circ\eta_{12}^\phi.
\end{gather*}

By \eqref{eq: basic curvature identity},
\begin{gather*}
0 = R^M(X,Y)\cdot \phi - \Theta_{12}(X,Y)\kappa_{2}^*(f_{12})\cdot \phi.
\end{gather*}
Proceeding as in \cite{Friedrich}, set $Y=e_a$, multiply by $e_a$ and sum over $a$
\begin{gather*}
0 = \sum_a e_a\cdot R^M(X,e_a)\cdot \phi - \sum_a \Theta_{12}(X,e_a)e_a\cdot\kappa_{2*}^1(f_{12})\cdot \phi\\
\hphantom{0}{}= -{\rm Ric}(X)\cdot \phi -\hat\Theta_{12}(X)\cdot \kappa_{2*}^1(f_{12})\cdot \phi,
\end{gather*}
i.e.,
\begin{gather*}
{\rm Ric}(X)\cdot \phi = -\hat\Theta_{12}(X)\cdot \kappa_{2*}^1(f_{12})\cdot \phi.
\end{gather*}
Now consider, for $i<j$,
\begin{gather*}
{\rm Ric}_{ji} =\operatorname{Re}\big\langle {\rm Ric}(e_i)\cdot \phi , e_j \cdot \phi\big\rangle =
 -\sum_a \Theta_{12}(e_i,e_a) \operatorname{Re}\big\langle e_a \cdot\kappa_{2*}^1(f_{12})\cdot \phi, e_j\cdot \phi\big\rangle \\
\hphantom{{\rm Ric}_{ji}}{} =
 \sum_a \Theta_{12}(e_i,e_a) \operatorname{Re}\big\langle e_je_a \cdot\kappa_{2*}^1(f_{12})\cdot \phi, \phi\big\rangle
 = \sum_a \big(\hat\Theta_{12}\big)_{ai} \eta_{12}^\phi(e_j,e_a) \\
\hphantom{{\rm Ric}_{ji}}{}
 = \sum_a \big(\hat\Theta_{12}\big)_{ai} \big(\hat\eta_{12}^\phi\big)_{aj} =- \sum_a \big(\hat\eta_{12}^\phi\big)_{ja}\big(\hat\Theta_{12}\big)_{ai} = - \big(\hat\eta_{12}^\phi\hat\Theta_{12}\big)_{ji} ,
\end{gather*}
i.e.,
\begin{gather*}
{\rm Ric} = -\hat\eta_{12}^\phi\hat\Theta_{12}.
\end{gather*}
Since ${\rm Ric}$ is symmetric, $\hat\Theta_{12}$ and $\hat\eta_{12}^\psi$ commute and also commute with ${\rm Ric}$.
Thus,
\begin{gather*}\hat\Theta_{12}= {\rm Ric}\circ\hat\eta_{12}^\phi,\qquad \text{i.e.},
\qquad \Theta= -R^M\bullet\eta^\phi,\end{gather*}
so that the pair $(\phi,\theta)$ satisfies \eqref{eq:clifford-monopole-equations} with $E=-R^M\bullet$.

Furthermore, if we assume $M$ is Einstein
\begin{gather*}
\Theta ={{\rm R}\over 2n}\eta^\phi.
\end{gather*}
If we assume the scalar curvature ${\rm R}$ to be constant, we can say that
\begin{itemize}\itemsep=0pt
\item if ${\rm R}$ is non-negative, the pair $\Big(\sqrt{{R\over 2n}}\phi, \theta\Big)$ satisfies \eqref{eq:clifford-monopole-equations} with $E={\rm Id}_{\ext^2T^*M}$;
\item if ${\rm R}$ is non-positive, the pair $\Big(\sqrt{{-R\over 2n}}\phi, \theta\Big)$satisfies \eqref{eq:clifford-monopole-equations} with $E=-{\rm Id}_{\ext^2T^*M}$.
\end{itemize}

\begin{Remark} The last two statements remind us again of the Seiberg--Witten equations in dimension~4. We could even consider the positive/negative spinor decomposition and the self-dual/anti-self-dual 2-form decomposition to arrive at two sets of equations: the Seiberg--Witten equations and an analogous pair of equations for negative spinors and anti-self-dual forms.
\end{Remark}

\subsection*{Acknowledgements}

The first named author would like to thank H.~Baum and U.~Bruzzo for their hospitality, as well as the International Centre for Theoretical Physics, the Institut des Hautes \'Etudes Scientifiques, and the Scuola Internazionale Superiore di Studi Avanzati for their hospitality and support. The authors would also like to thank Alexander Quintero for his very valuable insights, as well as the anonymous referees for their comments which helped
improve the paper.

The first author was partially supported by grants from CONACyT, LAISLA (CONACyT-CNRS), INFN-Italy and IMU Berlin Einstein Foundation, and the second author was partially supported by grants from CONACyT and LAISLA (CONACyT-CNRS).

\pdfbookmark[1]{References}{ref}
\LastPageEnding

\end{document}